\documentclass[a4paper]{article}
\usepackage[top=30mm,bottom=30mm,left=25mm,right=25mm]{geometry}
\usepackage{amsmath,amssymb,amsthm}
\usepackage{mathtools}
\usepackage{placeins}
\usepackage{graphicx}
\usepackage{tikz}
\usepackage{xcolor}
\usepackage{authblk}
\usepackage{hyperref}
\usepackage[style=alphabetic,backend=biber]{biblatex}
\mathtoolsset{showonlyrefs=true}
\theoremstyle{plain}
\newtheorem{theorem}{Theorem}[section]
\newtheorem{proposition}[theorem]{Proposition}
\newtheorem{corollary}[theorem]{Corollary}
\newtheorem{lemma}[theorem]{Lemma}
\newtheorem{assumption}[theorem]{Assumption}
\theoremstyle{definition}
\newtheorem{definition}[theorem]{Definition}
\theoremstyle{remark}
\newtheorem{remark}[theorem]{Remark}
\makeatletter
\@addtoreset{equation}{section}

\makeatother
\newcommand{\emails}[1]{\begin{NoHyper}\let\thefootnote\relax\footnotetext{#1}\end{NoHyper}}

\begin{document}
\title{Stationary bistable fronts on a periodic necklace graph:\\ existence, stability and multiplicity%
	\thanks{YM and CS were partially supported by the NSF Materials Research Science and Engineering Center (DMR-2309034).
		TO was supported by JSPS KAKENHI Grant Numbers JP24K06854.
	}}
\author[1]{Yoichiro Mori}
\author[2]{Toshiyuki Ogawa}
\author[1,3]{Chanoknun Sintavanuruk}
\affil[1]{Department of Mathematics, University of Pennsylvania, Philadelphia, PA 19104, USA}
\affil[2]{Graduate School of Advanced Mathematical Sciences, Meiji University, Tokyo, Japan}
\affil[3]{Department of Physiology, Faculty of Medicine Siriraj Hospital, Mahidol University, Bangkok, Thailand}

\date{}
\maketitle
\emails{\textit{Email addresses:} \texttt{y1mori@sas.upenn.edu} (Y.~Mori),
	\texttt{togw@meiji.ac.jp} (T.~Ogawa),
	\texttt{aidensin@sas.upenn.edu} (C.~Sintavanuruk).
	\textit{Corresponding author:} C.~Sintavanuruk.}

\begin{abstract}
	We study stationary fronts for a bistable reaction--diffusion equation
	with the McKean reaction on a periodic necklace metric graph.
	For profiles that agree on each pair of parallel edges,
	the problem reduces to a scalar equation on the real line
	with a periodic, piecewise constant diffusion coefficient.
	The diffusion contrast of the necklace is not small,
	so the coefficient cannot be treated as a small perturbation of a constant.
	Because the McKean reaction is piecewise linear,
	the equation is affine between threshold crossings and coefficient interfaces,
	and a front is determined explicitly by its crossing points.
	We find an explicit interval of reaction thresholds, the pinning interval,
	such that stationary fronts, monotone or not, exist exactly for thresholds in this interval.
	Up to translations by the spatial period,
	there are exactly two decreasing fronts at each interior threshold
	and one front at each endpoint.
	For each fixed interior threshold,
	one decreasing front is Lyapunov stable and the other is unstable
	in the $W^{1,\infty}(\mathbb{R})$ norm,
	among mild solutions with one continuously moving threshold crossing.
	The crossing phases of the two decreasing fronts also determine nonmonotone fronts.
	For each interior threshold and odd $m\geq3$,
	choose $m$ periods of the coefficient that are sufficiently far apart.
	There are then exactly $2^m$ fronts with one threshold crossing in each chosen period and none elsewhere.
	Varying the chosen periods gives infinitely many inequivalent fronts with $m$ crossings.
\end{abstract}



\section{Introduction}
\label{sec:introduction}

Spatial heterogeneity can pin a bistable reaction--diffusion front,
so that the transition between two stable states remains stationary.
Haudin et al.\ \cite{haudinDrivenFrontPropagation2009}
observed such fronts over an interval of applied voltages
in a nematic liquid-crystal cell with periodic optical forcing.
Geometry alone can also arrest a front.
In experiments with \textit{Xenopus} egg extracts,
apoptotic activation waves stopped where a narrow channel opened into a chamber,
although the regions remained diffusively connected
\cite{liGeometricalCompartmentalizationTrigger2025}.
For spatially uniform reactions and diffusion,
narrow connections can support stable nonconstant equilibria on bounded domains
\cite[Section~6]{matanoAsymptoticBehaviorStability1979},
and narrow passages or periodic obstacles can prevent invasion in unbounded domains
\cite{berestyckiFrontBlockingPropagation2016,ducasseBlockingInvasionReactiondiffusion2018}.
Whether a geometry blocks a front depends on its length scales relative to the front width.
For corrugated channels, perturbative calculations and simulations relate propagation failure
to the ratio of the channel period to the front width
\cite{martensFrontPropagationThreedimensional2015}.

Branching introduces a similar question for networks of channels or tissue strands.
Experiments on cardiac-cell strands with repeated side branches found slower impulse propagation
when the branches were longer or more closely spaced
\cite{kuceraSlowConductionCardiac1998}.
Metric graphs represent such networks by edges of prescribed lengths
joined at vertices, with continuity and Kirchhoff flux balance at each junction.
For bistable reactions, propagation and blocking have been studied on graphs
with one or two junctions \cite{JM1,JM2,IJM}
and on symmetric tree-like graphs with successive junctions
\cite{jimboFrontPropagationBistable2024}.
In the latter case, increasing the distance between neighboring junctions
can change blocking into propagation.
For a bounded finite graph core with attached half-lines,
Matano and Jimbo \cite{matanoFrontPropagationGeneral2025a}
proved existence and uniqueness of an entire solution arriving from one half-line
and convergence to a stationary limit.
The value of this limit at each attachment vertex determines propagation or blocking
on the corresponding outgoing half-line.
Monobe and Morita \cite{monobeSpatialPatternsStable2025}
constructed nonconstant stationary patterns on graphs assembled from identical stars with sufficiently long edges.
These results do not determine, for a periodically branching graph,
which geometries pin a bistable front.

\begin{figure}[htbp]
    \centering
\begin{tikzpicture}[x=1cm,y=1cm,line width=0.8pt,
    every node/.style={font=\small}]
    \draw (0,0)--(0.65,0)
          (3.15,0)--(5.25,0)
          (7.75,0)--(9.85,0)
          (12.35,0)--(13,0);
    \foreach \left/\right in {0.65/3.15,5.25/7.75,9.85/12.35} {
        \draw (\left,0)
            .. controls (\left+0.5,0.8) and (\right-0.5,0.8) .. (\right,0);
        \draw (\left,0)
            .. controls (\left+0.5,-0.8) and (\right-0.5,-0.8) .. (\right,0);
        \fill (\left,0) circle (1.7pt);
        \fill (\right,0) circle (1.7pt);
    }
    \node at (-0.35,0) {$\cdots$};
    \node at (13.35,0) {$\cdots$};
    \node[above] at (4.2,0.08) {$L_1$};
    \node[above] at (6.5,0.64) {$L_2/2$};
    \node[below] at (6.5,-0.64) {$L_2/2$};
\end{tikzpicture}
    \caption{A necklace metric graph.
    Single edges have length $L_1$, and each edge in a parallel pair has length $L_2/2$.
    Doubling the coordinate on each parallel edge gives the scalar period
    $L=L_1+L_2$ and diffusion coefficient $p=4$ on the corresponding intervals.}
    \label{fig:necklace_graph}
\end{figure}
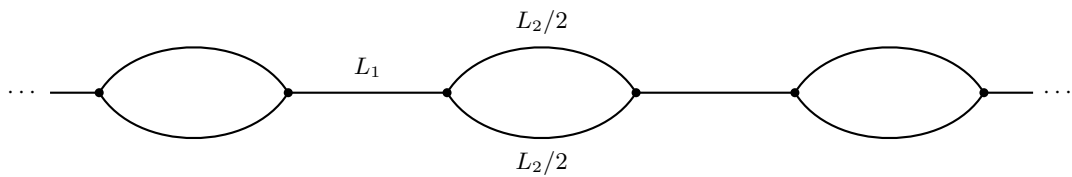

We consider a necklace metric graph,
a periodic succession of single edges and pairs of parallel edges
(Figure~\ref{fig:necklace_graph}).
Branching and reconnection introduce heterogeneity through the Kirchhoff conditions,
even though the equation is identical on every edge.
The same geometry supports localized solutions in other equations.
Maier, Reichel, and Schneider \cite{maierBreatherSolutionsSemilinear2023}
used spectral gaps of the necklace Laplacian to construct breathers
for a semilinear Klein--Gordon equation.
For the focusing cubic nonlinear Schr\"odinger equation on the necklace,
Le Coz, Pelinovsky, and Schneider \cite{leCozTravelingWavesPeriodic2026}
used spatial dynamics and a local center-stable manifold
to construct traveling modulating pulses with small oscillatory tails.
Here we determine how the edge lengths and reaction threshold
control the existence and number of stationary bistable fronts.

We restrict to profiles that agree on each pair of parallel edges.
Leaving the single-edge coordinate unchanged and doubling the parallel-edge coordinate
turns the stationary graph equation into
\begin{equation}
	\label{pblm}
	(p(x)u_x)_x+f(u,a)=0,
	\qquad x\in\mathbb{R},
\end{equation}
with limits
\begin{equation}
	\label{BC}
	\lim_{x\to-\infty}u(x)=1,
	\qquad
	\lim_{x\to+\infty}u(x)=0.
\end{equation}
We allow the more general $L$-periodic coefficient
\begin{equation}
	\label{func_p}
	p(x)=
	\begin{cases}
		1,        & x\in[kL,kL+L_1),     \\
		1/\alpha, & x\in[kL+L_1,(k+1)L),
	\end{cases}
	\qquad k\in\mathbb{Z},
\end{equation}
where $L=L_1+L_2$, $L_1,L_2>0$, and $0<\alpha<1$.
The necklace corresponds to $\alpha=1/4$,
with single-edge length $L_1$ and length $L_2/2$ on each parallel edge.
The weak formulation imposes continuity of $u$ and the flux $pu_x$;
Section~\ref{sec:switchingODE} gives the reduction and its spatial-dynamics formulation.

Our explicit results concern the McKean reaction
\begin{equation*}
	f_{\mathrm{MK}}(u,a)=
	\begin{cases}
		-u,  & u<a,     \\
		1-u, & u\geq a,
	\end{cases}
	\qquad a\in(0,1).
\end{equation*}
It is a piecewise linear counterpart of smooth bistable reactions
such as the Nagumo reaction $f_{\mathrm{Ng}}(u,a)=u(1-u)(u-a)$.
The spatial-dynamics construction in Section~\ref{sec:switchingODE}
also applies to continuous reactions satisfying Assumptions~\ref{cubic} and~\ref{ass:local_flow}.
We use the McKean reaction to determine pinning boundaries
at fixed diffusion contrast and edge lengths,
without a small-amplitude or small-period approximation.

To explain this choice, consider first the homogeneous equation $u_t=u_{xx}+f(u,a)$.
For smooth bistable reactions it admits a decreasing traveling front
with a unique speed and a profile unique up to translation
\cite{fifePhasePlaneDiscussion1981}.
For the McKean reaction these facts follow by matching the explicit solutions at $u=a$.
At the balanced parameter $a=1/2$, the reactions considered here satisfy
$\int_0^1f(u,1/2)\,du=0$.
The speed then vanishes, and the translates of the stationary profile
form a continuous family of stationary fronts.
Periodic variation of $p$ selects particular phases from this family.
For a smooth balanced reaction and a small perturbation $p=1+\varepsilon g$,
with real-valued $L$-periodic $g\in L^\infty(\mathbb{R})$ and $\varepsilon\|g\|_{L^\infty}<1$,
a formal expansion in $\varepsilon$ gives a first-order solvability condition
whose zeros are candidate phases.
Turning these candidates into stationary fronts requires a further existence argument.
Keener \cite{keenerHomogenizationPropagationBistable2000}
used averaging and a Melnikov calculation to predict propagation failure
for rapidly varying diffusion.
Near a stationary homogeneous front,
Dowdall, LeBlanc, and Lutscher \cite{dowdallInvasionPinningPeriodically2018}
used a center-manifold reduction to describe pinning
under small smooth periodic perturbations of diffusion or reaction.
Bastiaansen, Doelman, and Kaper \cite{bastiaansenMultifrontDynamicsSpatially2026}
extended Melnikov-type pinning calculations to patterns of several fronts
under a small smooth periodic heterogeneity, as discussed below.
For \eqref{func_p}, however, writing $p=1+\varepsilon g$
with $g$ the indicator of the high-diffusion intervals gives
$\varepsilon=\alpha^{-1}-1$,
so the necklace value $\alpha=1/4$ has $\varepsilon=3$.
The McKean reaction permits a direct calculation in this regime:
between coefficient interfaces and threshold crossings,
the stationary equation is affine with constant coefficients.

Propagation failure in bistable lattices has been studied since Keener \cite{Keener1987}.
Hupkes, Pelinovsky, and Sandstede \cite{HPS}
proved its absence under their normal-family hypotheses
when a smooth family of translated monotone stationary fronts exists.
For smooth reactions, Scheel and Van Vleck
\cite[Section~4]{scheelLatticeDifferentialEquations2009}
related lattice pinning to heterogeneous diffusion through a singular limit
of the stationary spatial return map.
They described pinned fronts by transverse heteroclinic intersections
and depinning by the unfolding of a quadratic tangency.

Piecewise-linear reactions permit more explicit lattice calculations.
For the McKean lattice, F\'{a}th \cite{Fa}
computed an explicit pinning interval and a monotone front unique
up to integer translations at each interior parameter.
Humphries, Moore, and Van Vleck \cite{HumphriesMooreVanVleck2011}
constructed McKean fronts with inhomogeneous discrete diffusion
and studied propagation failure when the diffusion coefficient
differs from its otherwise constant value at a single lattice connection.
Explicit McKean pinning boundaries are also available for layer solutions
on regular bi-infinite discrete trees
\cite{stehlikExplicitPropagationReversal2024}.
For continuous three-piece linear reactions,
Elmer constructed stationary lattice fronts \cite{Elmer2006}
and analyzed their stability \cite{Elmer2007}.
Chen, Mori, Ogawa, and Sintavanuruk \cite{ChenMoriOgawaSintavanuruk}
give an explicit construction of successive polygonal extensions
of the stable and unstable manifolds for a continuous three-piece reaction.

Since the McKean reaction is discontinuous at $u=a$,
one must ask whether the value assigned to $f_{\mathrm{MK}}$ at $u=a$ matters.
On a lattice it can.
A monotone lattice front can take the value $a$ at one site,
and the stationary equation at that site then holds only for a particular assigned value
\cite{ChenMoriOgawaSintavanuruk}.
In the necklace problem, the fronts cross $u=a$
at isolated points with nonzero flux (Lemma~\ref{lem:mckean_switching}),
so the value of $f_{\mathrm{MK}}$ at $a$ does not enter the weak equation
satisfied by these profiles.

Continuity of $u$ and $pu_x$ gives a continuous spatial flow
on $\mathbb{R}^2\setminus\{(a,0)\}$.
Its period map $\Psi^s$ advances the state $(u,pu_x)$
from $x=s$ to $x=s+L$.
A front corresponds to a heteroclinic orbit from $(1,0)$ to $(0,0)$,
independently of the sampling phase $s$
(Proposition~\ref{prop:heteroclinic_representation}).
A decreasing front joins two half-line solutions, converging to $(1,0)$ and to $(0,0)$,
at its threshold crossing, and their fluxes must agree there
(Proposition~\ref{prop:threshold_matching}).
Since the McKean equation is affine on each side of the threshold,
the states of these half-line solutions at a phase $s$ lie on two straight lines,
the unstable eigendirection through $(1,0)$ and the stable eigendirection through $(0,0)$
of the linearized period map.
The two lines meet where $u=b(s)$, for an explicit periodic function $b$, the matching function.
Hence a decreasing front crosses the threshold at phase $s$ precisely when $b(s)=a$
(Lemma~\ref{lem:matching_function}).
Thus the continuous spatial coordinate preserves a manifold-intersection description,
while the affine evolution makes the pinning boundary computable.

For fixed $(\alpha,L_1,L_2)$, the minimum and maximum of $b$
are explicit thresholds $a_-(\alpha,L_1,L_2)$ and $a_+(\alpha,L_1,L_2)$,
and Theorem~\ref{thm:monotone_classification} shows that a decreasing front exists precisely when
\begin{equation*}
	a_-(\alpha,L_1,L_2)\leq a\leq a_+(\alpha,L_1,L_2).
\end{equation*}
We call $[a_-,a_+]$ the \emph{pinning interval}.
Up to period translations, there are exactly two decreasing fronts at each interior threshold
and one at either endpoint.
Theorem~\ref{thm:all_front_existence} proves that no front of any kind exists outside this interval
and that every front at an endpoint is decreasing;
hence the same interval governs the existence of all fronts, monotone or not.
At an endpoint the two decreasing fronts coincide,
and the invariant manifolds meet in the nonsmooth contact
characterized by Lemma~\ref{lem:boundary_contact}.
For the evolution equation
\begin{equation}
	\label{pRD}
	u_t=(p(x)u_x)_x+f_{\mathrm{MK}}(u,a),
	\qquad x\in\mathbb{R},
\end{equation}
the interior front crossing in $(0,L_1)$ is Lyapunov stable
and the one crossing in $(L_1,L)$ is unstable
in the $W^{1,\infty}(\mathbb{R})$ norm
(Theorems~\ref{thm:stable_front_stability} and~\ref{thm:unstable_front_instability}).
These conclusions hold in the class of mild solutions with a unique threshold crossing
whose location is continuous in time.
We establish the required local existence, uniqueness and comparison theory in this class.
For $\alpha=1/4$, the coordinate reduction also identifies this evolution
with the necklace equation restricted to equal profiles on each pair of parallel edges.
The norm equivalence in Section~\ref{subsec:graph_transfer}
therefore transfers the stability conclusions to this graph solution class.

The two decreasing fronts also supply the choices for constructing nonmonotone fronts.
Choose an odd number $m\geq3$ of periods of the coefficient, far apart from one another,
and place one threshold crossing in each, alternately downward and upward.
The phase of each downward crossing lies near the crossing phase of one of the two decreasing fronts at threshold $a$.
The phase of each upward crossing lies near that of one of the two increasing fronts at threshold $a$,
which are the complements of the decreasing fronts at threshold $1-a$.
Theorem~\ref{thm:separated_fronts} gives exactly one front for each of the $2^m$ resulting choices,
and these are all the fronts with one crossing in each chosen period and none elsewhere.
Varying the chosen periods gives infinitely many inequivalent fronts for each odd $m$.
Consequently, nonmonotone fronts exist precisely when the two decreasing fronts coexist, that is, for $a_-<a<a_+$.
For thresholds a fixed distance away from $a_-$ and $a_+$,
the required distance between the chosen periods can be chosen independently of $m$
(Corollary~\ref{cor:uniform_central_fronts}).
At $a=1/2$, a symmetric family with three crossings has an explicit sharp minimum distance between the chosen periods
(Proposition~\ref{prop:balanced_minimum_gap}).
The same principle, placing well-separated fronts independently at single-front pinning positions,
underlies the perturbative construction of Bastiaansen, Doelman, and Kaper
\cite{bastiaansenMultifrontDynamicsSpatially2026}
for the Allen--Cahn equation with a small, smooth, periodic topographic heterogeneity.
There the pinning positions are given to leading order by simple zeros of a Melnikov function,
and the small amplitude also lets them compute the leading terms of the small eigenvalues of the linearization,
which decide the spectral stability of sufficiently separated patterns.
In Theorem~\ref{thm:separated_fronts}, the diffusion contrast is fixed,
and the only small quantity is the interaction between separated crossings.
The theorem also shows that the $2^m$ fronts exhaust the fronts with one crossing in each chosen period and none elsewhere,
and Theorem~\ref{thm:all_front_existence} excludes nonmonotone fronts outside $(a_-,a_+)$.
We do not prove stability of the nonmonotone fronts.
Section~\ref{sec:discussion} discusses what is expected and what a proof would need.

Section~\ref{sec:switchingODE} develops the spatial return map.
Section~\ref{sec:McKean} classifies the decreasing fronts and proves that the pinning interval governs all fronts.
Section~\ref{sec:stability} establishes temporal stability and instability,
and Section~\ref{sec:nonmonotone} constructs the nonmonotone fronts.
The appendices contain the proofs for the switching flow and its local geometry,
the evolution estimates, and the calculation behind Proposition~\ref{prop:balanced_minimum_gap}.

\section{Periodically switching spatial dynamics}\label{sec:switchingODE}

We express stationary fronts as heteroclinic orbits of a period map.
We impose the following conditions on the reaction.
The McKean reaction and smooth bistable reactions such as the Nagumo reaction satisfy them.

\begin{assumption}
	\label{cubic}
	For each $a\in(0,1)$,
	the function $f(\cdot,a):\mathbb{R}\to\mathbb{R}$ is locally bounded
	and satisfies the following conditions:
	\begin{itemize}
		\item
		      $f(\cdot,a)$ is continuously differentiable on $\mathbb{R}\setminus\{a\}$;
		      it may be discontinuous at $u=a$.
		\item
		      $f(u,a)>0$ for $u\in(-\infty,0)\cup(a,1)$.
		\item
		      $f(u,a)<0$ for $u\in(0,a)\cup(1,\infty)$.
		\item
		      $f_u(0,a)<0$ and $f_u(1,a)<0$.
		\item
		      $f(1-u,1/2)=-f(u,1/2)$ for all $u\neq1/2$.
	\end{itemize}
\end{assumption}

We write $f(u)=f(u,a)$ when the parameter $a$ is fixed.
Since $p$ is discontinuous,
we interpret \eqref{pblm} in the weak sense:
$u\in H^1_{\mathrm{loc}}(\mathbb{R})$ is a weak solution if
\begin{equation}
	\int_{\mathbb{R}}\bigl(-p(x)u_x\varphi_x(x)+f(u,a)\varphi(x)\bigr)\,dx=0
\end{equation}
for every $\varphi\in C_c^\infty(\mathbb{R})$.
The weak formulation gives the following interface conditions.

\begin{proposition}
	\label{c1gap}
	The following statements are equivalent for $u\in H^1_{\mathrm{loc}}(\mathbb{R})$.
	\begin{enumerate}
		\item [(1)] $u$ is a weak solution of \eqref{pblm} with \eqref{func_p}.
		\item [(2)] $u\in C^0(\mathbb{R})$; on each coefficient interval,
		      $u_x$ has a locally absolutely continuous representative and
		      \begin{equation*}
			      \begin{aligned}
				      u_{xx}+f(u,a)&=0,
				      &&x\in(kL,kL+L_1), \\
				      \frac{1}{\alpha}u_{xx}+f(u,a)&=0,
				      &&x\in(kL+L_1,(k+1)L),
			      \end{aligned}
		      \end{equation*}
		      almost everywhere, for every $k\in\mathbb{Z}$;
		      and the one-sided limits of $u_x$ satisfy the interface conditions
		      \begin{equation*}
			      \begin{aligned}
				      \frac{1}{\alpha}u_x(kL^-)&=u_x(kL^+), \\
				      u_x\bigl((kL+L_1)^-\bigr)&=\frac{1}{\alpha}u_x\bigl((kL+L_1)^+\bigr),
			      \end{aligned}
			      \qquad k\in\mathbb{Z}.
		      \end{equation*}
	\end{enumerate}
	If $f(\cdot,a)$ is continuous, the equations in {\it (2)} hold classically.
\end{proposition}

To obtain \eqref{pblm} from the necklace graph,
restrict to stationary profiles that agree on each pair of parallel edges.
For $j\in\mathbb{Z}$,
write $u^{2j}$ for the profile on the single edge
$\Omega_1=(0,L_1)$
and $u^{2j+1}$ for the common profile on the following pair of edges,
each identified with $\Omega_2=(0,L_2/2)$.
The edge equations are
\begin{equation}
	\label{necklace}
	\begin{aligned}
		u^{2j}_{xx}+f(u^{2j})&=0,
		&&x\in\Omega_1, \\
		u^{2j+1}_{xx}+f(u^{2j+1})&=0,
		&&x\in\Omega_2.
	\end{aligned}
\end{equation}
For each $j\in\mathbb{Z}$,
continuity at the vertices requires
\begin{equation*}
	u^{2j}(L_1)=u^{2j+1}(0),
	\qquad
	u^{2j-1}(L_2/2)=u^{2j}(0).
\end{equation*}
The profiles also satisfy the Kirchhoff conditions
\begin{equation}
	\label{NK_BC}
	\begin{aligned}
		u^{2j}_x(L_1)&=2u^{2j+1}_x(0), \\
		2u^{2j-1}_x(L_2/2)&=u^{2j}_x(0).
	\end{aligned}
\end{equation}
These equations and vertex conditions are equivalent to \eqref{pblm} with $\alpha=1/4$.
Indeed, on each pair of parallel edges,
introduce the stretched coordinate $x=2y$,
where $y\in(0,L_2/2)$ is the edge coordinate,
and set $u(x)=u^{2j+1}(x/2)$,
up to the corresponding translation of the coordinate origin.
Then the equation on the doubled-edge portion becomes
\begin{equation*}
	4u_{xx}+f(u)=0,
\end{equation*}
and the Kirchhoff conditions \eqref{NK_BC} become continuity of the flux
$p u_x$ with $p=4$.

More generally, if each pair is replaced by $n\geq2$ parallel edges of equal length $\ell$,
the stretching $x=ny$ turns the Kirchhoff conditions into continuity of $pu_x$ with $p=n^2$,
that is, $\alpha=1/n^2$ and $L_2=n\ell$.

The same stretching applies to the evolution equation.
Denote the necklace graph by $\Gamma$.
If $u_\Gamma(y,t)$ is the common value of a graph profile on $\Gamma$ at the coordinate $y$ of a parallel edge
and $u(x,t)=u_\Gamma(x/2,t)$, then $\partial_tu_\Gamma=\partial_y^2u_\Gamma+f(u_\Gamma,a)$ becomes
\begin{equation*}
    u_t=4u_{xx}+f(u,a),
    \qquad \partial_yu_\Gamma(y,t)=2u_x(2y,t),
\end{equation*}
and on a single edge the coordinate is unchanged.
Thus, for graph profiles that agree on each pair of parallel edges,
the necklace evolution becomes the scalar equation $u_t=(pu_x)_x+f(u,a)$ with $\alpha=1/4$,
which is \eqref{pRD} for the McKean reaction.
Section~\ref{subsec:graph_transfer} shows that this identification extends
to the mild solutions used in Section~\ref{sec:stability}.

For the scalar problem \eqref{pblm},
use the continuous flux as the second phase-plane coordinate:
\[
	v:=p(x)u_x.
\]
By Proposition~\ref{c1gap}, a solution of \eqref{pblm}
satisfies the periodically switching system
\begin{equation}
	\label{E1}
	\left\{
	\begin{aligned}
		u' & = v,     \\
		v' & = -f(u),
	\end{aligned}
	\right. \qquad (x\in(k L, kL+L_1))
\end{equation}
and
\begin{equation}
	\label{E2}
	\left\{
	\begin{aligned}
		u' & = \alpha v, \\
		v' & = -f(u),
	\end{aligned}
	\right. \qquad (x\in(kL+L_1,(k+1)L)),
\end{equation}
where $k\in\mathbb{Z}$.
The equilibria $(0,0)$ and $(1,0)$ are hyperbolic
for both systems by Assumption~\ref{cubic}.

At each coefficient interface $kL$ or $kL+L_1$,
we switch between \eqref{E1} and \eqref{E2}
while keeping both $u$ and $v$ continuous.
Solutions are understood to be locally absolutely continuous curves
satisfying the differential equations almost everywhere.

To define their flows, we also assume the following.

\begin{assumption}[Local flow]
	\label{ass:local_flow}
	The autonomous systems \eqref{E1} and \eqref{E2}
	admit unique local solutions on $\mathbb{R}^2\setminus\{(a,0)\}$
	depending continuously on time and initial data,
	and no nonconstant solution reaches $(a,0)$.
\end{assumption}

If $f(\cdot,a)\in C^1(\mathbb{R})$, then $f(a,a)=0$ by Assumption~\ref{cubic},
so $(a,0)$ is an equilibrium of both systems,
and Assumption~\ref{ass:local_flow} follows from local Lipschitz continuity of the vector fields.
For the McKean reaction, Lemma~\ref{lem:mckean_switching} establishes it.
Denote the local flows of \eqref{E1} and \eqref{E2} by $\phi^1(t)$ and $\phi^2(t)$, respectively.
Under Assumption~\ref{ass:local_flow}, a nonconstant solution with $u=a$ at some point has $v\neq0$ there,
so it crosses $u=a$ transversally, and the points where $u=a$ are isolated.

For each $s\in[0,L)$,
let $\mathcal{D}^s\subset\mathbb{R}^2\setminus\{(a,0)\}$
denote the set of initial values
for which the corresponding solution of the switching system
exists over one full period $[s,s+L]$.
This domain can be smaller than $\mathbb{R}^2\setminus\{(a,0)\}$,
because a spatial trajectory may blow up before completing one period.
For example,
when $f=f_{\mathrm{Ng}}(\cdot,1/2)$,
the solution of \eqref{E1} with
$u(0)=R>1$ and $v(0)=R(R-1)/\sqrt{2}$
satisfies $u'=u(u-1)/\sqrt{2}$
and blows up at $x=T_R:=\sqrt{2}\log(R/(R-1))$.
For sufficiently large $R$,
we have $T_R<L_1$,
so blow-up occurs before the first coefficient switch.

We define the period map
$\Psi^s:\mathcal{D}^s\to\mathbb{R}^2$
by assigning to each initial value at $x=s$
the corresponding solution value at $x=s+L$.
Explicitly,
\begin{equation}
	\label{Lmap-s}
	\Psi^s:=
	\begin{cases}
		\phi^1(s)\circ\phi^2(L_2)\circ\phi^1(L_1-s),
		 & s\in[0,L_1), \\
		\phi^2(s-L_1)\circ\phi^1(L_1)\circ\phi^2(L-s),
		 & s\in[L_1,L).
	\end{cases}
\end{equation}
In particular,
\[
	\Psi^0=\phi^2(L_2)\circ\phi^1(L_1),
	\qquad
	\Psi^{L_1}=\phi^1(L_1)\circ\phi^2(L_2).
\]
The points $(0,0)$ and $(1,0)$ belong to $\mathcal{D}^s$
for every $s\in[0,L)$.

\begin{lemma}\label{lem:hyperbolicity}
	Under Assumption~\ref{cubic},
	both $(0,0)$ and $(1,0)$ are hyperbolic fixed points of $\Psi^s$
	for every $s\in[0,L)$.
\end{lemma}

\begin{proof}
	Fix $a\in(0,1)$ and $u_*\in\{0,1\}$.
	By Assumption~\ref{cubic},
	$f(u_*,a)=0$ and $c:=-f_u(u_*,a)>0$.
	Since $f(\cdot,a)$ is $C^1$ near $u_*$,
	the finite-time flow maps in \eqref{Lmap-s}
	are defined and $C^1$ for initial data sufficiently close to $(u_*,0)$.
	Hence $\Psi^s$ is $C^1$ near this fixed point
	for every $s\in[0,L)$.
	The linearizations of \eqref{E1} and \eqref{E2} at $(u_*,0)$ are
	\[
		A_1=
		\begin{pmatrix}
			0 & 1 \\
			c & 0
		\end{pmatrix},
		\qquad
		A_2=
		\begin{pmatrix}
			0 & \alpha \\
			c & 0
		\end{pmatrix}.
	\]
	For $s=0$,
	the linearized period map is
	\[
		M=e^{L_2A_2}e^{L_1A_1}.
	\]
	Since $\operatorname{tr}A_1=\operatorname{tr}A_2=0$,
	we have $\det M=1$.
	A direct calculation gives
	\[
		\operatorname{tr}M
		=
		2\cosh(\sqrt{c}L_1)
		\cosh(\sqrt{\alpha c}L_2)
		+
		\left(
		\sqrt{\alpha}+\frac{1}{\sqrt{\alpha}}
		\right)
		\sinh(\sqrt{c}L_1)
		\sinh(\sqrt{\alpha c}L_2)
		>2.
	\]
	Thus $M$ has two positive real eigenvalues
	$\lambda^{-1}<1<\lambda$.

	For a general phase $s\in[0,L)$,
	let $B(s)$ denote the evolution matrix of the linearized system
	from $x=0$ to $x=s$:
	\begin{equation*}
		B(s)=
		\begin{cases}
			e^{sA_1},
			 & s\in[0,L_1), \\
			e^{(s-L_1)A_2}e^{L_1A_1},
			 & s\in[L_1,L).
		\end{cases}
	\end{equation*}
	By periodicity,
	evolution from $0$ to $s+L$
	can be expressed as either $D\Psi^s(u_*,0)B(s)$
	or $B(s)M$,
	according to whether we split the interval at $s$ or at $L$.
	Since $B(s)$ is invertible,
	we obtain
	\begin{equation*}
		D\Psi^s(u_*,0)=B(s)MB(s)^{-1}.
	\end{equation*}
	Thus the linearized period maps are conjugate,
	and their eigenvalues are $\lambda^{-1}$ and $\lambda$,
	independently of $s$.
\end{proof}

Since $\Psi^s$ need not be differentiable away from the equilibria,
we define its invariant manifolds through convergence of sampled states.

\begin{definition}[Stable and unstable manifolds]
	\label{def:invariant_manifolds}
	Fix $s\in[0,L)$.
	The \emph{unstable manifold} of $(1,0)$ for $\Psi^s$ is the set of states $Z(s)$
	of solutions $Z$ of the switching system on $(-\infty,s]$
	with $Z(s-nL)\to(1,0)$ as $n\to\infty$.
	The \emph{stable manifold} of $(0,0)$ for $\Psi^s$ is the set of states $Z(s)$
	of solutions $Z$ on $[s,\infty)$ with $Z(s+nL)\to(0,0)$ as $n\to\infty$.
	The unstable manifold of $(0,0)$ and the stable manifold of $(1,0)$
	are defined in the same way, with the two equilibria exchanged.
\end{definition}

For a $C^1$ reaction these are the usual global stable and unstable manifolds
of the hyperbolic fixed points in Lemma~\ref{lem:hyperbolicity}.
For the McKean reaction they need not be smooth curves,
and we use the names for the sets in Definition~\ref{def:invariant_manifolds}.

Sampling a stationary front at intervals of length $L$
gives a heteroclinic orbit of $\Psi^s$.
Conversely, such an orbit determines the profile between its sampled states.

\begin{proposition}
	\label{prop:heteroclinic_representation}
	The following statements are equivalent.
	\begin{enumerate}
		\item
		      The problem \eqref{pblm} with \eqref{BC}
		      has a solution.
		\item
		      There exists $s\in[0,L)$
		      and a sequence $\{Z_n\}_{n\in\mathbb{Z}}$,
		      with $Z_n\in\mathcal{D}^s$,
		      such that
		      \[
			      Z_{n+1}=\Psi^s(Z_n),
			      \qquad
			      \lim_{n\to-\infty}Z_n=(1,0),
			      \qquad
			      \lim_{n\to+\infty}Z_n=(0,0).
		      \]
		\item
		      For every $s\in[0,L)$,
		      there exists a sequence satisfying the conditions in {\it (2)}.
	\end{enumerate}
\end{proposition}

\begin{proof}
	Suppose first that {\it (1)} holds.
	For any $s\in[0,L)$,
	define $Z_n:=\bigl(u(s+nL),v(s+nL)\bigr)$,
	where $v$ denotes the continuous representative of the flux $pu_x$.
	Since a front is nonconstant, Assumption~\ref{ass:local_flow},
	applied on successive coefficient intervals, shows that its state $(u,v)$ never equals $(a,0)$.
	Hence $Z_n\in\mathcal D^s$ for every $n$.
	By periodicity and the definition of the period map,
	\[
		Z_{n+1}=\Psi^s(Z_n).
	\]
	The far-field condition gives
	\[
		u(s+nL)\to1
		\quad\text{as }n\to-\infty,
		\qquad
		u(s+nL)\to0
		\quad\text{as }n\to+\infty.
	\]
	Moreover,
	$v=pu_x$ satisfies $v'=-f(u)$
	and is uniformly continuous near both infinities.
	If $v$ failed to tend to zero at either infinity,
	uniform continuity would give intervals tending to that infinity,
	of a fixed positive length,
	on which $v$ has a constant sign and $|v|$ is bounded below by a positive constant.
	Since $u'=p^{-1}v$ and $p^{-1}\geq\alpha$,
	the change in $u$ across each such interval would be bounded away from zero,
	contradicting its convergence.
	Thus $v(x)\to0$ as $x\to\pm\infty$,
	and {\it (3)} holds.

	The implication {\it (3)}$\Rightarrow${\it (2)} is immediate.

	Finally,
	suppose that {\it (2)} holds.
	For each $n\in\mathbb{Z}$,
	let $Z^{(n)}$ be the solution of \eqref{E1}--\eqref{E2}
	on $[s+nL,s+(n+1)L]$
	with initial value $Z^{(n)}(s+nL)=Z_n$.
	Existence throughout this interval follows from
	$Z_n\in\mathcal{D}^s$ and the periodicity of $p$.
	By the definition of the period map,
	\begin{equation*}
		Z^{(n)}(s+(n+1)L)
		=\Psi^s(Z_n)
		=Z_{n+1}
		=Z^{(n+1)}(s+(n+1)L).
	\end{equation*}
	These pieces define a global trajectory $(u,v)$
	by setting $(u(x),v(x))=Z^{(n)}(x)$
	for $x\in[s+nL,s+(n+1)L]$.
	The resulting pair is locally absolutely continuous,
	with $v=pu_x$ and $v'=-f(u)$ almost everywhere.
	Integration by parts,
	using the continuity of $v$ at every interface,
	shows that $u$ is a weak solution of \eqref{pblm}.

	To obtain the far-field limits \eqref{BC},
	we control each trajectory segment by its initial state.
	Fix $E_q=(q,0)$ with $q\in\{0,1\}$.
	Since $f(q)=0$ and $f$ is continuously differentiable near $q$,
	choose $r,C>0$ such that
	$|f(z)|\leq C|z-q|$ whenever $|z-q|\leq r$,
	and $\|p^{-1}\|_{L^\infty(\mathbb{R})}\leq C$.

	For all sufficiently large positive $n$ when $q=0$,
	or sufficiently negative $n$ when $q=1$,
	we have $|Z_n-E_q|\leq (r/2)e^{-CL}$.
	Fix such an $n$ and write $x_n:=s+nL$.
	The equations \eqref{E1}--\eqref{E2} give
	$|(Z^{(n)})'(x)|\leq C|Z^{(n)}(x)-E_q|$ for almost every $x\in[x_n,x_n+L]$
	before the first exit from $B_r(E_q)$.
	Gr\"onwall's inequality therefore gives
	\begin{equation*}
		|Z^{(n)}(x)-E_q|
		\leq e^{C(x-x_n)}|Z_n-E_q|,
	\end{equation*}
	on the same interval.
	A first exit at $x_*\leq x_n+L$ would imply
	$r\leq e^{CL}|Z_n-E_q|\leq r/2$,
	a contradiction.
	Hence
	\begin{equation*}
		\sup_{x\in[x_n,x_n+L]}
		|(u(x),v(x))-E_q|
		\leq e^{CL}|Z_n-E_q|.
	\end{equation*}

	Applying this estimate with $E_q=(1,0)$ as $n\to-\infty$
	and with $E_q=(0,0)$ as $n\to+\infty$,
	we obtain $u(x)\to1$ as $x\to-\infty$
	and $u(x)\to0$ as $x\to+\infty$.
	Therefore, {\it (1)} holds.
\end{proof}

\begin{lemma}
	\label{lem:symmetry}
	Assume $a=1/2$,
	and let $\mathcal R(u,v):=(1-u,v)$.
	For each $s\in[0,L)$,
	let $\bar s\in[0,L)$ be the unique phase satisfying
	\begin{equation*}
		\bar s\equiv L_1-s\pmod L.
	\end{equation*}
	Then $\mathcal R$ maps the unstable manifold of $(1,0)$
	for $\Psi^s$ onto the stable manifold of $(0,0)$
	for $\Psi^{\bar s}$.
	In particular,
	these manifolds are symmetric with respect to $u=1/2$
	for the same period map when
	$s=L_1/2$ or $s=L_1+L_2/2$.
\end{lemma}

\begin{proof}
	Fix $s\in[0,L)$
	and set $c:=s+\bar s$.
	Since $c-L_1\in L\mathbb Z$,
	periodicity and \eqref{func_p} give
	\begin{equation*}
		p(c-x)=p(L_1-x)=p(x)
		\qquad\text{for almost every }x.
	\end{equation*}
	Assumption~\ref{cubic} gives $f(1-u)=-f(u)$ for $u\neq1/2$,
	and for every nonconstant solution, Assumption~\ref{ass:local_flow} makes the points where $u=1/2$ isolated.

	Let $Z_0$ belong to the unstable manifold of $(1,0)$
	for $\Psi^s$.
	The trajectory $Z=(u,v)$ with $Z(s)=Z_0$
	exists on $(-\infty,s]$
	and satisfies $Z(s-nL)\to(1,0)$ as $n\to\infty$.
	Define
	\begin{equation*}
		\widetilde Z(x)
		:=
		\mathcal R Z(c-x)
		=
		\bigl(1-u(c-x),v(c-x)\bigr),
		\qquad x\geq \bar s.
	\end{equation*}
	This curve is locally absolutely continuous,
	and its components satisfy
	\begin{equation*}
		\begin{aligned}
			\widetilde u'(x)
			 & =u'(c-x)
			=\frac{v(c-x)}{p(c-x)}
			=\frac{\widetilde v(x)}{p(x)}, \\
			\widetilde v'(x)
			 & =-v'(c-x)
			=f(u(c-x))
			=-f(\widetilde u(x))
		\end{aligned}
	\end{equation*}
	almost everywhere.
	Thus $\widetilde Z$ solves the switching system,
	with both components continuous at the coefficient interfaces.
	Furthermore,
	\begin{equation*}
		\widetilde Z(\bar s+nL)
		=\mathcal R Z(s-nL)
		\longrightarrow(0,0)
		\qquad\text{as }n\to\infty.
	\end{equation*}
	Hence $\mathcal R Z_0$ belongs to the stable manifold of $(0,0)$
	for $\Psi^{\bar s}$.
	The reverse inclusion follows by applying the same reflection again.

	Finally,
	$\bar s=s$ precisely when $2s\equiv L_1\pmod L$.
	The two solutions in $[0,L)$ are
	$s=L_1/2$ and $s=L_1+L_2/2$.
\end{proof}

By Proposition~\ref{prop:heteroclinic_representation},
the intersections of the unstable manifold of $(1,0)$
and the stable manifold of $(0,0)$ for $\Psi^s$
are the states at $x=s$ of stationary fronts.
For reactions smooth in $(u,a)$,
transverse heteroclinic intersections persist
under sufficiently small changes in $(a,L_1,L_2)$.
A local mechanism for their disappearance is a \emph{heteroclinic tangency},
at which the stable and unstable curves meet with a common tangent line.
If the contact is quadratic
and unfolds nondegenerately with respect to one parameter,
two nearby transverse intersections coalesce
and disappear as that parameter passes through its critical value.

For the McKean reaction, the invariant curves may fail to be differentiable,
so we compare the relative positions of their branches
in Section~\ref{sec:McKean}.
For fixed $\alpha$,
Theorem~\ref{thm:monotone_classification}
characterizes the existence region for decreasing stationary fronts
in $(a,L_1,L_2)$-space.
Lemma~\ref{lem:boundary_contact} describes the contact of the invariant branches
at its boundaries $a=a_-$ and $a=a_+$,
which is a nonsmooth counterpart of the tangency described above.
Section~\ref{sec:stability} relates these intersections to temporal stability.

\section{Stationary fronts for the McKean reaction}\label{sec:McKean}

We take $f=f_{\mathrm{MK}}(\cdot,a)$ with $a\in(0,1)$,
$0<\alpha<1$, and $L_1,L_2>0$.
A decreasing front crosses the threshold $u=a$ exactly once.
On the left of the crossing it is a solution of the switching system converging to $(1,0)$,
and on the right one converging to $(0,0)$.
Away from the threshold the switching system is affine,
so these two half-line solutions are explicit.
At a phase $s$, their states lie on two straight lines,
through $(1,0)$ and through $(0,0)$,
in eigendirections of the linearized period map (Section~\ref{subsec:equilibrium_tails}).
Joining the two half-line solutions at $u=a$ gives a front precisely when their fluxes agree there,
that is, when the two lines meet on the threshold line $u=a$.
The lines meet where $u=b(s)$, for an explicit periodic function $b$, the matching function.
Hence a decreasing front crosses the threshold at phase $s$ precisely when $b(s)=a$
(Section~\ref{subsec:matching_function}).
The range and monotonicity of $b$ give the existence interval and the number of decreasing fronts
(Section~\ref{subsec:decreasing_classification}),
and Section~\ref{subsec:all_stationary_fronts} shows that the same interval governs all fronts, monotone or not.
Section~\ref{subsec:bounded_solutions} collects bounded solutions of the linear equation used in Sections~\ref{sec:stability} and~\ref{sec:nonmonotone}.

We call two profiles $u$ and $\widetilde u$
\emph{$L$-translationally equivalent}
if $u(x)=\widetilde u(x+kL)$ for every $x\in\mathbb R$
and some $k\in\mathbb Z$.
Theorem~\ref{thm:monotone_classification} counts the decreasing fronts
up to this equivalence.
The local ordering of their invariant branches, established at the end of this section,
will be used in the stability analysis of Section~\ref{sec:stability}.

\subsection{Switching flow and explicit evolution}
\label{subsec:switching_flow}

For a weak stationary solution,
the flux $v=pu_x$ satisfies $v'=-f_{\mathrm{MK}}(u,a)$ in the distributional sense,
so by Proposition~\ref{c1gap} both $u$ and $v$ are continuous,
including at coefficient interfaces and threshold crossings.

On an interval where $p^{-1}=\kappa\in\{1,\alpha\}$,
the two regimes are
\begin{equation*}
	u'=\kappa v,
	\qquad
	v'=
	\begin{cases}
		u,   & u<a, \\
		u-1, & u>a.
	\end{cases}
\end{equation*}
At a threshold crossing $(a,v_0)$ with $v_0\neq0$,
we continue into $u>a$ if $v_0>0$
and into $u<a$ if $v_0<0$,
with both components continuous.
The regimes are reversed when evolving backward.

\begin{lemma}
	\label{lem:mckean_switching}
	Let $a\in(0,1)$ and $f=f_{\mathrm{MK}}(\cdot,a)$.
	For each fixed $\kappa\in\{1,\alpha\}$,
	the preceding rule defines a continuous flow on
	$\Omega_a:=\mathbb{R}^2\setminus\{(a,0)\}$,
	defined for all positive and negative times.
	Every threshold crossing is transverse.
	No solution of the switching system passes through $(a,0)$.
\end{lemma}

The proof, including continuity at threshold crossings,
is given in Appendix~\ref{app:mckean_switching}.

The McKean reaction satisfies Assumption~\ref{cubic},
and Lemma~\ref{lem:mckean_switching} gives Assumption~\ref{ass:local_flow},
so the results of Section~\ref{sec:switchingODE} apply to it.
Since every solution exists for all $x$,
the period maps \eqref{Lmap-s} are defined on $\mathcal{D}^s=\Omega_a$ for every $s\in[0,L)$.

For the explicit evolution formulas, set
\begin{equation*}
	H_L:=\{(u,v)\in\mathbb{R}^2:u<a\},
	\qquad
	H_R:=\{(u,v)\in\mathbb{R}^2:u>a\}.
\end{equation*}
On an interval where $p=1$,
the two affine systems are
\begin{align}
	u' & =v,\qquad v'=u,
	   &                   & (u,v)\in H_L,
	\label{E1left}                         \\
	u' & =v,\qquad v'=u-1,
	   &                   & (u,v)\in H_R.
	\label{E1right}
\end{align}
On an interval where $p=\alpha^{-1}$,
they are
\begin{align}
	u' & =\alpha v,\qquad v'=u,
	   &                          & (u,v)\in H_L,
	\label{E2left}                                \\
	u' & =\alpha v,\qquad v'=u-1,
	   &                          & (u,v)\in H_R.
	\label{E2right}
\end{align}

Extend each affine vector field to $\mathbb{R}^2$,
and denote its flow by $\phi_l^j$ or $\phi_r^j$,
according to the side of the threshold,
where $j=1,2$ identifies the coefficient interval.
The linear flows are represented by the matrices
\begin{equation}
	\label{time_t}
	\phi_l^1(t)=
	\begin{pmatrix}
		\cosh t & \sinh t \\
		\sinh t & \cosh t
	\end{pmatrix}
\end{equation}
and
\begin{equation}
	\label{time_t_alpha}
	\phi_l^2(t)=
	\begin{pmatrix}
		\cosh(\sqrt{\alpha}t)
		 & \sqrt{\alpha}\sinh(\sqrt{\alpha}t) \\
		\alpha^{-1/2}\sinh(\sqrt{\alpha}t)
		 & \cosh(\sqrt{\alpha}t)
	\end{pmatrix},
	\qquad t\in\mathbb{R}.
\end{equation}
Writing $E_1:=(1,0)^{\mathsf T}$,
the corresponding affine flows are
\begin{equation*}
	\phi_r^j(t)Z
	=E_1+\phi_l^j(t)(Z-E_1),
	\qquad j\in\{1,2\}.
\end{equation*}
These formulas do not depend on $a$;
the threshold determines which affine regime applies.

Let $Z_0\in\Omega_a$, $j\in\{1,2\}$, and $t\neq0$.
If $\phi^j(\tau)Z_0\in H_L$
for every $\tau$ strictly between $0$ and $t$,
then $\phi^j(t)Z_0=\phi_l^j(t)Z_0$.
The analogous statement holds with $H_R$ and $\phi_r^j$.
This follows from uniqueness for the affine initial-value problem
and continuity at the endpoints.

\subsection{Equilibrium tails and the straight parts of the invariant manifolds}
\label{subsec:equilibrium_tails}

An \emph{equilibrium tail} is a solution of the switching system
on a half-line that converges to an equilibrium at infinity.
We call it stable if it approaches $(0,0)$ as $x\to+\infty$,
and unstable if it approaches $(1,0)$ as $x\to-\infty$.
A tail that stays on one side of the threshold
solves one of the affine systems of Section~\ref{subsec:switching_flow} on its whole half-line.
Such tails are therefore described by the linear system
\begin{equation}
	\label{eq:linear_tail_system}
	y'=p(x)^{-1}z,
	\qquad z'=y.
\end{equation}
For $s\in\mathbb R$, let $A_s$ be its \emph{monodromy matrix} based at $x=s$,
which maps the state of a solution at $x=s$ to its state at $x=s+L$.
It is obtained from \eqref{Lmap-s}
by replacing each $\phi^j$ with $\phi_l^j$; in particular,
\begin{equation*}
	A_0=\phi_l^2(L_2)\phi_l^1(L_1).
\end{equation*}
For initial data sufficiently close to either equilibrium,
the switched trajectory remains in the corresponding affine regime
throughout one period.
Consequently,
\begin{equation*}
	A_s=D\Psi^s(0,0)=D\Psi^s(1,0).
\end{equation*}
The calculation in Lemma~\ref{lem:hyperbolicity}, with $c=1$,
gives the phase-independent eigenvalues $\lambda^{-1}$ and $\lambda$,
where $\lambda>1$.
The explicit matrices \eqref{time_t}--\eqref{time_t_alpha}
show that $A_s$ has strictly positive entries.
Since $(A_s)_{12}>0$, no eigenvector has zero first component.

A solution of \eqref{eq:linear_tail_system} that is multiplied by a constant over each period
is called a \emph{Floquet solution}.
Let $Y_-$ and $Y_+$ be the solutions of \eqref{eq:linear_tail_system} on $\mathbb R$
whose states at $x=0$ are stable and unstable eigenvectors of $A_0$, respectively, normalized by
\begin{equation}
	\label{eq:normalized_floquet_solutions}
	Y_\pm=(y_\pm,z_\pm)^{\mathsf T},
	\qquad y_+(0)=y_-(0)=1.
\end{equation}
Their second components are the continuous fluxes $z_\pm=py_\pm'$ almost everywhere.
Periodicity of $p$ and uniqueness give
\begin{equation}
	\label{eq:tail_floquet_relations}
	Y_-(x+L)=\lambda^{-1}Y_-(x),
	\qquad
	Y_+(x+L)=\lambda Y_+(x).
\end{equation}
Thus $Y_-$ decays exponentially as $x\to+\infty$,
and $Y_+$ decays exponentially as $x\to-\infty$.
The tails below are affine images of $Y_-$ and $Y_+$,
and the same two solutions reappear in the Green kernel of Section~\ref{subsec:bounded_solutions}.

The first components are positive on $\mathbb R$.
Indeed, if $y_-(x_0)=0$ for some $x_0$, then
\begin{equation*}
	0=[y_-z_-]_{x_0}^{\infty}
	=\int_{x_0}^{\infty}
	\bigl(p^{-1}z_-^2+y_-^2\bigr)\,dx.
\end{equation*}
This would force $Y_-\equiv0$ by uniqueness,
contradicting $y_-(0)=1$.
Hence $y_->0$ on $\mathbb R$.
The same argument on $(-\infty,x_0]$ gives $y_+>0$.
Integrating $z_\pm'=y_\pm$ from the respective decay endpoints gives
\begin{equation*}
	z_-(x)=-\int_x^{\infty}y_-(r)\,dr<0,
	\qquad
	z_+(x)=\int_{-\infty}^x y_+(r)\,dr>0.
\end{equation*}
Since $y_\pm'=p^{-1}z_\pm$,
$y_-$ is strictly decreasing and $y_+$ is strictly increasing.
We may therefore define the positive functions
\begin{equation}
	\label{eq:floquet_slopes}
	r_+(s):=\frac{z_+(s)}{y_+(s)},
	\qquad
	r_-(s):=-\frac{z_-(s)}{y_-(s)},
	\qquad s\in\mathbb R.
\end{equation}
Equivalently, $r_+=py_+'/y_+$ and $r_-=-py_-'/y_-$.

\begin{lemma}[Equilibrium tails]
	\label{lem:equilibrium_tails}
	Fix $a\in(0,1)$ and $s\in\mathbb R$.
	The vectors
	\begin{equation*}
		\boldsymbol w_s=\frac{Y_+(s)}{y_+(s)}=\begin{pmatrix}1\\r_+(s)\end{pmatrix},
		\qquad
		\boldsymbol v_s=\frac{Y_-(s)}{y_-(s)}=\begin{pmatrix}1\\-r_-(s)\end{pmatrix}
	\end{equation*}
	are unstable and stable eigenvectors of $A_s$, respectively.
	The initial states at $x=s$ whose forward trajectories
	remain in $u<a$ for $x>s$
	and converge to $(0,0)$ as $x\to+\infty$
	are precisely
	\begin{equation*}
		\{\mu\boldsymbol v_s:\mu\leq a\}.
	\end{equation*}
	The initial states whose backward trajectories
	remain in $u>a$ for $x<s$
	and converge to $(1,0)$ as $x\to-\infty$
	are precisely
	\begin{equation*}
		\{E_1+(\mu-1)\boldsymbol w_s:\mu\geq a\},
		\qquad E_1=(1,0)^{\mathsf T}.
	\end{equation*}
	In particular,
	there is a unique tail of each kind with $u(s)=a$.
	These two tails are strictly decreasing on their respective half-lines
	and have strictly negative flux.
\end{lemma}

\begin{proof}
	The period identities \eqref{eq:tail_floquet_relations}
	show that $Y_-(s)$ and $Y_+(s)$ are stable and unstable eigenvectors of $A_s$.
	Their positive first components permit the stated rescaling.
	For $\mu\leq a$,
	the curve $\mu Y_-(x)/y_-(s)$ on $[s,\infty)$
	therefore remains in $u<a$ for $x>s$
	and converges to $(0,0)$.
	For $\mu\geq a$,
	the curve $E_1+(\mu-1)Y_+(x)/y_+(s)$ on $(-\infty,s]$
	remains in $u>a$ for $x<s$
	and converges to $(1,0)$.
	Both curves solve the corresponding McKean system
	on their respective half-lines.

	Conversely,
	any trajectory with either stated regime and limit
	solves the corresponding affine system throughout that half-line.
	Convergence of its sampled states forces its initial displacement
	to lie in the stable eigenspace of $A_s$ in the forward case,
	or the unstable eigenspace in the backward case.
	Its initial $u$-coordinate then fixes the multiplier,
	giving exactly the two sets above.
	Taking $\mu=a$ gives the unique threshold tails.
	They are positive multiples of $Y_-$ and affine images $E_1-cY_+$ with $c>0$,
	so the signs of $z_\pm$ and the monotonicity of $y_\pm$ give the last assertion.
\end{proof}

The two sets in Lemma~\ref{lem:equilibrium_tails}
are the parts of two straight lines,
\begin{equation}
	\label{eq:eigenlines}
	\ell^{\mathrm s}_s:=\{(u,-r_-(s)u):u\in\mathbb R\},
	\qquad
	\ell^{\mathrm u}_s:=\{(u,-r_+(s)(1-u)):u\in\mathbb R\},
\end{equation}
the stable eigendirection through $(0,0)$ and the unstable eigendirection through $(1,0)$.
The stable manifold of $(0,0)$ for $\Psi^s$ contains the part of $\ell^{\mathrm s}_s$ with $u\leq a$,
and the unstable manifold of $(1,0)$ contains the part of $\ell^{\mathrm u}_s$ with $u\geq a$.
We call these the \emph{straight parts} of the two manifolds.
The numbers $r_+(s)$ and $-r_-(s)$ are the slopes of $\ell^{\mathrm u}_s$ and $\ell^{\mathrm s}_s$.
The lines themselves do not depend on $a$; the threshold determines only where the straight parts end.
By the argument of Lemma~\ref{lem:equilibrium_tails}, the states at $x=s$
whose backward trajectories remain in $u<a$ for $x<s$ and converge to $(0,0)$ as $x\to-\infty$
are the points $\mu\boldsymbol w_s$ with $\mu\leq a$.
We call this half-line the straight part of the unstable manifold of $(0,0)$.

\subsection{Threshold fluxes and the matching function}
\label{subsec:matching_function}

Let $v^{\mathrm s}(s,a)$ and $v^{\mathrm u}(s,a)$ denote the fluxes at $x=s$
of the stable and unstable equilibrium tails with $u(s)=a$,
respectively.
They are the $v$-coordinates of $\ell^{\mathrm s}_s$ and $\ell^{\mathrm u}_s$ at $u=a$:
\begin{equation}
	\label{eq:threshold_flux_asymmetric}
	v^{\mathrm s}(s,a)=-a\,r_-(s),
	\qquad
	v^{\mathrm u}(s,a)=-(1-a)\,r_+(s).
\end{equation}

\begin{proposition}[Characterization by threshold fluxes]
	\label{prop:threshold_matching}
	Fix $a\in(0,1)$.
	The nonincreasing weak solutions of \eqref{pblm} with \eqref{BC},
	up to translations by integer multiples of $L$,
	correspond bijectively to the phases $s\in[0,L)$ satisfying
	\begin{equation}
		\label{eq:threshold_matching}
		v^{\mathrm u}(s,a)=v^{\mathrm s}(s,a).
	\end{equation}
	Every such solution is strictly decreasing
	and has a unique threshold crossing.
	The representative associated with $s$ is uniquely determined
	by $u(s)=a$.
\end{proposition}

\begin{proof}
	Suppose first that \eqref{eq:threshold_matching} holds.
	Lemma~\ref{lem:equilibrium_tails}
	gives a unique unstable tail on $(-\infty,s]$
	and a unique stable tail on $[s,\infty)$,
	both taking the value $a$ at $s$.
	Their fluxes at $s$ are $v^{\mathrm u}(s,a)$ and $v^{\mathrm s}(s,a)$,
	so joining them gives a continuous profile with continuous flux.
	The joined profile belongs to $H^1_{\mathrm{loc}}(\mathbb R)$
	and solves the stationary equation on either side of $s$.
	Integration by parts across $s$
	shows that it is a weak solution on $\mathbb R$.
	Lemma~\ref{lem:equilibrium_tails} gives the limits \eqref{BC}
	and strict decrease on both half-lines,
	hence on all of $\mathbb R$.

	Conversely,
	let $u$ be a nonincreasing weak solution with the limits \eqref{BC}.
	Continuity gives a point $\rho$ with $u(\rho)=a$.
	If there were two such points,
	monotonicity would force $u\equiv a$ between them,
	with flux identically zero,
	contrary to Lemma~\ref{lem:mckean_switching}.
	Thus $\rho$ is unique,
	with $u>a$ on $(-\infty,\rho)$ and $u<a$ on $(\rho,\infty)$.
	As in the proof of Proposition~\ref{prop:heteroclinic_representation},
	the continuous flux $pu_x$ tends to zero at both infinities.
	Write $\rho=s+kL$ with $s\in[0,L)$ and $k\in\mathbb Z$.
	By periodicity,
	$u(\cdot+kL)$ is a solution whose threshold crossing is at $s$.
	Its restrictions to the two half-lines
	are therefore the tails in Lemma~\ref{lem:equilibrium_tails}.
	Continuity of its flux gives \eqref{eq:threshold_matching}.

	Uniqueness of the tails gives uniqueness of the front
	with its crossing fixed at $s$.
	Two fronts are $L$-translationally equivalent
	precisely when their crossing phases agree.
	This proves the claimed bijection.
\end{proof}

The two lines \eqref{eq:eigenlines} have slopes of opposite signs,
so they meet at exactly one point.
Equating $-r_+(s)(1-u)$ and $-r_-(s)u$ shows that this point has first coordinate
\begin{equation}
	\label{eq:matching_function}
	b(s):=\frac{r_+(s)}{r_+(s)+r_-(s)},
	\qquad s\in\mathbb R.
\end{equation}
At $u=a$, the height of the unstable line above the stable line is
\begin{equation}
	\label{eq:flux_mismatch}
	v^{\mathrm u}(s,a)-v^{\mathrm s}(s,a)
	=\bigl(r_+(s)+r_-(s)\bigr)\bigl(a-b(s)\bigr),
\end{equation}
by \eqref{eq:threshold_flux_asymmetric} and \eqref{eq:matching_function}.
Thus \eqref{eq:threshold_matching} holds precisely when $b(s)=a$,
that is, when the two lines meet on the threshold line $u=a$.
The common point is then the threshold state of the front,
where the straight parts of the two invariant manifolds meet.
For each phase $s$, $b(s)$ is the only threshold
for which a decreasing front crosses the threshold at $s$.
We call $b$ the \emph{matching function}.
Lemma~\ref{lem:step_profiles} gives a second description of it,
as the value $b(\rho)$ at $\rho$ of the bounded solution with forcing $\mathbf 1_{\{x<\rho\}}$.

\subsection{Classification of decreasing stationary fronts}
\label{subsec:decreasing_classification}

To count the solutions of $b(s)=a$, we need the dependence of the slopes $r_\pm$ on the phase.
For $s\in\mathbb R$, let $\bar s\in[0,L)$ satisfy $\bar s\equiv L_1-s\pmod L$,
as in Lemma~\ref{lem:symmetry}.

\begin{lemma}[Slopes of the eigendirections]
	\label{lem:floquet_slopes}
	Let $0<\alpha<1$.
	The functions $r_\pm$ are continuous and $L$-periodic,
	and continuously differentiable on each coefficient interval, where
	\begin{equation}
		\label{eq:riccati}
		r_+'=1-\kappa r_+^2,
		\qquad
		r_-'=-1+\kappa r_-^2,
		\qquad \kappa:=p^{-1}.
	\end{equation}
	Moreover, $r_-(s)=r_+(\bar s)$ and $1<r_\pm(s)<\alpha^{-1/2}$ for every $s\in\mathbb R$.
\end{lemma}

\begin{proof}
	For $0\leq s\leq L$, the explicit linear flows give
	\begin{equation*}
		Y_\pm(s)=
		\begin{cases}
			\phi_l^1(s)Y_\pm(0),                  & 0\leq s\leq L_1, \\
			\phi_l^2(s-L_1)\phi_l^1(L_1)Y_\pm(0), & L_1\leq s\leq L.
		\end{cases}
	\end{equation*}
	Hence $r_\pm$ are continuous and continuously differentiable on each coefficient interval.
	By \eqref{eq:tail_floquet_relations}, each quotient in \eqref{eq:floquet_slopes}
	takes the same value at $s$ and $s+L$, so $r_\pm$ are $L$-periodic.
	Differentiating the quotients with \eqref{eq:linear_tail_system} gives \eqref{eq:riccati}.

	The reflection $x\mapsto L_1-x$ preserves $p$ almost everywhere.
	Thus $x\mapsto\bigl(y_+(L_1-x),-z_+(L_1-x)\bigr)$ solves \eqref{eq:linear_tail_system},
	and it decays as $x\to+\infty$ because $Y_+$ decays as $x\to-\infty$.
	By \eqref{eq:tail_floquet_relations}, the solutions that decay as $x\to+\infty$
	are the multiples of $Y_-$; comparing first components at $x=0$ gives
	\begin{equation}
		\label{eq:floquet_reflection}
		y_-(x)=\frac{y_+(L_1-x)}{y_+(L_1)},
		\qquad
		z_-(x)=-\frac{z_+(L_1-x)}{y_+(L_1)}.
	\end{equation}
	Hence $r_-(s)=r_+(L_1-s)=r_+(\bar s)$.

	It remains to bound $r_+$.
	The positive equilibria of the scalar equation $r'=1-\kappa r^2$
	are $1$ on $(0,L_1)$ and $\alpha^{-1/2}$ on $(L_1,L)$.
	Uniqueness prevents a solution from crossing the corresponding equilibrium
	within either interval.
	If $r_+(0)\leq1$,
	then $r_+(0)\leq r_+(L_1)\leq1$
	and $r_+(L)>r_+(L_1)$,
	contradicting $r_+(L)=r_+(0)$.
	If $r_+(0)\geq\alpha^{-1/2}$, then $r_+(L_1)<r_+(0)$.
	When $r_+(L_1)\geq\alpha^{-1/2}$,
	the second interval gives $r_+(L)\leq r_+(L_1)<r_+(0)$;
	otherwise, it gives $r_+(L)<\alpha^{-1/2}\leq r_+(0)$.
	Both alternatives again contradict periodicity.
	Therefore $1<r_+(0)<\alpha^{-1/2}$.
	The scalar equations preserve this open interval on both parts of the period,
	so $1<r_+(s)<\alpha^{-1/2}$ for every $s$.
	The bound for $r_-$ follows from $r_-(s)=r_+(\bar s)$.
\end{proof}

The extreme values of the matching function are explicit.
Set
\begin{equation*}
	\begin{aligned}
		\tau
		 & :=2\cosh L_1\cosh(\sqrt{\alpha}L_2)
		+\bigl(\sqrt{\alpha}+\alpha^{-1/2}\bigr)
		\sinh L_1\sinh(\sqrt{\alpha}L_2),            \\
		d
		 & :=\bigl(\alpha^{-1/2}-\sqrt{\alpha}\bigr)
		\sinh L_1\sinh(\sqrt{\alpha}L_2),
	\end{aligned}
\end{equation*}
and define
\begin{equation}
	\label{eq:front_existence_thresholds}
	a_\pm(\alpha,L_1,L_2)
	:=\frac12\pm\frac{d}{2\sqrt{\tau^2-4}}.
\end{equation}
By Lemma~\ref{lem:hyperbolicity} with $c=1$, $\tau>2$.

\begin{lemma}[The matching function]
	\label{lem:matching_function}
	Let $0<\alpha<1$ and $L_1,L_2>0$.
	The function $b$ is continuous, $L$-periodic, and continuously differentiable
	on each coefficient interval, where
	\begin{equation}
		\label{eq:matching_derivative}
		b'=\frac{1-\kappa r_+r_-}{r_++r_-}.
	\end{equation}
	It is strictly decreasing on $[0,L_1]$ and strictly increasing on $[L_1,L]$,
	and its one-sided derivatives at $0$ and $L_1$ are nonzero.
	It satisfies $b(\bar s)=1-b(s)$,
	\begin{equation}
		\label{eq:matching_extremes}
		b(0)=\max b=a_+,
		\qquad
		b(L_1)=\min b=a_-,
	\end{equation}
	and $0<a_-<1/2<a_+<1$.
	In particular, for each $\theta\in(a_-,a_+)$
	the equation $b(s)=\theta$ has exactly two solutions modulo $L$,
	one in $(0,L_1)$ and one in $(L_1,L)$, and both are simple.
\end{lemma}

\begin{proof}
	Lemma~\ref{lem:floquet_slopes} gives continuity, periodicity and regularity,
	and \eqref{eq:riccati} gives
	\begin{equation*}
		b'=\frac{r_+'r_--r_+r_-'}{(r_++r_-)^2}
		=\frac{(r_++r_-)(1-\kappa r_+r_-)}{(r_++r_-)^2},
	\end{equation*}
	which is \eqref{eq:matching_derivative}.
	Since $1<r_+r_-<\alpha^{-1}$, the numerator $1-\kappa r_+r_-$ is negative where $\kappa=1$
	and positive where $\kappa=\alpha$.
	By continuity of $r_\pm$, the same strict signs hold for the one-sided limits at $0$ and $L_1$.
	This proves the monotonicity and the assertion on one-sided derivatives.
	The reflection identity $r_-(s)=r_+(\bar s)$ gives $b(\bar s)=r_-(s)/(r_-(s)+r_+(s))=1-b(s)$.

	For the extreme values, consider the monodromy matrix $A_0=\phi_l^2(L_2)\phi_l^1(L_1)$.
	Multiplying the matrices \eqref{time_t}--\eqref{time_t_alpha} gives
	\begin{equation*}
		\tau=\operatorname{tr}A_0,
		\qquad
		d=(A_0)_{22}-(A_0)_{11}>0.
	\end{equation*}
	The matrix $A_0$ has strictly positive entries and determinant $1$.
	Set $D:=\sqrt{\tau^2-4}$. Then
	\begin{equation*}
		D^2=d^2+4(A_0)_{12}(A_0)_{21}>d^2,
	\end{equation*}
	and the eigenvalues of $A_0$ are $\lambda^{\pm1}=(\tau\pm D)/2$.
	By \eqref{eq:normalized_floquet_solutions},
	$(1,r_+(0))^{\mathsf T}$ and $(1,-r_-(0))^{\mathsf T}$
	are eigenvectors of $A_0$ for $\lambda$ and $\lambda^{-1}$.
	The first rows of the eigenvalue equations give
	\begin{equation*}
		r_+(0)=\frac{\lambda-(A_0)_{11}}{(A_0)_{12}}=\frac{d+D}{2(A_0)_{12}},
		\qquad
		r_-(0)=\frac{(A_0)_{11}-\lambda^{-1}}{(A_0)_{12}}=\frac{D-d}{2(A_0)_{12}}.
	\end{equation*}
	Hence $b(0)=(D+d)/(2D)=a_+$.
	Since $\bar 0=L_1$, the reflection identity gives $b(L_1)=1-a_+=a_-$.
	The inequality $0<d<D$ gives $0<a_-<1/2<a_+<1$.
	The final assertion follows from strict monotonicity
	and the nonzero derivatives on the two open coefficient intervals.
\end{proof}

\begin{theorem}[Classification of decreasing stationary fronts]
	\label{thm:monotone_classification}
	Fix $0<\alpha<1$ and $L_1,L_2>0$,
	and let $a_\pm=a_\pm(\alpha,L_1,L_2)$ be given by \eqref{eq:front_existence_thresholds}.
	Then $0<a_-<1/2<a_+<1$.
	For $a\in(0,1)$,
	a nonincreasing weak solution of \eqref{pblm} with \eqref{BC}
	exists if and only if $a_-\leq a\leq a_+$.
	Up to translations by integer multiples of $L$,
	there are exactly two such fronts for $a_-<a<a_+$
	and exactly one at either endpoint.
	Their threshold-crossing phases are the solutions of $b(s)=a$ in $[0,L)$.
	In the interior,
	one front has its crossing phase in $(0,L_1)$
	and the other in $(L_1,L)$.
	At $a=a_-$ the crossing phase is $L_1$,
	and at $a=a_+$ it is $0$.
\end{theorem}

\begin{proof}
	By Proposition~\ref{prop:threshold_matching} and \eqref{eq:flux_mismatch},
	the decreasing fronts, up to translations by integer multiples of $L$,
	correspond bijectively to the phases $s\in[0,L)$ with $b(s)=a$.
	Lemma~\ref{lem:matching_function} shows that $b$ takes values in $[a_-,a_+]$,
	attains $a_-$ only at $L_1$ and $a_+$ only at $0$ modulo $L$,
	and attains each value in $(a_-,a_+)$ exactly once in each of $(0,L_1)$ and $(L_1,L)$.
\end{proof}

\begin{figure}[!htb]
	\centering
	\includegraphics[width=\textwidth]{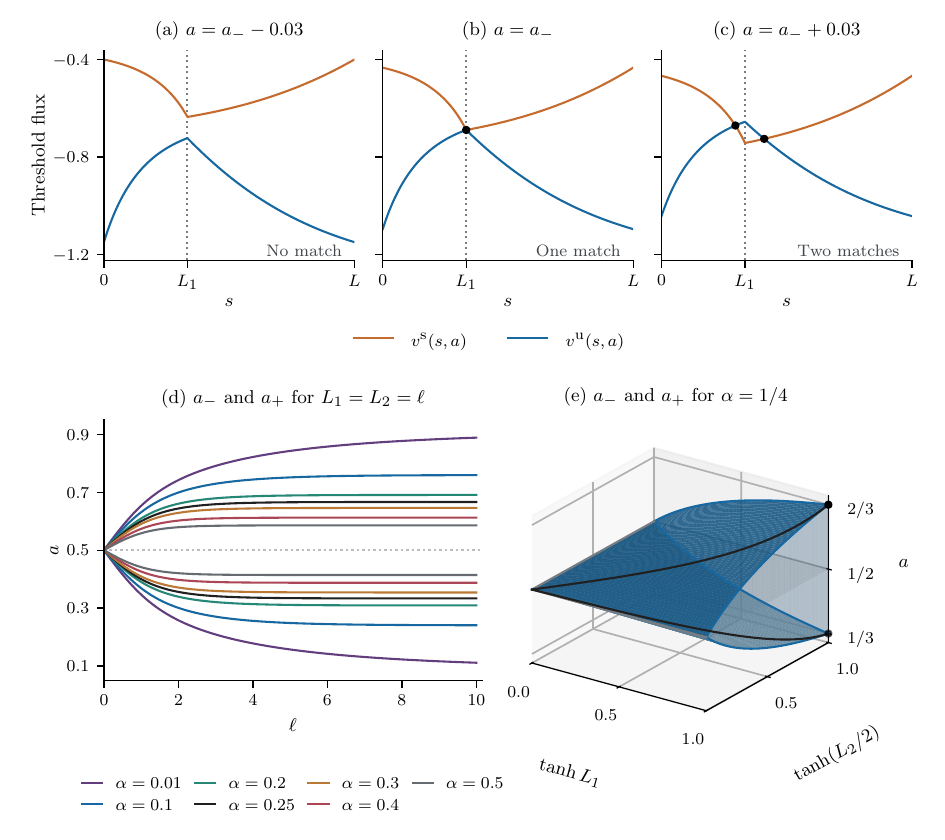}
	\caption{Threshold matching and the existence interval for decreasing fronts,
		counted up to translations by integer multiples of $L$.
		(a)--(c) Stable and unstable threshold fluxes $v^{\mathrm s}(s,a)$ (orange)
		and $v^{\mathrm u}(s,a)$ (blue), for the necklace value $\alpha=1/4$,
		with $L_1=0.8$ and $L_2=1.6$.
		The thresholds are $a_-\simeq0.385721$ and $a_+\simeq0.614279$.
		Dots mark flux matches, each of which determines one decreasing front.
		The dotted vertical line marks $s=L_1$.
		(d) The curves $a=a_\pm(\alpha,\ell,\ell)$
		for equal coefficient lengths $L_1=L_2=\ell$,
		with $0.01\leq\ell\leq10$ and the indicated values of $\alpha$.
		For each fixed $\alpha$ and $\ell$,
		there are two decreasing fronts between the two corresponding curve values,
		one at either value, and none outside the interval.
		(e) The pinning interval for $\alpha=1/4$ over $t_1=\tanh L_1$ and $t_2=\tanh(L_2/2)$.
		It is the region between the surfaces $a=a_-$ and $a=a_+$,
		which meet at $a=1/2$ on the edges $t_1=0$ and $t_2=0$.
		The corner $t_1=t_2=1$ gives the limits $1/3$ and $2/3$.
		The black curve is the equal-length case of (d).}
	\label{fig:monotone_summary}
\end{figure}

Figure~\ref{fig:monotone_summary} illustrates the flux matching
and the resulting existence interval.

For $a_-<a<a_+$,
let $u^P$ and $u^Q$ denote the two fronts
in Theorem~\ref{thm:monotone_classification},
with representatives chosen so that their unique threshold crossings satisfy
\begin{equation*}
	u^P(\rho_P)=u^Q(\rho_Q)=a,
	\qquad
	\rho_P\in(0,L_1),
	\quad
	\rho_Q\in(L_1,L).
\end{equation*}
These conditions fix the representatives uniquely.
When the threshold varies, we write $\rho_P(\theta)$ and $\rho_Q(\theta)$
for the corresponding crossing phases at threshold $\theta\in(a_-,a_+)$;
the argument is omitted when the threshold is fixed at $a$.

\begin{remark}
	\label{rem:balanced_phases}
	For $a=1/2$, the phases $s=L_1/2$ and $s=L_1+L_2/2$ satisfy $\bar s=s$,
	so $b(s)=1-b(s)=1/2$ there by Lemma~\ref{lem:matching_function}.
	Hence $\rho_P=L_1/2$ and $\rho_Q=L_1+L_2/2$.
\end{remark}

For $a=a_-$ and $a=a_+$,
Lemma~\ref{lem:boundary_contact} describes the intersection
of the unstable manifold of $(1,0)$
and the stable manifold of $(0,0)$
for $\Psi^{L_1}$ and $\Psi^0$, respectively.

The width of the pinning interval depends on the two lengths through
$t_1:=\tanh L_1$ and $t_2:=\tanh(\sqrt\alpha L_2)$.
Put $k:=\alpha^{-1/2}-\alpha^{1/2}>0$ and $c:=\alpha^{-1/2}+\alpha^{1/2}$.
Dividing $d$ and $\sqrt{\tau^2-4}$ in \eqref{eq:front_existence_thresholds}
by $\cosh L_1\cosh(\sqrt\alpha L_2)$ and using $c^2-4=k^2$ gives
\begin{equation*}
a_+-a_-
=\frac{d}{\sqrt{\tau^2-4}}
=\frac{k}{\sqrt{k^2+4t_1^{-2}+4t_2^{-2}+4c\,(t_1t_2)^{-1}}}.
\end{equation*}
Hence, for fixed $\alpha$, the pinning interval widens strictly as either length increases,
and $a_-$ and $a_+$ move symmetrically about $1/2$.
After continuous extension to $(t_1,t_2)\in[0,1]^2$, the region between the graphs of $a_-$ and $a_+$ is dorayaki-shaped.\footnote{A dorayaki is a Japanese confection of two round pancakes joined along their rims.}
The two surfaces meet at $a=1/2$ along the edges where one length vanishes
and separate as both lengths grow (Figure~\ref{fig:monotone_summary}(e)).
When $L_1$ and $L_2$ both tend to infinity, independently of each other,
\begin{equation*}
a_-\longrightarrow\frac{\sqrt{\alpha}}{1+\sqrt{\alpha}},
\qquad
a_+\longrightarrow\frac{1}{1+\sqrt{\alpha}}.
\end{equation*}
When only one length tends to infinity, the width tends to $kt/(2+ct)$,
where $t$ is the value of $t_1$ or $t_2$ for the fixed length.
This limit is smaller than the width $(1-\sqrt\alpha)/(1+\sqrt\alpha)$ obtained when both lengths are long.
When either length tends to zero, the width tends to zero,
since it is at most $k\min\{t_1,t_2\}/2$,
and both thresholds tend to $1/2$.
For $\alpha=n^{-2}$, that is, for $n$ parallel edges, these limits have graph interpretations.
When both lengths are long, the vertices are far apart,
and near each of them the necklace looks like an isolated junction of one edge with $n$ edges.
With equal profiles on the $n$ edges,
a decreasing McKean front at a junction where it passes from the single edge into the $n$ edges
exists exactly for $a\in[1/(n+1),1/2)$,
and one passing from the $n$ edges into the single edge exists exactly for $a\in(1/2,n/(n+1)]$.
Each bead of the necklace has a junction of each orientation, one at either end.
The limiting interval $[1/(n+1),n/(n+1)]$ is the union of these two ranges and the balanced value $a=1/2$,
at which the fronts lie on long homogeneous edges away from the vertices.
Its endpoints are the critical flux balances $1-a=na$ and $n(1-a)=a$,
in agreement with the junction-blocking criterion of Jimbo and Morita \cite{JM1,JM2}.
For $\alpha=n^{-2}$ and $L_2\to\infty$ with $L_1$ fixed, the width tends to
\begin{equation*}
\mathcal W_n(L_1)=\frac{(n-n^{-1})\tanh L_1}{2+(n+n^{-1})\tanh L_1},
\end{equation*}
which increases strictly with $L_1$ from $0$ to $(n-1)/(n+1)$.
Near each single edge, the necklace then looks like the graph studied by Iwasaki, Jimbo, and Morita \cite{IJM},
in which one edge of length $L_1$ joins two bundles of $n$ half-lines.
For $a\neq1/2$ in $(1/(n+1),n/(n+1))$, the limit of the pinning interval contains $a$ exactly when $L_1\geq L_1^*(a)$,
where $\mathcal W_n(L_1^*(a))=|2a-1|$.
This is the McKean counterpart of the minimum connector length
that Iwasaki, Jimbo, and Morita obtain for smooth bistable reactions under their hypotheses on the reaction \cite[Theorem~1.2]{IJM}.
Figure~\ref{fig:monotone_summary}(d) shows the equal-length case $L_1=L_2=\ell$.

\subsection{A necessary condition for all stationary fronts}
\label{subsec:all_stationary_fronts}

Theorem~\ref{thm:monotone_classification} classifies the decreasing fronts.
We now show that no front of any kind exists for $a\notin[a_-,a_+]$,
and that at $a=a_\pm$ the decreasing front is the only one.
A front with several threshold crossings has flux of both signs,
so Proposition~\ref{prop:threshold_matching} does not apply to it.
Its first and last crossings, however, are still joined to equilibrium tails.
After the first crossing the front continues below the threshold,
and the stable eigendirection decides whether it can return to the threshold.
Before the last crossing the same holds above the threshold, backward in $x$.

\begin{theorem}[Necessary condition for stationary fronts]
	\label{thm:all_front_existence}
	Fix $0<\alpha<1$, $L_1,L_2>0$ and $a\in(0,1)$.
	If \eqref{pblm} with \eqref{BC} has a weak solution, then $a\in[a_-,a_+]$.
	If that solution is not monotone, then $a_-<a<a_+$.
	If $a=a_-$ or $a=a_+$, every such solution is,
	up to translation by an integer multiple of $L$,
	the decreasing front of Theorem~\ref{thm:monotone_classification}.
\end{theorem}

\begin{proof}
	Let $u$ be a weak stationary front and $v=pu_x$ its continuous flux.
	As in the proof of Proposition~\ref{prop:heteroclinic_representation},
	$v(x)\to0$ at both infinities,
	so the state $(u,v)$ tends to $(1,0)$ as $x\to-\infty$ and to $(0,0)$ as $x\to+\infty$.

	The threshold set $\{x:u(x)=a\}$ is nonempty and compact.
	By Lemma~\ref{lem:mckean_switching}, the flux is nonzero at each of its points,
	so $u$ is strictly monotone near each of them.
	Hence the threshold points are isolated, and there are finitely many.
	Let $s$ and $t$ be the first and the last.
	Then $u>a$ on $(-\infty,s)$ and $u<a$ on $(t,\infty)$,
	and Lemma~\ref{lem:equilibrium_tails} identifies the two exterior pieces as the threshold tails:
	\begin{equation*}
		\bigl(u(s),v(s)\bigr)=\bigl(a,-(1-a)r_+(s)\bigr),
		\qquad
		\bigl(u(t),v(t)\bigr)=\bigl(a,-a\,r_-(t)\bigr).
	\end{equation*}

	We first show that $b(s)\leq a$.
	Just to the right of $s$ the front lies below the threshold,
	where the switching system is the linear system \eqref{eq:linear_tail_system}.
	Decomposing the state at $s$ in the eigenvectors $\boldsymbol w_s$ and $\boldsymbol v_s$,
	and using $b(s)\bigl(r_+(s)+r_-(s)\bigr)=r_+(s)$,
	shows that, as long as $u<a$,
	\begin{equation*}
		\begin{pmatrix}u(x)\\v(x)\end{pmatrix}
		=\bigl(a-b(s)\bigr)\frac{Y_+(x)}{y_+(s)}
		+b(s)\frac{Y_-(x)}{y_-(s)}.
	\end{equation*}
	Suppose that $a<b(s)$.
	The first coefficient is then negative and the second positive.
	Since $z_+>0$ and $z_-<0$, the right-hand side has strictly negative flux for every $x\geq s$,
	so its first component stays below $a$ for all $x>s$.
	It therefore agrees with $u$ on $[s,\infty)$.
	The growing term $(a-b(s))y_+(x)/y_+(s)$ gives $u(x)\to-\infty$ as $x\to+\infty$,
	a contradiction.
	Hence $b(s)\leq a$.
	If $b(s)=a$, the right-hand side is the stable tail $aY_-(x)/y_-(s)$,
	which never returns to the threshold.
	Then $s=t$, and $u$ is the decreasing front with crossing at $s$.

	For the last crossing, set $w:=1-u$ and $q:=-v$ on the interval to the left of $t$ where $u>a$.
	There $(w,q)$ solves the linear system \eqref{eq:linear_tail_system},
	with $(w(t),q(t))=(1-a,a\,r_-(t))$.
	The same decomposition gives, as long as $u>a$,
	\begin{equation*}
		\begin{pmatrix}w(x)\\q(x)\end{pmatrix}
		=\bigl(1-b(t)\bigr)\frac{Y_+(x)}{y_+(t)}
		+\bigl(b(t)-a\bigr)\frac{Y_-(x)}{y_-(t)}.
	\end{equation*}
	Suppose that $a>b(t)$.
	The first coefficient is then positive and the second negative,
	so $q>0$ for every $x\leq t$.
	Hence $w<1-a$, that is, $u>a$, for all $x<t$,
	and the formula holds on $(-\infty,t]$.
	As $x\to-\infty$, the growing term $(b(t)-a)y_-(x)/y_-(t)$ gives $w(x)\to-\infty$,
	that is, $u(x)\to+\infty$, a contradiction.
	Hence $a\leq b(t)$.
	If $b(t)=a$, the formula is the unstable tail,
	which has no threshold crossing to the left of $t$.
	Then $s=t$, and $u$ is again the decreasing front.

	The bounds $b(s)\leq a\leq b(t)$ and \eqref{eq:matching_extremes} give
	\begin{equation}
		\label{eq:all_front_threshold_bounds}
		a_-\leq b(s)\leq a\leq b(t)\leq a_+.
	\end{equation}
	In particular, $a\in[a_-,a_+]$.
	If $u$ is not monotone, then it has more than one threshold crossing,
	so the equality cases above are excluded and $b(s)<a<b(t)$.
	This gives $a_-<a<a_+$.
	Finally, if $a=a_-$, then \eqref{eq:all_front_threshold_bounds} forces $b(s)=a$,
	and if $a=a_+$, it forces $b(t)=a$.
	In either case $u$ is the decreasing front,
	which is unique up to translation by an integer multiple of $L$
	by Theorem~\ref{thm:monotone_classification}.
\end{proof}

In the phase plane at phase $s$, the proof compares the threshold state $(a,v^{\mathrm u}(s,a))$,
which lies on the unstable eigendirection $\ell^{\mathrm u}_s$,
with the stable eigendirection $\ell^{\mathrm s}_s$.
A state just below the threshold that lies below $\ell^{\mathrm s}_s$ escapes to $u\to-\infty$,
one on $\ell^{\mathrm s}_s$ follows the stable tail,
and one above $\ell^{\mathrm s}_s$ returns to the threshold (Remark~\ref{rem:spargel_departure}).
Since the two lines meet where $u=b(s)$, the state $(a,v^{\mathrm u}(s,a))$ lies above $\ell^{\mathrm s}_s$ precisely when $a>b(s)$.

With Theorem~\ref{thm:monotone_classification}, which supplies a decreasing front for every $a\in[a_-,a_+]$,
Theorem~\ref{thm:all_front_existence} shows that stationary fronts exist exactly for $a\in[a_-,a_+]$.
The curves in Figure~\ref{fig:monotone_summary}(d) therefore bound the existence of all stationary fronts.
The count in the figure caption concerns decreasing fronts.
Theorem~\ref{thm:separated_fronts} constructs nonmonotone fronts at every interior threshold.

\subsection{Bounded solutions with prescribed forcing}
\label{subsec:bounded_solutions}

Sections~\ref{sec:stability} and~\ref{sec:nonmonotone} use bounded solutions of the linear equation $KV=q$, where
\begin{equation*}
	KV:=-(pV_x)_x+V
\end{equation*}
is understood in the distributional sense.
A stationary McKean front $u$ solves $Ku=\mathbf 1_{\{u\geq a\}}$.
If the set where $u\geq a$ is prescribed, this equation is linear in $u$.
In Section~\ref{sec:stability}, comparison profiles near $u^P$ are built from the solutions $V_\rho$ with forcing $\mathbf 1_{\{x<\rho\}}$,
and uniqueness of bounded solutions identifies limits.
In Section~\ref{sec:nonmonotone}, a front with prescribed crossings is an alternating sum of such solutions,
and the identity $\partial_\rho V_\rho(x)=G(x,\rho)$ for the Green kernel $G$ below controls the interaction of separated crossings.
The value $V_\rho(\rho)$ also gives a second description of the matching function.

The Floquet solutions \eqref{eq:normalized_floquet_solutions} satisfy $y_\pm'=p^{-1}z_\pm$ and $z_\pm'=y_\pm$,
so their Wronskian
\begin{equation*}
	\omega:=z_+(x)y_-(x)-y_+(x)z_-(x)
\end{equation*}
is independent of $x$; it is positive because $y_\pm>0$ and $z_-<0<z_+$.

\begin{lemma}[Positive bounded inverse]
	\label{lem:positive_bounded_inverse}
	Define
	\begin{equation}
		\label{eq:positive_green_kernel}
		G(x,y):=\frac{1}{\omega}\,
		y_+\bigl(\min\{x,y\}\bigr)\,
		y_-\bigl(\max\{x,y\}\bigr),
		\qquad x,y\in\mathbb R.
	\end{equation}
	Then $G$ is continuous and positive, and
	\begin{equation}
		\label{eq:green_decay}
		G(x,y)+|p(x)\partial_xG(x,y)|
		\leq C e^{-\gamma|x-y|},
		\qquad
		\gamma:=\frac{\log\lambda}{L}>0,
	\end{equation}
	for every $y$ and almost every $x$, where $C$ depends only on $\alpha,L_1,L_2$.
	For every real-valued $q\in L^\infty(\mathbb R)$, the function
	\begin{equation}
		\label{eq:bounded_green_representation}
		V(x)=\int_{\mathbb R}G(x,y)q(y)\,dy
	\end{equation}
	is the unique bounded weak solution of $KV=q$, and it belongs to $W^{1,\infty}(\mathbb R)$.
	If $q_1\leq q_2$ almost everywhere,
	the corresponding bounded solutions satisfy $V_1\leq V_2$ everywhere.
	If, in addition, $q_1<q_2$ on a set of positive measure,
	then $V_1<V_2$ everywhere.
\end{lemma}

\begin{proof}
	Continuity and positivity of $G$ follow from those of $y_\pm$.
	For fixed $y$, the function $G(\cdot,y)$ is continuous
	and solves the homogeneous equation away from $y$.
	Its flux is continuous at coefficient interfaces different from $y$,
	and at $x=y$ the flux jump is
	\begin{equation*}
		[p\partial_xG(\cdot,y)]_{y^-}^{y^+}
		=\frac{y_+(y)z_-(y)-z_+(y)y_-(y)}{\omega}=-1.
	\end{equation*}
	Thus $K_xG(x,y)=\delta_y$ in distributions,
	where $\delta_y$ is the unit point mass at $y$.

	The functions $e^{-\gamma x}Y_+(x)$
	and $e^{\gamma x}Y_-(x)$ are continuous and $L$-periodic.
	Their boundedness gives \eqref{eq:green_decay}.
	Consequently, \eqref{eq:bounded_green_representation}
	defines a bounded function with bounded weak derivative.
	Testing the kernel equation against a compactly supported smooth function
	and using Fubini's theorem gives $KV=q$.

	The difference of two bounded weak solutions solves $KW=0$.
	Its continuous state $(W,pW_x)$ therefore solves
	\eqref{eq:linear_tail_system} and is a linear combination of $Y_+$ and $Y_-$.
	Boundedness at $+\infty$ excludes the $Y_+$ term,
	and boundedness at $-\infty$ excludes the $Y_-$ term.
	This proves uniqueness.
	The comparison assertions follow from
	\begin{equation*}
		V_2(x)-V_1(x)
		=\int_{\mathbb R}G(x,y)\bigl(q_2(y)-q_1(y)\bigr)\,dy
	\end{equation*}
	and strict positivity of $G$.
\end{proof}

For the half-line forcing $\mathbf 1_{\{x<\rho\}}$,
the value of the bounded solution at $\rho$ is the matching function.

\begin{lemma}[Profiles with one prescribed interface]
	\label{lem:step_profiles}
	For $\rho\in\mathbb R$, let $V_\rho$ be the unique bounded weak solution of
	\begin{equation}
		\label{eq:step_profile_problem}
		KV_\rho=\mathbf 1_{\{x<\rho\}}.
	\end{equation}
	Then $0<V_\rho<1$, and $V_\rho$ is strictly decreasing,
	with limits $1$ and $0$ at $-\infty$ and $+\infty$, respectively.
	Its value at the interface is
	\begin{equation}
		\label{eq:step_profile_phase}
		V_\rho(\rho)=b(\rho),
	\end{equation}
	where $b$ is the matching function \eqref{eq:matching_function}.
\end{lemma}

\begin{proof}
	Lemma~\ref{lem:positive_bounded_inverse} gives
	\begin{equation}
		\label{eq:step_green_representation}
		V_\rho(x)=\int_{-\infty}^{\rho}G(x,y)\,dy.
	\end{equation}
	The bounded solutions for the constant forcings $0$ and $1$ are $0$ and $1$.
	Strict comparison therefore gives $0<V_\rho<1$.
	On the left of $\rho$, $V_\rho-1$ solves the homogeneous equation;
	on the right, $V_\rho$ does.
	As in the uniqueness part of the proof of Lemma~\ref{lem:positive_bounded_inverse}, the state of each is a linear combination of $Y_+$ and $Y_-$ on its half-line.
	By \eqref{eq:tail_floquet_relations}, $y_+$ decays and $y_-$ grows exponentially as $x\to-\infty$, and the reverse holds as $x\to+\infty$.
	Boundedness therefore excludes the $y_-$ term on the left of $\rho$ and the $y_+$ term on the right.
	The value $V_\rho(\rho)$ fixes the remaining coefficients, which gives
	\begin{equation}
		\label{eq:step_profile_tails}
		V_\rho(x)=
		\begin{cases}
			1-\bigl(1-V_\rho(\rho)\bigr)y_+(x)/y_+(\rho),&x\leq\rho,\\
			V_\rho(\rho)\,y_-(x)/y_-(\rho),&x\geq\rho.
		\end{cases}
	\end{equation}
	The signs $z_-<0<z_+$ prove strict decrease,
	and the period identities \eqref{eq:tail_floquet_relations} give the limits at infinity.
	By \eqref{eq:floquet_slopes}, the fluxes of the two pieces at $\rho$ are
	$-\bigl(1-V_\rho(\rho)\bigr)r_+(\rho)$ and $-V_\rho(\rho)r_-(\rho)$.
	Continuity of the flux at $\rho$ gives
	$\bigl(1-V_\rho(\rho)\bigr)r_+(\rho)=V_\rho(\rho)r_-(\rho)$,
	which is \eqref{eq:step_profile_phase}.
\end{proof}

By \eqref{eq:step_profile_phase},
the tail formulas \eqref{eq:step_profile_tails} hold with $V_\rho(\rho)$ replaced by $b(\rho)$.
Lemma~\ref{lem:matching_function} gives the range and monotonicity of $b$:
it decreases from $a_+$ at $\rho=0$ to $a_-$ at $\rho=L_1$
and increases back on $(L_1,L)$.

The profile $V_\rho$ is defined for every $\rho$,
and it solves the nonlinear McKean equation at threshold $a$
precisely when $b(\rho)=a$.

\subsection{Local geometry of the invariant-manifold intersections}
\label{subsec:local_intersections}

The comparison arguments of Section~\ref{sec:stability}
require the ordering of the invariant branches near a front.
Fix $(\alpha,L_1,L_2,a)$ with $a\in[a_-,a_+]$,
and let $s_0$ be the crossing phase of a front from
Theorem~\ref{thm:monotone_classification}.
The next two lemmas write the unstable and stable branches through its threshold state
as graphs over the state coordinate,
and relate their ordering to the derivative of the matching function at $s_0$. Their proofs are in Appendix~\ref{app:mckean_geometry}.

\begin{lemma}[Local ordering of the invariant branches]
	\label{lem:local_branch_ordering}
	Fix $a_-<a<a_+$ and $s_0\in\{\rho_P,\rho_Q\}$,
	and set
	\begin{equation*}
		v_*:=v^{\mathrm u}(s_0,a)=v^{\mathrm s}(s_0,a)<0,
		\qquad
		\kappa:=p(s_0)^{-1}.
	\end{equation*}
	The unstable manifold of $(1,0)$
	and the stable manifold of $(0,0)$ for $\Psi^{s_0}$
	admit local arcs through $(a,v_*)$
	that are graphs $v=h^u(\xi)$ and $v=h^s(\xi)$, respectively.
	Both graphs are continuous and of class $C^1$ on either side of $\xi=a$.
	Their difference is continuously differentiable across $a$ and satisfies
	\begin{equation}
		\label{eq:local_branch_gap}
		(h^s-h^u)'(a)
		=\frac{\partial_s(v^{\mathrm u}-v^{\mathrm s})(s_0,a)}{\kappa v_*}
		=-\frac{\bigl(r_+(s_0)+r_-(s_0)\bigr)b'(s_0)}{\kappa v_*}.
	\end{equation}
	Consequently,
	for all $\xi\neq a$ sufficiently close to $a$,
	\begin{equation}
		\label{eq:local_branch_ordering}
		(\xi-a)\bigl(h^s(\xi)-h^u(\xi)\bigr)
		\begin{cases}
			<0, & s_0=\rho_P, \\
			>0, & s_0=\rho_Q.
		\end{cases}
	\end{equation}
	After restricting these arcs if necessary,
	one forward iterate maps the stable arc into the straight portion
	of the stable manifold described in Lemma~\ref{lem:equilibrium_tails},
	and one backward iterate maps the unstable arc
	into the corresponding straight portion of the unstable manifold.
\end{lemma}

When $a=a_-$ or $a=a_+$, the crossing lies at a coefficient interface,
so the local ordering depends on the two one-sided derivatives.

\begin{lemma}[Intersections at the existence boundary]
	\label{lem:boundary_contact}
	Let $(a,s_0)=(a_-,L_1)$ or $(a_+,0)$,
	and set $v_*:=v^{\mathrm u}(s_0,a)=v^{\mathrm s}(s_0,a)<0$.
	The unstable manifold of $(1,0)$
	and the stable manifold of $(0,0)$ for $\Psi^{s_0}$
	admit local arcs through $(a,v_*)$
	that are continuous graphs $v=h^u(\xi)$ and $v=h^s(\xi)$, respectively.
	They satisfy $h^u(a)=h^s(a)=v_*$
	and are of class $C^1$ on either side of $\xi=a$.
	For all $\xi\neq a$ sufficiently close to $a$,
	\begin{equation}
		\label{eq:boundary_branch_ordering}
		h^s(\xi)-h^u(\xi)
		\begin{cases}
			>0, & a=a_-, \\
			<0, & a=a_+.
		\end{cases}
	\end{equation}
	The difference has nonzero one-sided derivatives of opposite signs at $a$.
\end{lemma}

Thus the two invariant branches meet without crossing at an existence endpoint,
and the graph of $h^s-h^u$ has a corner at $\xi=a$.
This contact between the two invariant-manifold arcs
is a nonsmooth analogue of the heteroclinic tangency
described in Section~\ref{sec:switchingODE}.
For the corresponding front $u$ with $u(s_0)=a$,
continuity of the flux and $v_*<0$ give
\begin{equation*}
	u_x(s_0^-)=\frac{v_*}{p(s_0^-)}<0,
	\qquad
	u_x(s_0^+)=\frac{v_*}{p(s_0^+)}<0.
\end{equation*}
Thus the front still crosses $u=a$ with nonzero one-sided derivatives
when $a=a_-$ or $a=a_+$.

\section{Stability of the decreasing fronts}\label{sec:stability}

Fix $0<\alpha<1$, $L_1,L_2>0$, and $a\in(a_-,a_+)$,
where $a_\pm$ are defined in \eqref{eq:front_existence_thresholds}.
For the McKean evolution equation \eqref{pRD},
we prove Lyapunov stability of $u^P$ and instability of $u^Q$
in $W^{1,\infty}(\mathbb{R})$.
These are the two fronts defined after Theorem~\ref{thm:monotone_classification},
with threshold crossings $\rho_P\in(0,L_1)$ and $\rho_Q\in(L_1,L)$.
Periodicity gives the same conclusions for their translates by integer multiples of $L$.

We work with mild solutions near a decreasing stationary front,
with one threshold crossing whose location varies continuously in time.
Prescribing this crossing path makes the McKean equation linear.
We obtain a solution by requiring the prescribed path
to coincide with the threshold crossing of the resulting profile.
After stating the local existence and comparison results,
we construct stationary sub- and supersolutions.
Near $u^P$ they are vertical shifts of the profiles $V_\rho$,
whose signs come from the monotonicity of the matching function.
Near $u^Q$ they are joined from pieces of the invariant manifolds,
whose local ordering is given in Section~\ref{subsec:local_intersections}.
They trap solutions near $u^P$ and force suitable perturbations away from $u^Q$.
The semigroup, single-crossing, and comparison proofs are given
in Appendix~\ref{app:temporal}.
The parameters remain fixed throughout this section,
and no uniform neighborhood size is asserted as $a\to a_\pm$.

\subsection{Local existence and comparison}

Fix a decreasing stationary front $u_*$ of \eqref{pblm} with \eqref{BC},
and let $\rho_*$ be its unique threshold-crossing point.
Write $u_0=u_*+w_0$ with $w_0\in W^{1,\infty}(\mathbb R)$.
For $T>0$ and a \textit{prescribed path} $\rho\in C([0,T])$,
consider the linear problem
\begin{equation}
	\label{eq:trial_RD}
	\partial_tu^\rho=(p\partial_xu^\rho)_x-u^\rho
	+\mathbf{1}_{\{x<\rho(t)\}},
	\qquad u^\rho(\cdot,0)=u_0.
\end{equation}
Let $\mathsf S(t)$ be the semigroup in Lemma~\ref{lem:semigroup_est}.
Since $u_*$ satisfies
$(p\partial_xu_*)_x-u_*+\mathbf{1}_{\{x<\rho_*\}}=0$ weakly,
the mild formulation for $w(t):=u^\rho(\cdot,t)-u_*$ is
\begin{equation}\label{eq:mild_w}
	\begin{aligned}
		w(t)
		 & =\mathsf S(t)w_0+\int_0^t \mathsf S(t-r)h_\rho(r)\,dr,                  \\
		h_\rho(r)
		 & :=\mathbf{1}_{\{x<\rho(r)\}}-\mathbf{1}_{\{x<\rho_*\}}.
	\end{aligned}
\end{equation}
Here $\mathsf S(t)w_0$ uses the $L^\infty(\mathbb{R})$ extension of the semigroup,
and the time integral is a Bochner integral in $L^2(\mathbb R)$,
since $h_\rho\in C([0,T];L^2(\mathbb R))$.
The prescribed path specifies the discontinuity of the forcing.
The equation agrees with \eqref{pRD} when $u^\rho(\cdot,t)>a$
to the left of $\rho(t)$ and $u^\rho(\cdot,t)<a$ to its right.

\begin{lemma}\label{lem:semigroup_est}
	The $L^2(\mathbb R)$ realization of $K=-\partial_x(p\partial_x)+1$ has domain
	\begin{equation*}
		\mathcal D(K)
		:=\{g\in H^1(\mathbb R):p g_x\in H^1(\mathbb R)\},
	\end{equation*}
	and $-K$ generates the strongly continuous semigroup $\mathsf S_2(t)=e^{-tK}$ with
	$\|\mathsf S_2(t)\|_{2\to2}\leq e^{-t}$ for $t\geq0$.
	There is a strongly continuous contraction semigroup $\mathsf S_1(t)$ on $L^1(\mathbb R)$
	that agrees with $\mathsf S_2(t)$ on $L^1(\mathbb{R})\cap L^2(\mathbb{R})$.
	Its Banach adjoint $\mathsf S_\infty(t):=\mathsf S_1(t)^*$
	is a contraction semigroup on $L^\infty(\mathbb R)$
	that agrees with $\mathsf S_2(t)$ on $L^\infty(\mathbb{R})\cap L^2(\mathbb{R})$.
	For each $g\in L^\infty(\mathbb R)$,
	the map $t\mapsto \mathsf S_\infty(t)g$ is continuous
	in the weak-$*$ topology $\sigma(L^\infty(\mathbb{R}),L^1(\mathbb{R}))$.
	We denote these consistent semigroups by $\mathsf S(t)$.
	There are constants $C_0,C_1>0$, depending only on $p$,
	such that the following estimates hold for every $t>0$:
	\begin{align}
		\|\mathsf S(t)g\|_{\infty}
		 & \leq e^{-t}\|g\|_{\infty},
		 &                                & g\in L^\infty(\mathbb{R}),
		\label{eq:S_contract}                                              \\
		\|p\partial_x\mathsf S(t)g\|_{\infty}
		 & \leq e^{-t}\|pg_x\|_{\infty},
		 &                                & g\in W^{1,\infty}(\mathbb{R}),
		\label{eq:DxS_contract}                                            \\
		\|\mathsf S(t)g\|_{\infty}
		 & \leq C_0e^{-t}t^{-1/2}\|g\|_1,
		 &                                & g\in L^1(\mathbb{R}),
		\label{eq:S_bdd}                                                   \\
		\|\partial_x\mathsf S(t)g\|_{\infty}
		 & \leq C_1e^{-t}t^{-3/4}\|g\|_2,
		 &                                & g\in L^2(\mathbb{R}).
		\label{eq:DxS_bdd}
	\end{align}
	For the preceding mild formulation \eqref{eq:mild_w}, define
	\begin{equation*}
		d_\rho(t):=\max_{0\leq r\leq t}|\rho(r)-\rho_*|,
		\qquad
		I_\gamma(t):=\int_0^t e^{-r}r^{-\gamma}\,dr,
		\quad 0<\gamma<1.
	\end{equation*}
	Then, for $0\leq t\leq T$,
	\begin{align}
		\|w(t)\|_\infty
		 & \leq e^{-t}\|w_0\|_\infty
		+C_0d_\rho(t)I_{1/2}(t),
		\label{eq:w_bdd}                                         \\
		\|w_x(t)\|_\infty
		 & \leq e^{-t}\|p^{-1}\|_\infty\|p\partial_xw_0\|_\infty
		+C_1d_\rho(t)^{1/2}I_{3/4}(t).
		\label{eq:wx_bdd}
	\end{align}
\end{lemma}

\begin{lemma}[Single threshold crossing]
	\label{lem:single_crossing}
	Let $R,\delta,M>0$ satisfy
	\begin{equation}
		\label{eq:R_restriction}
		\begin{aligned}
			u_*(x)   & \geq a+2\delta &                                          & (x\leq\rho_*-R), \\
			u_*(x)   & \leq a-2\delta &                                          & (x\geq\rho_*+R), \\
			\partial_xu_*(x) & \leq-M
			       &                & \text{for a.e. }x\in[\rho_*-R,\rho_*+R].
		\end{aligned}
	\end{equation}
	Let $C_0,C_1$ be as in Lemma~\ref{lem:semigroup_est}.
	Suppose $T>0$ and $0<M_*<M$ satisfy
	\begin{equation}
		\label{eq:w}
		\|w_0\|_\infty+C_0R I_{1/2}(T)<\delta
	\end{equation}
	and
	\begin{equation}
		\label{eq:wm}
		\|p^{-1}\|_\infty\|p\partial_xw_0\|_\infty
		+C_1R^{1/2}I_{3/4}(T)\leq M_*.
	\end{equation}
	Define
	\begin{equation*}
		X_{R,T}:=\{\rho\in C([0,T]):
		\|\rho-\rho_*\|_{C([0,T])}\leq R\}.
	\end{equation*}
	For every $\rho\in X_{R,T}$ and $t\in[0,T]$,
	the mild solution $u^\rho$ of \eqref{eq:trial_RD}
	has a unique threshold-crossing point $\widetilde\rho(t)$.
	These points define a continuous path in $X_{R,T}$.
	The map
	\begin{equation*}
		\Phi_{u_0}:X_{R,T}\to X_{R,T},
		\qquad \Phi_{u_0}(\rho):=\widetilde\rho,
	\end{equation*}
	is a contraction if
	\begin{equation}
		\label{eq:contract_est}
		\frac{C_0 I_{1/2}(T)}{M-M_*}<1.
	\end{equation}
\end{lemma}

\begin{proposition}[Local existence and uniqueness]
	\label{prop:loc_exist_unique}
	Let $R,\delta,M$ satisfy \eqref{eq:R_restriction}, and suppose
	\begin{equation*}
		\|w_0\|_\infty<\delta,
		\qquad
		\|p^{-1}\|_\infty\|p\partial_xw_0\|_\infty<M.
	\end{equation*}
	There exist $T>0$, a mild solution $u$ of \eqref{pRD}
	with initial datum $u_0=u_*+w_0$, and $\rho\in X_{R,T}$ such that
	\begin{equation*}
		u(\rho(t),t)=a,
		\qquad
		u(x,t)>a\quad(x<\rho(t)),
		\qquad
		u(x,t)<a\quad(x>\rho(t))
	\end{equation*}
	for every $t\in[0,T]$.
	The solution is unique among mild solutions with the same initial datum
	and a crossing path in $X_{R,T}$ satisfying these sign conditions.
\end{proposition}

\begin{proof}
	Choose $M_*$ strictly between
	$\|p^{-1}\|_\infty\|p\partial_xw_0\|_\infty$ and $M$.
	Since $I_{1/2}(T)$ and $I_{3/4}(T)$ tend to zero as $T\downarrow0$,
	conditions \eqref{eq:w}, \eqref{eq:wm}, and \eqref{eq:contract_est}
	hold for sufficiently small $T>0$.
	The space $X_{R,T}$ is closed in $C([0,T])$ and hence complete.
	Lemma~\ref{lem:single_crossing} and the contraction mapping theorem
	give a unique fixed point $\rho$ of $\Phi_{u_0}$.
	For $u:=u^\rho$, the prescribed path is its actual threshold-crossing path,
	so \eqref{eq:trial_RD} agrees with \eqref{pRD} almost everywhere in space.
	Any other mild solution with the stated sign properties
	satisfies \eqref{eq:mild_w} with its own crossing path.
	That path is another fixed point of $\Phi_{u_0}$,
	so it coincides with $\rho$ and the solutions agree.
\end{proof}

The comparison principle uses positivity of the semigroup.

\begin{lemma}[Positivity preservation]
	\label{lem:semigroup_positive}
	For $q\in\{1,2,\infty\}$ and $t\geq0$,
	the semigroup $\mathsf S(t)$ preserves nonnegativity on $L^q(\mathbb{R})$.
	Consequently, if $f,g\in L^q(\mathbb{R})$ and $f\leq g$ almost everywhere,
	then $\mathsf S(t)f\leq \mathsf S(t)g$ almost everywhere.
\end{lemma}

For $\rho\in X_{R,T}$ and $f-u_*\in W^{1,\infty}(\mathbb{R})$,
denote the solution of \eqref{eq:trial_RD} with initial datum $f$ by
\begin{equation}
	\label{eq:prescribed_duhamel}
	\mathcal{M}[\rho,f](t)
	:=u_*+\mathsf S(t)(f-u_*)
	+\int_0^t \mathsf S(t-r)
	\bigl(\mathbf{1}_{\{x<\rho(r)\}}-\mathbf{1}_{\{x<\rho_*\}}\bigr)\,dr.
\end{equation}

\begin{proposition}[Local comparison principle]
	\label{prop:loc_comparison}
	Let $R,\delta,M>0$ satisfy \eqref{eq:R_restriction},
	and fix $T>0$ and $0<M_*<M$.
	Let $u^-,u^+$ be bounded functions on $\mathbb{R}\times[0,T]$,
	with $u^\pm(\cdot,t)$ continuous for every $t\in[0,T]$.
	Write $u_0^\pm:=u^\pm(\cdot,0)$ and $w_0^\pm:=u_0^\pm-u_*$,
	and suppose $w_0^\pm\in W^{1,\infty}(\mathbb{R})$.
	Assume \eqref{eq:w} and \eqref{eq:wm} hold
	with $w_0=w_0^-$ and with $w_0=w_0^+$,
	and that \eqref{eq:contract_est} holds.
	Suppose there are paths $\rho^\pm\in X_{R,T}$ such that,
	for every $t\in[0,T]$,
	\begin{equation*}
		\begin{aligned}
			u^\pm(\rho^\pm(t),t) & =a,                       \\
			u^\pm(x,t)           & >a  &  & (x<\rho^\pm(t)), \\
			u^\pm(x,t)           & <a  &  & (x>\rho^\pm(t)).
		\end{aligned}
	\end{equation*}
	Assume also the mild sub- and supersolution inequalities for \eqref{pRD},
	\begin{equation*}
		u^-(\cdot,t)\leq\mathcal{M}[\rho^-,u_0^-](t),
		\qquad
		\mathcal{M}[\rho^+,u_0^+](t)\leq u^+(\cdot,t),
	\end{equation*}
	for every $t\in[0,T]$.
	If $u_0^-\leq u_0^+$, then $u^-(x,t)\leq u^+(x,t)$
	for all $(x,t)\in\mathbb{R}\times[0,T]$.
\end{proposition}

\subsection{Lyapunov stability}
\label{subsec:stability}

In the local theory, take $u_*=u^P$ and $\rho_*=\rho_P$.
We construct nearby stationary sub- and supersolutions
by vertically shifting the profiles $V_\rho$ from Lemma~\ref{lem:step_profiles}.
Recall that $V_\rho$ is the bounded solution of
$KV_\rho=\mathbf{1}_{\{x<\rho\}}$ and that $V_{\rho_P}=u^P$.
To estimate $V_\rho-u^P$ in $W^{1,\infty}(\mathbb R)$
and verify the mild comparison inequalities after shifting,
we first represent this difference by the evolution semigroup.
The weak stationary equation for $u^P$ is
\begin{equation}
	Ku^P=\mathbf{1}_{\{x<\rho_P\}}.
\end{equation}
For $\rho$ near $\rho_P$, set
\begin{equation*}
	h_\rho:=\mathbf{1}_{\{x<\rho\}}-\mathbf{1}_{\{x<\rho_P\}},
	\qquad
	z_\rho:=\int_0^\infty \mathsf S(t)h_\rho\,dt.
\end{equation*}
Since $h_\rho\in L^1(\mathbb{R})\cap L^2(\mathbb{R})$
and $\|\mathsf S_2(t)\|_{2\to2}\leq e^{-t}$,
the integral converges in $L^2(\mathbb{R})$.
The resolvent formula gives
$z_\rho=K^{-1}h_\rho\in\mathcal D(K)$ and $Kz_\rho=h_\rho$.
The semigroup estimates show that $u^P+z_\rho$ is bounded.
Combining the equations for $u^P$ and $z_\rho$ gives
\begin{equation}
	\label{eq:prescribed_interface_eq}
	K(u^P+z_\rho)=\mathbf{1}_{\{x<\rho\}}
\end{equation}
in the weak sense.
Uniqueness in Lemma~\ref{lem:positive_bounded_inverse}
therefore gives the representation
\begin{equation}
	\label{eq:prescribed_interface}
	V_\rho=u^P+z_\rho
\end{equation}
for the same prescribed-interface profile $V_\rho$.
Its value at the prescribed interface need not equal $a$;
write
\begin{equation}
	\label{eq:g_rho}
	g(\rho):=b(\rho)-a=V_\rho(\rho)-a.
\end{equation}
The sign of $g$ determines the direction of the vertical shift
needed to make $V_\rho$ attain the threshold at $\rho$.
Since $b(\rho_P)=a$ and $b$ is strictly decreasing on $(0,L_1)$
(Lemma~\ref{lem:matching_function}),
for $\rho$ sufficiently close to $\rho_P$ we have
\begin{equation}
	\label{eq:g_sign}
	g(\rho)>0 \quad\text{if }\rho<\rho_P,
	\qquad
	g(\rho)<0 \quad\text{if }\rho>\rho_P.
\end{equation}

To make $\rho$ the threshold crossing, shift $V_\rho$ vertically:
\begin{equation}
	\label{eq:Brho}
	B_\rho:=V_\rho-g(\rho).
\end{equation}
Then
\begin{equation}\label{eq:Brho_threshold}
	B_\rho(\rho)=a.
\end{equation}
Since $g(\rho)$ is constant in $x$,
\eqref{eq:prescribed_interface_eq} gives
\begin{equation}
	\label{eq:Brho_prescribed}
	-KB_\rho+\mathbf{1}_{\{x<\rho\}}=g(\rho).
\end{equation}

\begin{lemma}\label{lem:Brho_close}
	As $\rho\to\rho_P$,
	\begin{equation}
		\|B_\rho-u^P\|_{W^{1,\infty}(\mathbb{R})}
		\leq
		C\left(
		|\rho-\rho_P|
		+
		|\rho-\rho_P|^{1/2}
		\right)
		\longrightarrow 0.
	\end{equation}
	In particular, for $\rho$ sufficiently close to $\rho_P$,
	$B_\rho$ has a unique threshold crossing, located at $x=\rho$.
\end{lemma}

\begin{proof}
	For $h_\rho$ defined above,
	\begin{equation*}
		\|h_\rho\|_1=|\rho-\rho_P|,
		\qquad
		\|h_\rho\|_2=|\rho-\rho_P|^{1/2}.
	\end{equation*}
	Lemma~\ref{lem:semigroup_est} therefore gives
	\begin{equation*}
		\|z_\rho\|_\infty
		\leq C_0\sqrt{\pi}\,|\rho-\rho_P|,
		\qquad
		\|\partial_xz_\rho\|_\infty
		\leq C_1\Gamma(1/4)|\rho-\rho_P|^{1/2}.
	\end{equation*}
	The derivative estimate follows by differentiating the integral
	in distributions, using the integrability of $e^{-t}t^{-3/4}$.
	Since $u^P(\rho_P)=a$,
	\begin{equation*}
		|g(\rho)|
		\leq |z_\rho(\rho)|+|u^P(\rho)-u^P(\rho_P)|
		\leq \bigl(C_0\sqrt{\pi}+\|u_x^P\|_\infty\bigr)|\rho-\rho_P|.
	\end{equation*}
	The stated estimate follows from $B_\rho-u^P=z_\rho-g(\rho)$.

	Lemma~\ref{lem:step_profiles} shows that $V_\rho$ is strictly decreasing.
	The vertical shift $B_\rho$ is therefore strictly decreasing as well.
	Since $B_\rho(\rho)=a$, it has the asserted unique threshold crossing,
	with $B_\rho(x)>a$ for $x<\rho$ and $B_\rho(x)<a$ for $x>\rho$.
\end{proof}

\begin{lemma}[Stationary comparison profiles]
	\label{lem:Brho_trap}
	Let $\rho_-<\rho_P<\rho_+$ be sufficiently close to $\rho_P$.
	Then:
	\begin{enumerate}
		\item $B_{\rho_-}$ and $B_{\rho_+}$ are respectively
		      a stationary mild subsolution and a stationary mild supersolution
		      of \eqref{pRD}.
		\item $B_{\rho_-}<u^P<B_{\rho_+}$ on $\mathbb{R}$.
		\item There exists $\eta>0$ such that $B_{\rho_-}<u<B_{\rho_+}$
		      for every bounded continuous $u$ with $\|u-u^P\|_\infty<\eta$.
	\end{enumerate}
\end{lemma}

\begin{proof}
	For a fixed $\rho$ near $\rho_P$,
	also denote the constant path $t\mapsto\rho$ by $\rho$.
	The semigroup property and convergence of the integral defining $z_\rho$
	give, in $L^2(\mathbb{R})$,
	\begin{equation*}
		\mathsf S(t)z_\rho=\int_t^\infty \mathsf S(r)h_\rho\,dr,
		\qquad
		\mathsf S(t)z_\rho+\int_0^t \mathsf S(t-r)h_\rho\,dr=z_\rho.
	\end{equation*}
	Let $\mathbf{1}$ denote the constant function equal to one.
	Use $u_*=u^P$ in the definition of $\mathcal M$ in \eqref{eq:prescribed_duhamel}.
	Since $B_\rho-u^P=z_\rho-g(\rho)\mathbf{1}$,
	consistency of the semigroup realizations gives
	\begin{align*}
		\mathcal M[\rho,B_\rho](t)
		 & =u^P+\mathsf S(t)z_\rho-g(\rho)\mathsf S_\infty(t)\mathbf{1}
		+\int_0^t \mathsf S(t-r)h_\rho\,dr                                      \\
		 & =B_\rho+g(\rho)\bigl(\mathbf{1}-\mathsf S_\infty(t)\mathbf{1}\bigr).
	\end{align*}
	Positivity and \eqref{eq:S_contract} imply
	\begin{equation*}
		0\leq \mathsf S_\infty(t)\mathbf{1}\leq e^{-t}\mathbf{1}.
	\end{equation*}
	Together with $g(\rho_-)>0$ and $g(\rho_+)<0$
	from \eqref{eq:g_sign}, this yields
	\begin{equation}
		\label{eq:Brho_mild}
		B_{\rho_-}\leq\mathcal M[\rho_-,B_{\rho_-}](t),
		\qquad
		\mathcal M[\rho_+,B_{\rho_+}](t)\leq B_{\rho_+},
		\qquad t\geq0.
	\end{equation}
	These inequalities hold for the continuous spatial representatives.
	Together with the crossing properties from Lemma~\ref{lem:Brho_close},
	they give the stationary mild sub- and supersolutions
	required by Proposition~\ref{prop:loc_comparison}.

	For $\rho>\rho_P$, positivity also gives
	\begin{equation}
		\label{eq:Brho_upper}
		B_\rho-u^P
		=\int_0^\infty \mathsf S(t)\mathbf{1}_{(\rho_P,\rho)}\,dt-g(\rho)
		\geq-g(\rho)>0.
	\end{equation}
	For $\rho<\rho_P$, similarly,
	\begin{equation}
		\label{eq:Brho_lower}
		u^P-B_\rho
		=\int_0^\infty \mathsf S(t)\mathbf{1}_{(\rho,\rho_P)}\,dt+g(\rho)
		\geq g(\rho)>0.
	\end{equation}
	The asserted ordering follows.
	Choosing $0<\eta<\min\{g(\rho_-),-g(\rho_+)\}$
	then gives the final assertion.
\end{proof}

\begin{theorem}[Lyapunov stability]
	\label{thm:stable_front_stability}
	The front $u^P$ is Lyapunov stable in $W^{1,\infty}(\mathbb{R})$.
	More precisely, there exists $R>0$ such that, for every $\varepsilon>0$,
	there is $\eta>0$ with the following property.
	Every initial datum $u_0$ satisfying
	\begin{equation*}
		\|u_0-u^P\|_{W^{1,\infty}}<\eta
	\end{equation*}
	admits a global mild solution $u$ of \eqref{pRD} with
	\begin{equation*}
		\|u(\cdot,t)-u^P\|_{W^{1,\infty}}<\varepsilon,
		\qquad t\geq0.
	\end{equation*}
	There is a continuous path
	$\rho:[0,\infty)\to[\rho_P-R,\rho_P+R]$ such that
	\begin{equation*}
		u(\rho(t),t)=a,
		\qquad u(x,t)>a\quad(x<\rho(t)),
		\qquad u(x,t)<a\quad(x>\rho(t)).
	\end{equation*}
	The solution is unique among mild solutions with the same initial datum
	and a continuous crossing path in this interval satisfying these sign conditions.
\end{theorem}

\begin{proof}
	The comparison profiles bound both the solution and its crossing path.
	The semigroup estimates then control its derivative,
	allowing the local construction to continue on intervals of a fixed length.

	Fix $R,\delta,M>0$ satisfying \eqref{eq:R_restriction} with $u_*=u^P$,
	and fix $0<M_*<M$.
	Set $C_p:=\|p^{-1}\|_\infty\|p\|_\infty$.
	Given $\varepsilon>0$, for $0<r<R$ set
	\begin{equation*}
		\rho_-:=\rho_P-r,
		\qquad
		\rho_+:=\rho_P+r.
	\end{equation*}
	For all sufficiently small $r$,
	Lemma~\ref{lem:Brho_trap} shows that $B_{\rho_-}$ is a stationary mild subsolution
	and $B_{\rho_+}$ is a stationary mild supersolution, with
	$B_{\rho_-}<u^P<B_{\rho_+}$.
	Define
	\begin{equation*}
		D_r:=\max_{\nu\in\{-,+\}}\|B_{\rho_\nu}-u^P\|_\infty,
		\qquad
		A_r:=C_pC_1\Gamma(1/4)r^{1/2}.
	\end{equation*}
	Since $\partial_x(B_\rho-u^P)=\partial_xz_\rho$,
	the derivative estimate in the proof of Lemma~\ref{lem:Brho_close} gives
	\begin{equation*}
		\|p^{-1}\|_\infty
		\|p\partial_x(B_{\rho_\pm}-u^P)\|_\infty\leq A_r.
	\end{equation*}
	Moreover, $D_r,A_r\to0$ as $r\downarrow0$.
	We may therefore fix $r$ sufficiently small that
	\begin{equation}
		\label{eq:stable_r_choice}
		D_r<\min\left\{\frac{\delta}{2},\frac{\varepsilon}{2}\right\},
		\qquad
		A_r<\min\left\{\frac{M_*}{4},\frac{\varepsilon}{4}\right\}.
	\end{equation}

	By Lemma~\ref{lem:Brho_trap}, there is $\eta_0>0$ such that
	$B_{\rho_-}<f<B_{\rho_+}$ whenever $f$ is bounded and continuous
	and $\|f-u^P\|_\infty<\eta_0$.
	Choose
	\begin{equation}
		\label{eq:stable_eta_choice}
		0<\eta<\min\left\{\eta_0,\frac{M_*}{4C_p},
		\frac{\varepsilon}{4C_p}\right\}.
	\end{equation}
	Since $I_{1/2}(T),I_{3/4}(T)\to0$ as $T\downarrow0$,
	we may also fix $T>0$ such that
	\begin{equation}
		\label{eq:stable_T_choice}
		C_0R I_{1/2}(T)<\frac{\delta}{2},
		\qquad
		C_1R^{1/2}I_{3/4}(T)<\frac{M_*}{2},
		\qquad
		\frac{C_0 I_{1/2}(T)}{M-M_*}<1.
	\end{equation}

	Let $u_0$ satisfy $\|u_0-u^P\|_{W^{1,\infty}}<\eta$,
	and put $w_0:=u_0-u^P$.
	Then $B_{\rho_-}<u_0<B_{\rho_+}$ and
	\begin{equation*}
		\|w_0\|_\infty\leq D_r,
		\qquad
		\|p^{-1}\|_\infty\|p\partial_xw_0\|_\infty
		\leq C_p\eta<\frac{M_*}{4}.
	\end{equation*}
	These estimates and \eqref{eq:stable_T_choice}
	verify \eqref{eq:w}, \eqref{eq:wm}, and \eqref{eq:contract_est}
	for the initial datum $u_0$.
	For the initial data $B_{\rho_-}$ and $B_{\rho_+}$,
	the initial-data terms in \eqref{eq:w} and \eqref{eq:wm}
	are bounded by $D_r$ and $A_r$, respectively,
	so the same conditions hold.
	Their constant crossing paths $t\mapsto\rho_\pm$ belong to $X_{R,T}$.
	The construction in Proposition~\ref{prop:loc_exist_unique}
	therefore gives a mild solution $u$ on $[0,T]$
	with a unique continuous threshold-crossing path $\rho\in X_{R,T}$.
	Applying Proposition~\ref{prop:loc_comparison}
	to compare $u$ with each stationary profile yields
	\begin{equation}
		\label{eq:stable_trapping}
		B_{\rho_-}\leq u(\cdot,t)\leq B_{\rho_+},
		\qquad 0\leq t\leq T.
	\end{equation}
	In particular, $u(\rho_-,t)\geq a$ and $u(\rho_+,t)\leq a$.
	The sign conditions on either side of the unique crossing imply
	\begin{equation}
		\label{eq:stable_crossing}
		\rho_-\leq\rho(t)\leq\rho_+,
		\qquad |\rho(t)-\rho_P|\leq r.
	\end{equation}

	Write $w(t):=u(\cdot,t)-u^P$.
	To restart the local construction at time $T$,
	we also need the derivative bound in \eqref{eq:wm}.
	In \eqref{eq:mild_w}, $h_\rho(s)$ is, up to sign, the indicator of the interval between $\rho(s)$ and $\rho_P$.
	The crossing bound \eqref{eq:stable_crossing} therefore gives
	\begin{equation*}
		\|h_\rho(s)\|_2=|\rho(s)-\rho_P|^{1/2}\leq r^{1/2},
		\qquad 0\leq s\leq T.
	\end{equation*}
	Apply \eqref{eq:DxS_contract} to the initial-data term
	and \eqref{eq:DxS_bdd} to the time integral.
	Since $I_{3/4}(t)\leq\Gamma(1/4)$, this gives
	\begin{align}
		\|p^{-1}\|_\infty\|p w_x(t)\|_\infty
		 & \leq e^{-t}\|p^{-1}\|_\infty\|p\partial_xw_0\|_\infty
		+C_pC_1r^{1/2}I_{3/4}(t)\notag                           \\
		 & \leq C_p\eta+A_r
		<\min\left\{\frac{M_*}{2},\frac{\varepsilon}{2}\right\},
		\qquad 0\leq t\leq T.
		\label{eq:stable_flux_bound}
	\end{align}

	We now extend the solution inductively.
	Suppose that, for some integer $n\geq1$, it is defined on $[0,nT]$,
	satisfies \eqref{eq:stable_trapping} and \eqref{eq:stable_crossing} there,
	and has the mild representation \eqref{eq:mild_w} with initial perturbation $w_0$.
	The preceding derivative estimate then holds throughout $[0,nT]$.
	At $t=nT$, we consequently have
	\begin{equation*}
		\|w(nT)\|_\infty\leq D_r<\frac{\delta}{2},
		\qquad
		\|p^{-1}\|_\infty\|p w_x(nT)\|_\infty<\frac{M_*}{2}.
	\end{equation*}
	Together with \eqref{eq:stable_T_choice}, these inequalities verify
	\eqref{eq:w}, \eqref{eq:wm}, and \eqref{eq:contract_est}
	with $w_0$ replaced by $w(nT)$.
	The stationary profiles $B_{\rho_-}$ and $B_{\rho_+}$
	continue to satisfy the same initial-data conditions,
	and \eqref{eq:stable_trapping} gives
	$B_{\rho_-}\leq u(\cdot,nT)\leq B_{\rho_+}$.
	Local existence and comparison therefore extend $u$ to $[0,(n+1)T]$,
	preserving the pointwise bounds and the crossing bound.
	The two crossing paths agree at $nT$, since they are crossings
	of the same spatial profile there, so the extended path is continuous.
	The semigroup property combines the mild formulas on the two intervals
	into \eqref{eq:mild_w} with the original initial datum $w_0$.
	Thus the induction hypotheses hold on $[0,(n+1)T]$.
	Since $T$ is fixed, this constructs the solution for all $t\geq0$.
	Applying local uniqueness on successive intervals of length $T$
	gives uniqueness in the single-crossing class stated in the theorem.

	Finally, \eqref{eq:stable_trapping} and \eqref{eq:stable_flux_bound}
	now hold for all $t\geq0$.
	They give
	\begin{equation*}
		\|w(t)\|_\infty\leq D_r<\frac{\varepsilon}{2},
		\qquad
		\|w_x(t)\|_\infty
		\leq\|p^{-1}\|_\infty\|p w_x(t)\|_\infty
		<\frac{\varepsilon}{2}.
	\end{equation*}
	Hence $\|u(\cdot,t)-u^P\|_{W^{1,\infty}}<\varepsilon$
	for every $t\geq0$.
\end{proof}

\subsection{Lyapunov instability}
\label{subsec:instability}

The vertical shifts used for $u^P$ do not give comparison profiles near $u^Q$.
For $\rho>\rho_Q$ close to $\rho_Q$, we have $g(\rho)=b(\rho)-a>0$,
because $b$ increases on $(L_1,L)$.
The shifted profile $V_\rho-g(\rho)$ equals $a>u^Q(\rho)$ at $x=\rho$,
whereas $V_\rho-g(\rho)-u^Q\to-g(\rho)<0$ as $x\to\pm\infty$.
Hence the two profiles are not ordered.
The case $\rho<\rho_Q$ is symmetric.
Instead, we join half-line pieces of the two invariant manifolds at $\rho_Q$.
The local ordering in Lemma~\ref{lem:local_branch_ordering}
gives the sign of the resulting flux jump,
which yields a stationary mild subsolution above $u^Q$
and a stationary mild supersolution below it.
A solution starting at the subsolution is nondecreasing in time
for as long as local comparison applies,
and we use this to force it out of a fixed neighborhood of $u^Q$.

\begin{lemma}[Comparison profiles for instability]
	\label{lem:unstable_subsuper}
	For every sufficiently small $\sigma>0$, there are a stationary mild
	subsolution $\underline u_\sigma$ and a stationary mild supersolution
	$\overline u_\sigma$ of \eqref{pRD}, both strictly decreasing, such that
	\begin{equation*}
		\overline u_\sigma<u^Q<\underline u_\sigma
		\quad\text{on }\mathbb{R}.
	\end{equation*}
	Their unique threshold crossings $\rho_\sigma^+$ and $\rho_\sigma^-$,
	respectively, satisfy
	\begin{equation*}
		\rho_\sigma^-<\rho_Q<\rho_\sigma^+,
		\qquad \rho_\sigma^\pm\to\rho_Q\quad(\sigma\downarrow0).
	\end{equation*}
	Moreover, $\underline u_\sigma-u^Q$ and $\overline u_\sigma-u^Q$ belong to $H^1(\mathbb R)$, and
	\begin{equation*}
		\|\underline u_\sigma-u^Q\|_{W^{1,\infty}}
		+\|\overline u_\sigma-u^Q\|_{W^{1,\infty}}
		\longrightarrow0\qquad(\sigma\downarrow0).
	\end{equation*}
\end{lemma}

Each profile joins, at $x=\rho_Q$, the unstable and stable half-line solutions
through the points of the local branches of Lemma~\ref{lem:local_branch_ordering}
with first coordinate $a\pm\sigma$.
Their only flux jump is at $\rho_Q$, and its sign,
given by \eqref{eq:local_branch_ordering},
makes the profiles a subsolution and a supersolution.
The proof is given in Appendix~\ref{app:temporal}.

\begin{theorem}[Lyapunov instability]
	\label{thm:unstable_front_instability}
	There exist $R>0$ and $\varepsilon_*>0$ such that, for every $\eta>0$,
	there are an initial datum $u_0$, a time $t_\eta>0$,
	and a mild solution $u$ of \eqref{pRD} on $[0,t_\eta]$
	with $u(\cdot,0)=u_0$ and
	\begin{equation*}
		\|u_0-u^Q\|_{W^{1,\infty}}<\eta,
		\qquad
		\|u(\cdot,t_\eta)-u^Q\|_{W^{1,\infty}}\geq\varepsilon_*.
	\end{equation*}
	The solution has a continuous crossing path
	$\rho:[0,t_\eta]\to[\rho_Q-R,\rho_Q+R]$ with
	\begin{equation*}
		u(\rho(t),t)=a,
		\qquad u(x,t)>a\quad(x<\rho(t)),
		\qquad u(x,t)<a\quad(x>\rho(t)).
	\end{equation*}
	It is unique among mild solutions with the same initial datum
	and a continuous crossing path in this interval satisfying these sign conditions.
	Thus $u^Q$ is Lyapunov unstable in $W^{1,\infty}(\mathbb{R})$
	within this single-crossing class.
\end{theorem}

\begin{proof}
	We choose a neighborhood in which $u^Q$ is the only stationary profile.
	A solution starting at $\underline u_\sigma>u^Q$ is nondecreasing in time
	for as long as it remains in this neighborhood.
	If it remained there for all time, we would obtain a stationary limit
	strictly above $u^Q$, contradicting this isolation.

	Choose $R,\delta,M>0$ satisfying \eqref{eq:R_restriction} with $u_*=u^Q$,
	with $R$ small enough that $J:=[\rho_Q-R,\rho_Q+R]$
	contains no threshold crossing of the decreasing fronts of Theorem~\ref{thm:monotone_classification}
	or their integer-period translates other than $\rho_Q$.
	Fix $0<M_*<M$, set $C_p:=\|p^{-1}\|_\infty\|p\|_\infty$, and choose
	\begin{equation*}
		0<\varepsilon<\min\left\{\frac{\delta}{2},\frac{M_*}{2C_p}\right\},
		\qquad \varepsilon_*:=\frac{\varepsilon}{2}.
	\end{equation*}
	Define the closed neighborhood
	\begin{equation*}
		\mathcal N:=\{f\in W^{1,\infty}(\mathbb{R}):
		\|f-u^Q\|_{W^{1,\infty}}\leq\varepsilon\}.
	\end{equation*}
	Every $f\in\mathcal N$ satisfies $f>a$ to the left of $J$,
	$f<a$ to its right, and $f_x\leq-M/2$ almost everywhere on $J$.
	It therefore has a unique threshold crossing in the interior of $J$.
	Choose $T>0$ such that
	\begin{equation*}
		C_0R I_{1/2}(T)<\frac{\delta}{2},
		\qquad C_1R^{1/2}I_{3/4}(T)<\frac{M_*}{2},
		\qquad \frac{C_0 I_{1/2}(T)}{M-M_*}<1.
	\end{equation*}
	For $w_0=f-u^Q$ with $f\in\mathcal N$,
	we have $\|w_0\|_\infty\leq\varepsilon<\delta/2$ and
	\begin{equation*}
		\|p^{-1}\|_\infty\|p\partial_xw_0\|_\infty
		\leq C_p\varepsilon<\frac{M_*}{2}.
	\end{equation*}
	Thus \eqref{eq:w}, \eqref{eq:wm}, and \eqref{eq:contract_est}
	hold on the same interval $[0,T]$ for every initial datum in $\mathcal N$.
	The local existence and comparison results apply with these fixed constants.

	We next prove that
	\begin{equation}
		\label{Uu-isolate}
		\{V\in\mathcal N:V\text{ is a weak solution of }\eqref{pblm}\}
		=\{u^Q\}.
	\end{equation}
	A stationary $V\in\mathcal N$ has a single threshold crossing $\tau\in J$,
	so it satisfies $KV=\mathbf{1}_{\{x<\tau\}}$ weakly.
	Since $V$ is bounded, Lemma~\ref{lem:positive_bounded_inverse} gives $V=V_\tau$.
	Lemma~\ref{lem:step_profiles} shows that this profile is strictly decreasing,
	with limits $1$ and $0$ at the respective infinities.
	Theorem~\ref{thm:monotone_classification} and the choice of $J$
	then give $V=u^Q$.

	Let $\eta>0$ be given.
	Choose $\sigma>0$ sufficiently small that $u_0:=\underline u_\sigma$ satisfies
	\begin{equation*}
		\|u_0-u^Q\|_{W^{1,\infty}}<\min\{\eta,\varepsilon_*\}.
	\end{equation*}
	The local theory gives a solution on $[0,T]$ with a continuous crossing path in $J$.
	If its distance from $u^Q$ is at least $\varepsilon_*$ at some time,
	the conclusion of this theorem follows.
	Otherwise, $u(\cdot,T)\in\mathcal N$, so the same local construction applies
	on $[T,2T]$.
	Continuing in this way either gives the required time $t_\eta$
	or produces a global solution satisfying
	\begin{equation}
		\label{eq:unstable_stay_close}
		\|u(\cdot,t)-u^Q\|_{W^{1,\infty}}<\varepsilon_*,
		\qquad t\geq0.
	\end{equation}
	The solutions and their crossing paths agree at the interval endpoints
	by local uniqueness.
	Using the semigroup identity $\mathsf S(t+s)=\mathsf S(t)\mathsf S(s)$,
	the mild formulas on the two adjacent time intervals combine to give
	\eqref{eq:mild_w} on the extended time interval
	with the original initial perturbation $w_0$.
	We show that the global alternative is impossible.

	Since $u_0$ is a stationary mild subsolution,
	local comparison on successive intervals of length $T$ gives
	\begin{equation*}
		u_0\leq u(\cdot,t),\qquad t\geq0.
	\end{equation*}
	Fix $s>0$.
	The equation \eqref{pRD} is autonomous, so $u(\cdot,t)$ and $u(\cdot,t+s)$
	are solutions with ordered initial data $u_0\leq u(\cdot,s)$.
	All their initial profiles at later times belong to $\mathcal N$.
	Applying local comparison successively on intervals of length $T$ gives
	\begin{equation*}
		u(x,t)\leq u(x,t+s),\qquad x\in\mathbb{R},\quad t\geq0.
	\end{equation*}
	Thus $u$ is nondecreasing in time at every spatial point.

	The bound \eqref{eq:unstable_stay_close} gives a pointwise limit
	$u_\infty(x):=\lim_{t\to\infty}u(x,t)$ and
	\begin{equation*}
		|u(x,t)-u(y,t)|
		\leq C|x-y|,
		\qquad x,y\in\mathbb{R},\quad t\geq0,
	\end{equation*}
	where $C:=\|u_x^Q\|_\infty+\varepsilon_*$ is independent of $t$.
	The limit $u_\infty$ has the same Lipschitz bound.
	On each compact interval $I_0\subset\mathbb{R}$,
	the Arzel\`a--Ascoli theorem therefore gives relative compactness in $C(I_0)$.
	Pointwise convergence identifies every subsequential limit with $u_\infty$, so
	\begin{equation}
		\label{eq:loc_unif_U}
		u(\cdot,t)\longrightarrow u_\infty
		\quad\text{locally uniformly as }t\to\infty.
	\end{equation}
	Writing $w(t):=u(\cdot,t)-u^Q$ and $w_\infty:=u_\infty-u^Q$,
	we have $\|w_\infty\|_\infty\leq\varepsilon_*$.
	Passing to the limit in
	$|w(x,t)-w(y,t)|\leq\varepsilon_*|x-y|$
	also gives $\|\partial_xw_\infty\|_\infty\leq\varepsilon_*$.
	Hence $\|u_\infty-u^Q\|_{W^{1,\infty}}\leq2\varepsilon_*=\varepsilon$,
	and $u_\infty\in\mathcal N$.

	The crossing path is nondecreasing as well.
	Indeed, monotonicity in time gives
	$u(\rho(t),t+s)\geq u(\rho(t),t)=a$.
	Since $u(\cdot,t+s)$ has the unique threshold crossing $\rho(t+s)$
	and satisfies
	$u(x,t+s)<a$ for $x>\rho(t+s)$,
	we must have
	\begin{equation*}
		\rho(t)\leq\rho(t+s).
	\end{equation*}
	Since $\rho(t)\in J$, it converges to some $\rho_\infty\in J$.
	Local uniform convergence and continuity of $u_\infty$ imply $u_\infty(\rho_\infty)=a$.
	As $u_\infty\in\mathcal N$, this is its unique threshold crossing,
	with $u_\infty>a$ to the left and $u_\infty<a$ to the right.

	It remains to prove that $u_\infty$ is stationary.
	Lemma~\ref{lem:unstable_subsuper} gives
	$w_0:=u_0-u^Q\in H^1(\mathbb{R})\subset L^2(\mathbb{R})$.
	Apply \eqref{eq:mild_w} with reference front $u^Q$ and reference crossing $\rho_Q$.
	Since $\|h_\rho(t)\|_2=|\rho(t)-\rho_Q|^{1/2}\leq R^{1/2}$,
	the contraction estimate on $L^2(\mathbb{R})$ gives
	\begin{equation*}
		\|w(t)\|_2\leq e^{-t}\|w_0\|_2+(1-e^{-t})R^{1/2}.
	\end{equation*}
	We also have $0<w_0\leq w(t)\uparrow w_\infty$ pointwise.
	Fatou's lemma therefore gives $w_\infty\in L^2(\mathbb{R})$,
	and dominated convergence gives
	\begin{equation*}
		w(t)\longrightarrow w_\infty\quad\text{in }L^2(\mathbb{R}).
	\end{equation*}
	Set $h_\infty:=\mathbf{1}_{\{x<\rho_\infty\}}
		-\mathbf{1}_{\{x<\rho_Q\}}$.
	For each fixed $s>0$, the mild formula starting at time $t$ reads
	\begin{equation*}
		w(t+s)=\mathsf S(s)w(t)+\int_0^{s} \mathsf S(s-r)h_\rho(t+r)\,dr.
	\end{equation*}
	Here, the $L^2(\mathbb{R})$ contraction estimate for $\mathsf S(t)$ gives
	\begin{align*}
		\left\|
		\int_0^{s} \mathsf S(s-r)
		\bigl(h_\rho(t+r)-h_\infty\bigr)\,dr
		\right\|_2
		 & \leq
		s\,|\rho(t)-\rho_\infty|^{1/2}
		\longrightarrow0.
	\end{align*}
	Passing to the limit in $L^2(\mathbb{R})$ in the mild formula
	and changing the integration variable gives
	\begin{equation*}
		w_\infty=\mathsf S(s)w_\infty+\int_0^{s} \mathsf S(r)h_\infty\,dr.
	\end{equation*}
	Consequently,
	\begin{equation*}
		\frac{\mathsf S(s)w_\infty-w_\infty}{s}
		=-\frac1{s}\int_0^{s} \mathsf S(r)h_\infty\,dr
		\longrightarrow-h_\infty
		\quad\text{in }L^2(\mathbb{R})\quad(s\downarrow0).
	\end{equation*}
	By the definition of the generator, $w_\infty\in\mathcal D(K)$
	and $Kw_\infty=h_\infty$.
	Adding the stationary equation for $u^Q$ yields
	\begin{equation*}
		Ku_\infty=\mathbf{1}_{\{x<\rho_\infty\}}.
	\end{equation*}
	The sign properties of $u_\infty$ identify this as its weak stationary McKean equation.
	Thus \eqref{Uu-isolate} gives $u_\infty=u^Q$.
	But comparison gives $u_\infty\geq u_0=\underline u_\sigma>u^Q$, a contradiction.

	Thus the solution reaches distance at least $\varepsilon_*$ from $u^Q$
	at a finite time $t_\eta>0$.
	Since $R$ and $\varepsilon_*$ are independent of $\eta$,
	this proves the assertion.
\end{proof}

\subsection{Transfer to the necklace graph}
\label{subsec:graph_transfer}

Section~\ref{sec:switchingODE} showed that, for $\alpha=1/4$,
the coordinate change $x=2y$ on the parallel edges turns the necklace evolution,
for graph profiles $u_\Gamma$ that agree on each pair of parallel edges,
into the scalar equation \eqref{pRD}.
The identification also holds for the mild solutions used above.
The coordinate change is an isometry from $L^2(\mathbb R)$
onto the paired-edge subspace of $L^2(\Gamma)$ and preserves the diffusion form,
because the two parallel-edge integrals become one integral over the doubled interval,
with derivative coefficient $4$.
It therefore identifies the diffusion semigroups.
It also preserves the $L^1$ and $L^\infty$ norms,
so, by density and duality, it identifies the consistent extensions
on $L^1(\mathbb R)$ and $L^\infty(\mathbb R)$
with their counterparts on the paired-edge subspaces of
$L^1(\Gamma)$ and $L^\infty(\Gamma)$.
Both the initial term and the forcing integral in the mild formulation \eqref{eq:mild_w}
therefore transfer under this change.

Take the norm on $W^{1,\infty}(\Gamma)$ to be the sum
of the supremum of the absolute edge values
and the essential supremum of the absolute edge derivatives.
Then
\begin{equation}
    \label{eq:graph_norm_equivalence}
    \|u\|_{W^{1,\infty}(\mathbb R)}
    \leq\|u_\Gamma\|_{W^{1,\infty}(\Gamma)}
    \leq2\|u\|_{W^{1,\infty}(\mathbb R)},
\end{equation}
and the same estimates hold for differences of profiles.
Consequently, Theorems~\ref{thm:stable_front_stability} and~\ref{thm:unstable_front_instability}
hold for the necklace equation with perturbations
that preserve equality on each pair of parallel edges.

\section{Nonmonotone stationary fronts}
\label{sec:nonmonotone}

For every interior threshold $a_-<a<a_+$,
Theorem~\ref{thm:monotone_classification} gives two decreasing fronts, $u^P$ and $u^Q$,
and Theorem~\ref{thm:all_front_existence} shows that no front of any kind exists outside $[a_-,a_+]$
and that only decreasing fronts exist at its endpoints.
This section shows that nonmonotone fronts exist at every interior threshold
and counts those whose threshold crossings are far apart.
We prescribe the crossings, solve the linear equation with the corresponding alternating forcing,
and require the solution to take the value $a$ at the prescribed points
(Section~\ref{subsec:alternating_matching}).
When the crossings lie in widely separated periods,
each downward crossing lies close to the crossing phase of one of the two decreasing fronts at threshold $a$,
each upward crossing lies close to one of the two such phases at threshold $1-a$,
and these choices can be made independently
(Section~\ref{subsec:nonmonotone_fronts}).
Section~\ref{subsec:central_parameters} determines the sharp minimum slot gap
for one symmetric family at $a=1/2$,
and Section~\ref{subsec:global_manifold_illustration} shows the corresponding picture
for the period map.

\subsection{Matching alternating interfaces}
\label{subsec:alternating_matching}

We express a profile with prescribed crossings through the bounded linear solutions
$V_\rho$ of Lemma~\ref{lem:step_profiles}, which satisfy
$KV_\rho=\mathbf 1_{\{x<\rho\}}$ and $V_\rho(\rho)=b(\rho)$.
The decreasing profile $V_\rho$ is a stationary front at threshold $a$ precisely when $b(\rho)=a$.
Its complement $1-V_\rho$ solves $K(1-V_\rho)=\mathbf 1_{\{x>\rho\}}$,
so it is an increasing front at threshold $a$ precisely when $b(\rho)=1-a$.
For several crossings, an alternating sum of the $V_\rho$ prescribes the forcing on the successive intervals,
and matching the value $a$ at the prescribed points gives a stationary front.

\begin{lemma}[Matching alternating interfaces]
	\label{lem:alternating_matching}
	Let $m\geq1$ be odd, let $\boldsymbol\rho=(\rho_1,\ldots,\rho_m)$ with $\rho_1<\cdots<\rho_m$, and set $\sigma_j:=(-1)^{j+1}$.
	Define
	\begin{equation}
		\label{eq:alternating_superposition}
		W_{\boldsymbol\rho}(x):=
		\sum_{j=1}^m\sigma_jV_{\rho_j}(x).
	\end{equation}
	There is a stationary McKean front whose threshold set is exactly $\{\rho_1,\ldots,\rho_m\}$ if and only if
	\begin{equation}
		\label{eq:multicrossing_matching}
		W_{\boldsymbol\rho}(\rho_i)=a,
		\qquad i=1,\ldots,m.
	\end{equation}
	That front is $W_{\boldsymbol\rho}$, and its crossing flux is negative for odd $i$ and positive for even $i$.
\end{lemma}

\begin{proof}
	Linearity gives
	\begin{equation}
		\label{eq:alternating_forcing}
		KW_{\boldsymbol\rho}
		=\sum_{j=1}^m\sigma_j\mathbf 1_{\{x<\rho_j\}}
		=\mathbf 1_{\mathcal A_{\boldsymbol\rho}},
		\qquad
		\mathcal A_{\boldsymbol\rho}:=(-\infty,\rho_1)
		\cup\bigcup_{k=1}^{(m-1)/2}(\rho_{2k},\rho_{2k+1}).
	\end{equation}
	Both $\mathcal A_{\boldsymbol\rho}$ and its complement contain open intervals.
	Strict comparison with the constant forcings $0$ and $1$ therefore gives $0<W_{\boldsymbol\rho}<1$ everywhere.

	Suppose that \eqref{eq:multicrossing_matching} holds, and write $W:=W_{\boldsymbol\rho}$ and $v:=pW_x$.
	On an interval $(\rho_i,\rho_{i+1})$ with odd $i$, \eqref{eq:alternating_forcing} gives $v'=W>0$ almost everywhere, so $v$ is strictly increasing.
	Since $W(\rho_i)=W(\rho_{i+1})=a$ and $W_x=p^{-1}v$,
	\begin{equation*}
		0=W(\rho_{i+1})-W(\rho_i)
		=\int_{\rho_i}^{\rho_{i+1}}p(x)^{-1}v(x)\,dx.
	\end{equation*}
	It follows that $v(\rho_i)<0<v(\rho_{i+1})$ and that $v$ has a unique zero in the interval.
	Thus $W$ first decreases and then increases, with $W<a$ throughout $(\rho_i,\rho_{i+1})$.
	For even $i$, the same argument uses $v'=W-1<0$.
	It gives $v(\rho_i)>0>v(\rho_{i+1})$ and $W>a$ throughout the interval.

	By \eqref{eq:alternating_forcing}, $W-1$ solves the homogeneous equation on $(-\infty,\rho_1)$, and $W$ solves it on $(\rho_m,\infty)$.
	As in the derivation of \eqref{eq:step_profile_tails}, $W-1$ and $W$ are linear combinations of $y_+$ and $y_-$ there, and boundedness excludes the $y_-$ term on the left and the $y_+$ term on the right.
	The values $W(\rho_1)=W(\rho_m)=a$ fix the remaining coefficients, so
	\begin{equation*}
		W(x)=1-(1-a)\frac{y_+(x)}{y_+(\rho_1)}
		\quad(x\leq\rho_1),
		\qquad
		W(x)=a\frac{y_-(x)}{y_-(\rho_m)}
		\quad(x\geq\rho_m).
	\end{equation*}
	These formulas give $W>a$ to the left of $\rho_1$ and $W<a$ to the right of $\rho_m$, with negative flux at both exterior crossings.
	They also cover the case $m=1$.
	The period identities \eqref{eq:tail_floquet_relations} give the limits $1$ and $0$ at the respective infinities.
	Consequently, the prescribed points are all the threshold crossings, and each has the asserted nonzero flux.
	The forcing in \eqref{eq:alternating_forcing} agrees almost everywhere with $\mathbf 1_{\{W\geq a\}}$, so $W$ solves the stationary McKean equation.

	Conversely, let $u$ be a stationary front with threshold set $\{\rho_1,\ldots,\rho_m\}$.
	By Lemma~\ref{lem:mckean_switching}, these crossings are transverse; their directions therefore alternate, starting downward because $u(-\infty)=1$.
	Hence $Ku=\mathbf 1_{\mathcal A_{\boldsymbol\rho}}$ almost everywhere.
	Uniqueness in Lemma~\ref{lem:positive_bounded_inverse} gives $u=W_{\boldsymbol\rho}$, and evaluating at the crossings gives \eqref{eq:multicrossing_matching}.
\end{proof}

\subsection{Separated fronts}
\label{subsec:nonmonotone_fronts}

For $\theta\in(a_-,a_+)$, write
$\rho_P(\theta)\in(0,L_1)$ and $\rho_Q(\theta)\in(L_1,L)$
for the two solutions of $b(s)=\theta$ (Lemma~\ref{lem:matching_function}).
For $\theta=a$ they are the crossing phases of $u^P$ and $u^Q$.
For $\theta=1-a$ they are the crossing phases of the two decreasing fronts at threshold $1-a$.
The complements of these two fronts are the two increasing fronts at threshold $a$,
with the same crossing phases.
The labels $P$ and $Q$ follow the stable and unstable decreasing fronts
of Theorems~\ref{thm:stable_front_stability} and~\ref{thm:unstable_front_instability}.
Complementation preserves these single-crossing conclusions,
so the labels also mark the stable and unstable increasing fronts.
They refer to the isolated fronts.
We do not study the stability of the nonmonotone fronts constructed below
(see Section~\ref{sec:discussion}).

A nonmonotone front with $m$ crossings alternates downward and upward transitions,
starting and ending downward, so $m$ is odd.
We place its crossings in prescribed periods and choose one isolated phase for each crossing.

\begin{definition}[Period slots and crossing words]
	\label{def:slots_words}
	Let $c\in\mathbb R$.
	The \emph{period slots} with cut $c$ are the intervals $[NL+c,(N+1)L+c)$, $N\in\mathbb Z$.
	For $c=0$ they are the periods $[NL,(N+1)L)$.
	Fix an odd $m\geq3$ and integers $N_1<\dots<N_m$, and set
	\begin{equation}
		\label{eq:separated_target_phases}
		\theta_i:=
		\begin{cases}
			a,   & i\text{ odd},  \\
			1-a, & i\text{ even}.
		\end{cases}
	\end{equation}
	A \emph{crossing word} is a sequence $X=(X_1,\ldots,X_m)\in\{P,Q\}^m$.
	The letter $X_i$ selects the isolated phase $\rho_{X_i}(\theta_i)$ for the crossing in slot $N_i$.
	When the slots matter, we write $X_i$ as $P_{N_i}$ or $Q_{N_i}$
	and the word as a product, for example $Q_0P_2P_4$.
\end{definition}

\begin{theorem}[Two choices in each separated crossing period]
	\label{thm:separated_fronts}
	Fix $0<\alpha<1$, $L_1,L_2>0$, $a\in(a_-,a_+)$, and an odd integer $m\geq3$.
	There are an integer $n_m\geq1$ and a constant $C_m>0$ with the following property.
	Let $N_1<\cdots<N_m$ be integers with
	\begin{equation*}
		n:=\min_{1\leq i<m}(N_{i+1}-N_i)\geq n_m.
	\end{equation*}
	Then there are exactly $2^m$ stationary fronts whose threshold set consists of one point
	in each period $[N_iL,(N_i+1)L)$ and no points elsewhere.
	They are indexed by the crossing words $X\in\{P,Q\}^m$.
	Writing the crossings of the front with word $X$ as $\rho_i=N_iL+s_i$, $s_i\in[0,L)$,
	and setting $s_i^*:=\rho_{X_i}(\theta_i)$, we have
	\begin{equation}
		\label{eq:separated_phase_estimate}
		\max_i|s_i-s_i^*|
		\leq C_m \lambda^{-n}.
	\end{equation}
	The constants depend only on $\alpha,L_1,L_2,a$ and $m$.
	In particular, there are infinitely many such fronts, inequivalent under translations by integer multiples of $L$.
\end{theorem}

\begin{proof}
	Since $a_-+a_+=1$, both $a$ and $1-a$ lie in $(a_-,a_+)$.
	Lemma~\ref{lem:matching_function} therefore gives exactly two simple roots of $b(s)=\theta_i$ for each $i$.
	We show that the equations for separated crossings have precisely one solution near each of the $2^m$ selections of these roots, that is, near each word.

	By \eqref{eq:tail_floquet_relations}, the functions $y_+(x)e^{-\gamma x}$ and $y_-(x)e^{\gamma x}$ are positive and $L$-periodic,
	with $\gamma$ as in \eqref{eq:green_decay}.
	Together with the tail formulas \eqref{eq:step_profile_tails}, this gives
	\begin{equation}
		\label{eq:step_tail_estimates}
		\left|V_\rho(x)-\mathbf 1_{\{x<\rho\}}\right|
		+|p(x)\partial_xV_\rho(x)|
		\leq C e^{-\gamma|x-\rho|}
	\end{equation}
	for every $\rho$ and almost every $x$.
	In addition, \eqref{eq:step_green_representation} and \eqref{eq:green_decay} give
	\begin{equation}
		\label{eq:step_interface_derivative}
		\partial_\rho V_\rho(x)=G(x,\rho),
		\qquad
		|\partial_\rho V_\rho(x)|\leq C e^{-\gamma|x-\rho|}.
	\end{equation}
	Both estimates are uniform in the interface phase.

	For $s=(s_1,\ldots,s_m)\in[0,L]^m$, put $\rho_i=N_iL+s_i$ and $F_{\boldsymbol N,i}(s):=W_{\boldsymbol\rho}(\rho_i)-a$.
	For $n\geq2$, these interfaces are ordered throughout the cube.
	Since $\sum_{j>i}\sigma_j=(1-\sigma_i)/2$, the matching equations become
	\begin{equation}
		\label{eq:multicrossing_remainder}
		\begin{aligned}
			F_{\boldsymbol N,i}(s)
			 & =\sigma_i\bigl(b(s_i)-\theta_i\bigr)+R_{\boldsymbol N,i}(s), \\
			R_{\boldsymbol N,i}(s)
			 & :=\sum_{j<i}\sigma_jV_{\rho_j}(\rho_i)
			+\sum_{j>i}\sigma_j\bigl(V_{\rho_j}(\rho_i)-1\bigr).
		\end{aligned}
	\end{equation}
	Because $|\rho_i-\rho_j|\geq(n-1)L$ for $i\neq j$ and $e^{-\gamma L}=\lambda^{-1}$, \eqref{eq:step_tail_estimates} gives
	\begin{equation}
		\label{eq:matching_c0_estimate}
		\sup_{s\in[0,L]^m}\|R_{\boldsymbol N}(s)\|_\infty
		\leq C_m \lambda^{-n}.
	\end{equation}

	Choose a small closed interval around each of the two roots of $b(s)=\theta_i$, lying in its open coefficient interval.
	Take all these intervals to have a common radius $\delta_0>0$, and let $I_X$ be their product for the word $X\in\{P,Q\}^m$.
	Write $F_{\infty,i}(s):=\sigma_i(b(s_i)-\theta_i)$.
	At the center $s^*$ of $I_X$, we have $F_\infty(s^*)=0$ and the invertible diagonal matrix
	\begin{equation*}
		J_*:=DF_\infty(s^*)
		=\operatorname{diag}\bigl(\sigma_i b'(s_i^*)\bigr)_{i=1}^m.
	\end{equation*}
	There are only finitely many words, so we may choose $\delta_0$ small enough that on every $I_X$,
	\begin{equation*}
		\sup_{s\in I_X}
		\|1-J_*^{-1}DF_\infty(s)\|\leq\frac14,
	\end{equation*}
	where $1$ is the identity on $\mathbb R^m$ and the matrix norm is induced by the maximum norm.

	For $j\neq i$, differentiating $V_{\rho_j}(\rho_i)$ with respect to $s_i$ gives $\partial_xV_{\rho_j}(\rho_i)$, and differentiating with respect to $s_j$ gives $G(\rho_i,\rho_j)$.
	The other phase derivatives vanish.
	These derivatives are continuous on each $I_X$, since the evaluation points remain inside their coefficient intervals.
	Equations~\eqref{eq:step_tail_estimates}--\eqref{eq:step_interface_derivative} therefore give
	\begin{equation}
		\label{eq:matching_c1_estimate}
		\|R_{\boldsymbol N}\|_{C^1(I_X)}\leq C_m \lambda^{-n},
	\end{equation}
	uniformly over the words.
	For large enough $n$, it follows that
	\begin{equation*}
		\sup_{s\in I_X}\|J_*^{-1}DR_{\boldsymbol N}(s)\|\leq\frac14,
		\qquad
		\|J_*^{-1}R_{\boldsymbol N}(s^*)\|_\infty\leq\frac{\delta_0}2.
	\end{equation*}
	Hence $T_{\boldsymbol N}(s):=s-J_*^{-1}F_{\boldsymbol N}(s)$ is a contraction of $I_X$ into itself, with Lipschitz constant at most $1/2$.
	Its unique fixed point solves $F_{\boldsymbol N}(s)=0$ and satisfies
	\begin{equation*}
		\|s-s^*\|_\infty
		\leq2\|J_*^{-1}R_{\boldsymbol N}(s^*)\|_\infty
		\leq C_m \lambda^{-n}.
	\end{equation*}
	Thus every word gives exactly one matching solution in its product of intervals, with the estimate \eqref{eq:separated_phase_estimate}.

	It remains to exclude matching solutions elsewhere in the selected periods.
	For each $i$, the continuous function $|b(s)-\theta_i|$ has a positive lower bound outside the interiors of its two chosen root intervals.
	Let $\eta>0$ be a common lower bound for these finitely many complements in $[0,L]$.
	Increase the gap $n$ so that the right-hand side of \eqref{eq:matching_c0_estimate} is less than $\eta$.
	If $F_{\boldsymbol N}(s)=0$, then \eqref{eq:multicrossing_remainder} implies $|b(s_i)-\theta_i|<\eta$ for every $i$.
	Consequently, each $s_i$ lies in one of its two root intervals, and $s$ lies in one of the $2^m$ products already considered.
	There are therefore exactly $2^m$ matching solutions in the full phase cube.
	Lemma~\ref{lem:alternating_matching}, including its converse, identifies these solutions with all fronts whose threshold set consists of one point in each selected period and no points elsewhere.
	Their phase vectors lie in disjoint products, so the fronts are distinct.
	They are also inequivalent under period translations, since a nonzero common shift would change their first selected period.

	Finally, take $N_i=(i-1)k$ with an integer $k\geq n_m$, fix a word, and let $k\to\infty$.
	The first crossing gap satisfies
	\begin{equation*}
		\rho_2-\rho_1
		=kL+s_2^*-s_1^*+\mathcal O(\lambda^{-k})\longrightarrow\infty.
	\end{equation*}
	Crossing gaps are invariant under a common translation, so these fronts represent infinitely many translation classes.
\end{proof}

\begin{remark}[Slots with a cut]
	\label{rem:cut}
	Theorem~\ref{thm:separated_fronts} also holds for the slots $[N_iL+c,(N_i+1)L+c)$
	with any cut $c$ not congruent modulo $L$ to one of the four isolated phases
	$\rho_P(a)$, $\rho_Q(a)$, $\rho_P(1-a)$, $\rho_Q(1-a)$.
	In the proof, use the representatives of these phases in $(c,c+L)$
	and replace the phase cube by $[c,c+L]^m$.
	The cut is then fixed, and $n_m$ and $C_m$ may also depend on it.
	A suitable cut keeps both choices for a crossing in one slot
	when they lie on either side of a coefficient interface.
\end{remark}

\begin{figure}[tbp]
    \centering
    \includegraphics[width=\textwidth]{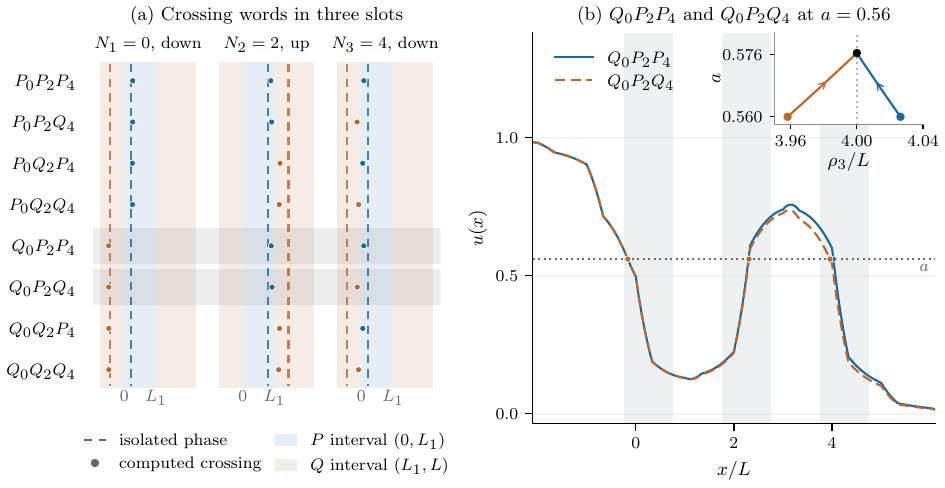}
    \caption{Crossing choices and a pair of colliding fronts,
        for $\alpha=1/4$, $L_1=0.8$, $L_2=1.6$, and $a=0.56$.
        (a) Eight computed crossing words in the slots $N_iL+[-0.6,1.8)$,
        with $(N_1,N_2,N_3)=(0,2,4)$.
        The horizontal coordinate in each column is the phase $x-N_iL$.
        A negative phase near zero represents a $Q$ crossing from the preceding coefficient period.
        The blue and orange regions give the coefficient intervals for $P$ and $Q$, respectively.
        Dashed lines mark the isolated phases $\rho_P(\theta_i)$ and $\rho_Q(\theta_i)$, and dots mark the crossings of the computed fronts.
        The middle crossing is upward, so its isolated phases solve $b=1-a$.
        (b) The fronts $Q_0P_2P_4$ and $Q_0P_2Q_4$, whose rows are highlighted in (a).
        Grey bands mark the three chosen slots.
        The inset follows their last crossings as $a$ increases.
        In this computation the two fronts meet at $\rho_3=4L$ when $a\approx0.576324$.
        All three crossing locations vary along each branch.
        These profiles and branches are computed from the full matching equations.
        Theorem~\ref{thm:separated_fronts} gives the exact count for sufficiently separated slots.}
    \label{fig:global_many_crossings}
\end{figure}

Figure~\ref{fig:global_many_crossings} displays all eight words for three crossings
in the slots with indices $0,2,4$, for $L_1=0.8$, $L_2=1.6$ (so $L=2.4$) and the cut $c=-0.6$.
The slots are
\begin{equation*}
	[-0.6,1.8),\qquad [4.2,6.6),\qquad [9.0,11.4).
\end{equation*}
A crossing at $x=4L-\delta$ with $0<\delta<0.6$ lies in slot $4$, with phase $-\delta$.
Its phase modulo $L$ is $L-\delta\in(L_1,L)$, so its letter is $Q_4$.
With the cut $c=0$ it would lie in slot $3$.
The subscript of a letter thus names the slot,
and the letter names the coefficient interval modulo $L$.
By Lemma~\ref{lem:alternating_matching}, the crossing sets of the fronts with $m$ crossings are the solutions $\rho_1<\cdots<\rho_m$ of the $m$ equations \eqref{eq:multicrossing_matching}.
Each $V_{\rho_j}$ in them is given by \eqref{eq:step_profile_tails} with $V_{\rho_j}(\rho_j)=b(\rho_j)$,
and this expression is explicit because the functions $y_\pm$ in it are.
Indeed, $y_\pm$ are the first components of $Y_\pm$, which on $[0,L]$ are obtained by applying the matrices \eqref{time_t} and \eqref{time_t_alpha} to the eigenvectors of $A_0$,
and \eqref{eq:tail_floquet_relations} extends them to all of $\mathbb R$.
The crossing locations in the figure are computed by a root finder for the equations \eqref{eq:multicrossing_matching}, starting from the isolated phases, from which they differ slightly.
The front is then given by \eqref{eq:alternating_superposition},
and on each interval between consecutive crossings and coefficient interfaces it is a constant plus an explicit combination of hyperbolic functions.

Together with Theorem~\ref{thm:all_front_existence},
Theorem~\ref{thm:separated_fronts} shows that nonmonotone stationary fronts exist
precisely when the two decreasing front classes coexist, that is, for $a_-<a<a_+$.
The estimate \eqref{eq:separated_phase_estimate} bounds the displacement of each crossing
caused by the other crossings.

In Theorem~\ref{thm:separated_fronts}, the minimum slot gap $n_m$ and the constant $C_m$ may depend on $m$.
For thresholds away from the endpoints of the pinning interval, they can be chosen independently of $m$.

\begin{corollary}[Uniform existence away from the endpoints of the pinning interval]
    \label{cor:uniform_central_fronts}
    Fix $0<\alpha<1$ and $L_1,L_2>0$.
    There are constants $C_{\mathrm{gap}},C>0$ and an integer $N_0\geq1$, depending only on $\alpha$, $L_1$ and $L_2$, with the following property.
    Let $N\geq N_0$ be an integer, let $m\geq3$ be odd, and let $N_1<\cdots<N_m$ be integers satisfying
    $N_{i+1}-N_i\geq N$.
    For every threshold
    \begin{equation}
        \label{eq:uniform_central_band}
        a\in[a_-+C_{\mathrm{gap}}\lambda^{-N},a_+-C_{\mathrm{gap}}\lambda^{-N}],
    \end{equation}
    the conclusions of Theorem~\ref{thm:separated_fronts} hold for the slots $N_1,\ldots,N_m$,
    with the phase estimate \eqref{eq:separated_phase_estimate} replaced by
    \begin{equation}
        \label{eq:uniform_central_phase_bound}
        \max_{1\leq i\leq m}|s_i-s_i^*|
        \leq C\lambda^{-N}.
    \end{equation}
\end{corollary}

\begin{proof}
    By Lemma~\ref{lem:matching_function}, the derivative \eqref{eq:matching_derivative} of $b$
    has nonzero one-sided limits at the coefficient interfaces,
    with the same strict signs as on the adjoining open intervals.
    Compactness of the two closed coefficient intervals therefore gives constants $\beta,M>0$ such that
    $\beta\leq|b'|\leq M$ on their respective one-sided extensions.
    The second derivatives are bounded on the same extensions.

    Put $\delta:=\min\{a-a_-,a_+-a\}$.
    Each root of $b(s)=a$ or $b(s)=1-a$ lies at distance at least $\delta/M$ from the ends of its coefficient interval.
    The root intervals in the proof of Theorem~\ref{thm:separated_fronts} may therefore be chosen with common radius
    \begin{equation*}
        \delta_*:=\min\{\delta_1,\delta/(4M)\},
    \end{equation*}
    where $\delta_1>0$ depends only on $\alpha,L_1,L_2$ and is small enough that
    $\|1-J_*^{-1}DF_\infty\|\leq1/4$ on every resulting product.
    This choice is uniform in the parameter and in $m$, since $J_*$ is diagonal and $\|J_*^{-1}\|\leq\beta^{-1}$.

    The off-diagonal terms also admit a bound independent of $m$.
    For $\rho_i=N_iL+s_i$ with $s_i\in[0,L]$, we have
    $|\rho_i-\rho_j|\geq|i-j|NL-L$.
    Summing the tail and derivative estimates in \eqref{eq:step_tail_estimates}--\eqref{eq:step_interface_derivative} gives
    \begin{equation*}
        \|R_{\boldsymbol N}\|_{C^0([0,L]^m)}
        +\|DR_{\boldsymbol N}\|_{L^\infty(I_X)}
        \leq C'\sum_{j=1}^{\infty}\lambda^{-jN}
        =\frac{C'\lambda^{-N}}{1-\lambda^{-N}}.
    \end{equation*}
    The derivative norm is the matrix norm induced by the maximum norm.
    The constant $C'$ absorbs the fixed factor $e^{\gamma L}=\lambda$ and is independent of $m$.

    Choose $N_0$ large enough to make the derivative contribution to the contraction less than $1/4$ and the displacement at the center less than $\delta_1/2$.
    Then choose $C_{\mathrm{gap}}$ large enough that $\delta\geq C_{\mathrm{gap}}\lambda^{-N}$ also makes this displacement less than $\delta/(8M)$.
    The contraction argument gives one solution in each product of root intervals and the bound \eqref{eq:uniform_central_phase_bound}.
    Outside these products, some coordinate has $|b(s_i)-\theta_i|\geq\beta\delta_*$.
    Increasing $C_{\mathrm{gap}},N_0$ if necessary makes the remainder strictly smaller than this bound, excluding any further matching solution.
    Lemma~\ref{lem:alternating_matching} completes the count.
\end{proof}

Thus a sufficient slot gap for a fixed interior threshold is
\begin{equation}
    \label{eq:central_sufficient_gap}
    N\geq\max\left\{N_0,
        \left\lceil\frac{\log(C_{\mathrm{gap}}/\min\{a-a_-,a_+-a\})}{\log\lambda}\right\rceil
        \right\}.
\end{equation}
In particular, one finite slot gap works for every odd crossing number throughout any closed subinterval of $(a_-,a_+)$.
For smaller gaps, the matching equations remain exact, but their solutions need not be small perturbations of independently chosen isolated phases.

\subsection{The balanced minimum gap}
\label{subsec:central_parameters}

Theorem~\ref{thm:separated_fronts} and Corollary~\ref{cor:uniform_central_fronts} give sufficient slot gaps only.
At $a=1/2$, the sharp minimum slot gap can be determined for reflection-symmetric fronts with the word $P_{-N}P_0P_N$.
Proposition~\ref{prop:balanced_minimum_gap} shows that such a front exists only for $N\geq N_*$,
and that, for fixed $L_1$ and $L_2$, $N_*$ depends on the diffusion contrast and grows without bound as $\alpha\uparrow1$.
Put $s_c:=L_1/2$ and $r_c:=r_+(s_c)$.
Since $\bar s_c=s_c$, Lemma~\ref{lem:floquet_slopes} gives $r_c=r_-(s_c)>1$,
and Remark~\ref{rem:balanced_phases} gives $\rho_P=s_c$ at $a=1/2$.
For $0\leq\delta\leq s_c$, the coefficient equals $1$ on $[s_c,s_c+\delta]$, so
\begin{equation*}
    Y(\delta):=\frac{y_+(s_c+\delta)}{y_+(s_c)}=\cosh\delta+r_c\sinh\delta.
\end{equation*}
Define
\begin{equation}
    \label{eq:balanced_AB}
    A(\delta):=\frac{(r_c-r_c^{-1})\sinh(2\delta)}{2Y(\delta)},
    \qquad
    B(\delta):=\cosh\delta+r_c^{-1}\sinh\delta.
\end{equation}

\begin{proposition}[A minimum slot gap at the balanced parameter]
    \label{prop:balanced_minimum_gap}
    Let $a=1/2$, and let $s_c$, $r_c$, $A$ and $B$ be as above.
    For an integer $N\geq1$, there is a unique stationary front at $a=1/2$
    satisfying $u(x)=1-u(L_1-x)$, with exactly three crossings in the periods
    $-N,0,N$ and all crossing phases in $(0,L_1)$, if and only if
    \begin{equation}
        \label{eq:balanced_gap_criterion}
        A(s_c)\lambda^N+B(s_c)\lambda^{-N}>1.
    \end{equation}
    Its crossings are
    \begin{equation*}
        \rho_1=s_c-NL+\delta_N,\qquad
        \rho_2=s_c,\qquad
        \rho_3=s_c+NL-\delta_N,
        \qquad 0<\delta_N<s_c.
    \end{equation*}
    Let $N_*$ be the least integer satisfying the strict inequality \eqref{eq:balanced_gap_criterion}.
    The integers satisfying it are precisely $N_*,N_*+1,\ldots$,
    and $\delta_N$ strictly decreases to zero on this set.
    If equality holds in \eqref{eq:balanced_gap_criterion} at an integer $N$,
    there is a front with the symmetry $u(x)=1-u(L_1-x)$ and $\delta_N=s_c$, whose two exterior crossings lie at coefficient interfaces.
    This front is excluded by the condition on the crossing phases, and then $N_*=N+1$.
    For fixed $L_1,L_2>0$, the minimum $N_*$ tends to infinity as $\alpha\uparrow1$.
\end{proposition}

The criterion follows by matching the bounded left tail to the middle crossing.
Reflection fixes this crossing at $s_c$ and places the exterior crossings at
$s_c\pm(NL-\delta)$.
Put $q:=\lambda^{-N}$.
Propagating the left-tail state from the first crossing to $s_c$
reduces the condition $u(s_c)=1/2$ at the middle crossing to the scalar equation
\begin{equation}
    \label{eq:balanced_scalar_matching}
    R_q(\delta):=\frac{A(\delta)}q+B(\delta)q=1.
\end{equation}
On $[0,s_c]$, the function $R_q$ is strictly increasing,
with $R_q(0)=q<1$ and $R_q(s_c)=A(s_c)/q+B(s_c)q$.
Thus \eqref{eq:balanced_gap_criterion} places its unique matching root in $(0,s_c)$;
equality gives $\delta=s_c$, with both exterior crossings at coefficient interfaces.
Appendix~\ref{app:balanced_minimum_gap} gives the calculation
and proves the remaining assertions.

For example, with $L_1=0.8$ and $L_2=1.6$, direct evaluation of
\eqref{eq:balanced_gap_criterion} gives $N_*=1$ at $\alpha=0.25$,
$N_*=2$ at $\alpha=0.5$, and $N_*=3$ at $\alpha=0.99$.
Thus even at the balanced threshold, the minimum gap for this symmetric word depends on the diffusion contrast.
For arbitrary words, Corollary~\ref{cor:uniform_central_fronts} gives a common sufficient bound.

\subsection{The fronts in the phase plane of the period map}
\label{subsec:global_manifold_illustration}

\begin{figure}[!htbp]
	\centering
	\includegraphics[width=\textwidth]{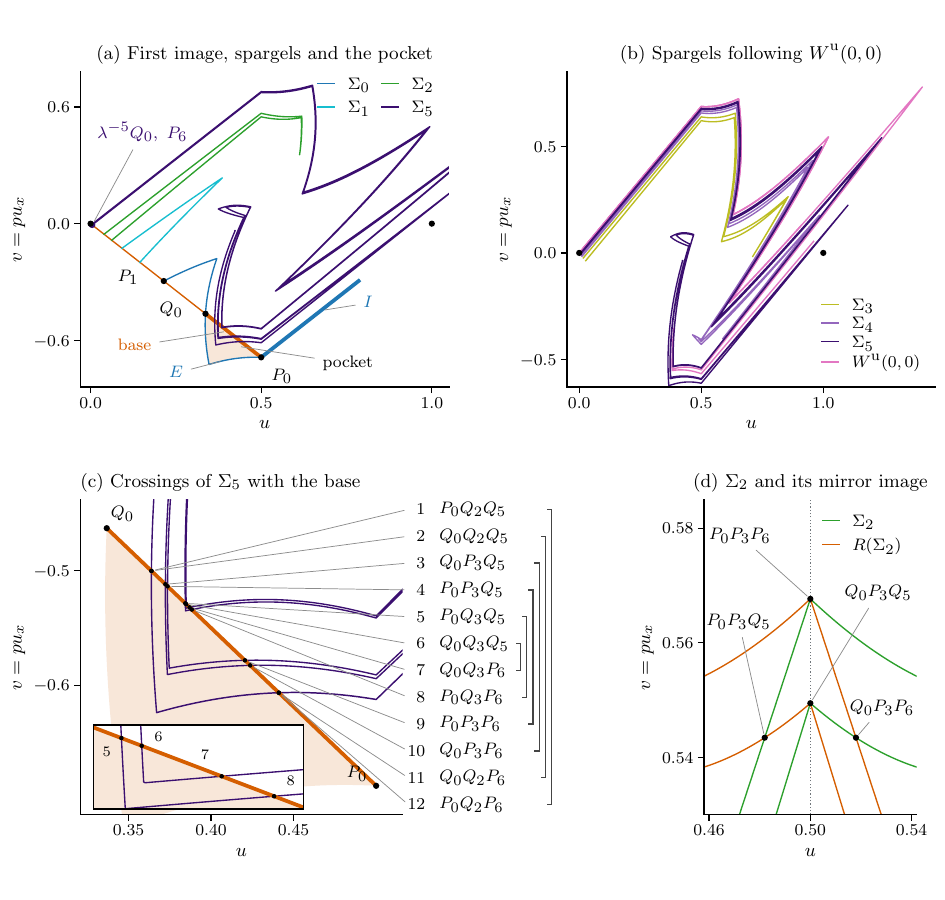}
	\caption{Computed period-map images of the unstable manifold of $(1,0)$
		for $\alpha=1/4$, $L_1=0.4$, $L_2=0.8$, $a=1/2$ and the sampling phase $s=L_1/2$.
		(a) The straight unstable segment $I$, its first image, consisting of the arc $E$ and the spargel $\Sigma_0$,
		the spargels $\Sigma_1$ and $\Sigma_2$, and the part of the spargel $\Sigma_5$ inside the frame, starting from its endpoints $\lambda^{-5}Q_0$ and $P_6$ near $(0,0)$.
		The shaded pocket is bounded by $E$ and the base $[Q_0,P_0]$ (Proposition~\ref{prop:pocket}).
		(b) The spargels $\Sigma_3$, $\Sigma_4$ and $\Sigma_5$ follow the computed part of the unstable manifold of $(0,0)$,
		consisting of its segment from $(0,0)$ along the unstable eigendirection to the threshold $u=a$, and the images of that segment under $\Psi$, $\Psi^2$ and $\Psi^3$ (Proposition~\ref{prop:pocket}(f)).
		(c) The twelve computed crossings of $\Sigma_5$ with the base, labelled by the crossing words of their fronts.
		The labels are numbered in the order of the crossings along the base from $Q_0$ to $P_0$,
		and brackets join the two crossings of each visit of the pocket.
		The inset magnifies crossings 5--8, which lie close together.
		The six visits are nested, and the two words of each visit differ only in their last letter.
		In this computation $\Sigma_3$ and $\Sigma_4$ do not cross the base.
		(d) The spargel $\Sigma_2$ and its mirror image $R(\Sigma_2)$ on the stable manifold of $(0,0)$.
		Of the four intersections shown, two lie on the threshold line and represent symmetric fronts,
		and the other two are exchanged by $R$.
		Words are obtained by tracing the computed trajectories to their equilibrium tails,
		backward in (c) and in both directions in (d), and shifting the first slot index to zero.
		Long mesh chords are omitted, so a visible break need not be a gap in the manifold.}
	\label{fig:global_first_return}
\end{figure}

\begin{figure}[!htbp]
	\centering
	\includegraphics[width=\textwidth]{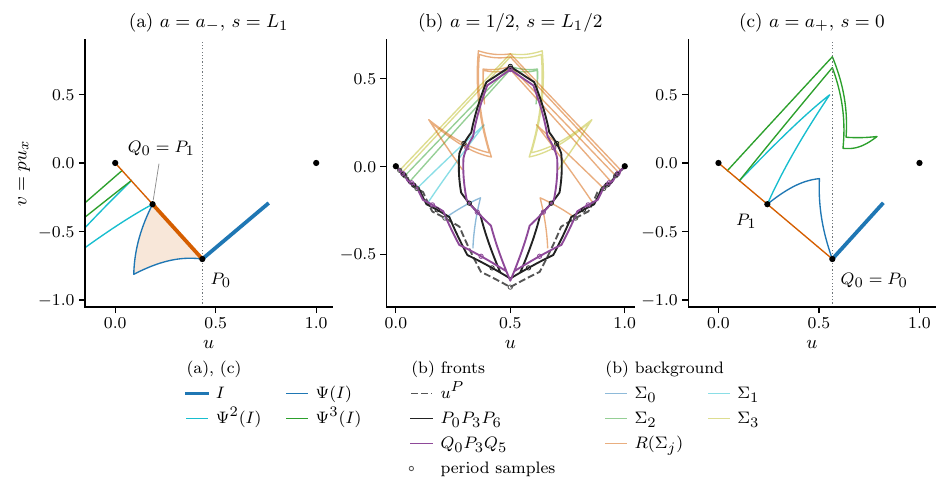}
	\caption{The phase plane of the period map at the endpoints of the pinning interval and at its centre,
		for $\alpha=1/4$, $L_1=0.4$ and $L_2=0.8$.
		Panels (a) and (c) show the straight unstable segment $I$ and its first three images,
		with the stable segment in orange.
		(a) At $a=a_-$, sampled at $s=L_1$,
		every unstable tail other than that of the decreasing front escapes,
		the pocket (shaded) is bounded by $\Psi(I)$ and the base, and $Q_0=P_1$.
		(b) At $a=1/2$, sampled at $s=L_1/2$, the continuous trajectories $x\mapsto(u(x),p(x)u_x(x))$
		of the decreasing front $u^P$ and of the symmetric fronts $P_0P_3P_6$ and $Q_0P_3Q_5$,
		with their samples at $s+kL$ marked by open circles.
		The trajectories have corners at the coefficient interfaces, where $u_x$ jumps and the flux is continuous.
		In the background are the spargels $\Sigma_0,\dots,\Sigma_3$
		and their mirror images $R(\Sigma_j)$ under $R(u,v)=(1-u,v)$.
		(c) At $a=a_+$, sampled at $s=0$, the base and the pocket collapse to $P_0=Q_0$.}
	\label{fig:global_endpoints}
\end{figure}

By Proposition~\ref{prop:heteroclinic_representation},
a front corresponds, at any sampling phase, to an orbit of the period map
from $(1,0)$ to $(0,0)$.
Its samples are common points of the unstable manifold of $(1,0)$ and the stable manifold of $(0,0)$,
and every common point of the two manifolds is a sample of a front.
This subsection describes how arcs of the unstable manifold of $(1,0)$ fold and meet the stable manifold of $(0,0)$,
which of these common points are samples of the fronts of Theorem~\ref{thm:separated_fronts},
and how the pinning interval appears in this picture.
Figures~\ref{fig:global_first_return} and~\ref{fig:global_endpoints} show the computed manifolds.

\paragraph{Setting.}
Fix $a\in(a_-,a_+)$, sample at the crossing phase $s:=\rho_P$ of $u^P$, and write $\Psi:=\Psi^s$.
We write states as $(u,v)=\mu_+\boldsymbol w_s+\mu_-\boldsymbol v_s$,
with the eigenvectors $\boldsymbol w_s=(1,r_+(s))$ and $\boldsymbol v_s=(1,-r_-(s))$ of Lemma~\ref{lem:equilibrium_tails}.
The stable eigendirection $\mu_+=0$ divides the plane into the \emph{outer side} $\mu_+<0$ and the \emph{inner side} $\mu_+>0$.
The threshold state of $u^P$ is $P_0:=a\boldsymbol v_s$, and $P_j:=\Psi^j(P_0)=\lambda^{-j}P_0$ for $j\geq0$.
The sample of $u^Q$ at $s+L$ is a point $Q_0$ of the straight stable segment
$\{\mu\boldsymbol v_s:0\leq\mu\leq a\}$, strictly between $P_1$ and $P_0$.
The straight part of the unstable manifold of $(1,0)$ from $\Psi^{-1}(P_0)$ to $P_0$ is the segment
\begin{equation}
	\label{eq:spargel_seed}
	I:=\bigl\{(1-c,-r_+(s)c):(1-a)/\lambda\leq c\leq1-a\bigr\}.
\end{equation}
Its image $\Psi(I)$ runs from $P_0$ through $Q_0$ to $P_1$.
Let $E$ be its arc from $P_0$ to $Q_0$ and $\Sigma_0$ its arc from $Q_0$ to $P_1$,
and set $\Sigma_j:=\Psi^j(\Sigma_0)$ for $j\geq0$.
We call the arcs $\Sigma_j$ \emph{spargels}.
The segment $B:=[Q_0,P_0]$ of the stable line is the \emph{base}.
Figure~\ref{fig:global_first_return}(a) shows these objects at $a=1/2$.

Below the threshold, the side of the stable eigendirection on which a state lies decides its future.
The following remark records this.
It is used to separate $\Psi(I)$ into $E$ and $\Sigma_0$
and to show that the region bounded by $E$ and the base contains no front.

\begin{remark}[Escape and return near $(0,0)$]
	\label{rem:spargel_departure}
	For a state $(u,v)$ at phase $s$ with $u<a$,
	write $(u,v)=\mu_+(1,r_+(s))+\mu_-(1,-r_-(s))$, that is,
	\begin{equation*}
		\mu_+=\frac{v+r_-(s)u}{r_+(s)+r_-(s)},\qquad \mu_-=u-\mu_+.
	\end{equation*}
	Until the solution reaches $u=a$, its state at $x\geq s$ is
	$\mu_+Y_+(x)/y_+(s)+\mu_-Y_-(x)/y_-(s)$.
	If $\mu_+<0<\mu_-$, the flux is negative,
	so the solution decreases and tends to $-\infty$ without reaching the threshold.
	If $\mu_+=0$ and $\mu_->0$, it is the stable equilibrium tail.
	If $\mu_+>0$ and $u>0$, the growing term forces arrival at $u=a$ in finite distance,
	with positive flux.
	For $\mu_-<0$ the flux is positive throughout,
	while for $\mu_-\geq0$ we have $u>0$ and $v'=u>0$ before arrival.
	While the solution remains below $a$ over a complete period,
	\eqref{eq:tail_floquet_relations} gives
	$(\mu_+,\mu_-)\mapsto(\lambda\mu_+,\lambda^{-1}\mu_-)$.
	When $\mu_+=0$ and $0<\mu_-<a$, the period samples stay on the straight stable segment,
	although the continuous trajectory between them need not lie on that line.
\end{remark}

\begin{proposition}[Front samples, spargels and the pocket]
	\label{prop:pocket}
	Let $a\in(a_-,a_+)$ and $s=\rho_P$.
	\begin{enumerate}
		\item[(a)] Every front has exactly one sample in the base.
		Fronts, up to translations by integer multiples of $L$, correspond bijectively to
		the common points of the unstable manifold of $(1,0)$ and the base.
		The endpoints $P_0$ and $Q_0$ correspond to $u^P$ and $u^Q$,
		and the points of the open base $(Q_0,P_0)$ to the nonmonotone fronts.
		\item[(b)] Every state on $I$ other than $\Psi^{-1}(P_0)$ crosses the threshold downward at a phase $\sigma\in[s,s+L)$.
		If $b(\sigma)>a$, its solution then escapes as in Remark~\ref{rem:spargel_departure},
		and after one period these states form the interior of $E$, which lies on the outer side and below the threshold.
		If $b(\sigma)<a$, its solution returns to the threshold,
		and after one period these states form the interior of $\Sigma_0$.
		\item[(c)] The arc $E$ meets the stable line only at its endpoints.
		Together with the base it bounds a region, the \emph{pocket},
		in which every state escapes.
		In particular, the pocket contains no sample of a front.
		\item[(d)] Let $u$ be a front of Theorem~\ref{thm:separated_fronts} with slot gap $n$.
		Its sample in the base lies within distance $C\lambda^{-n}$ of $P_0$ if its last letter is $P$,
		and within distance $C\lambda^{-n}$ of $Q_0$ if its last letter is $Q$.
		The constant $C$ depends only on $\alpha,L_1,L_2,a$ and $m$.
		Under the hypotheses of Corollary~\ref{cor:uniform_central_fronts},
		it can be chosen to depend only on $\alpha,L_1,L_2$.
		\item[(e)] For every odd $m\geq3$ there is $j_m$ such that,
		for every $j\geq j_m$, the spargel $\Sigma_j$ meets the open base in at least $2^{m-2}$ points.
		\item[(f)] Let $r\geq0$ be an integer, and let
		$\mathcal C_r:=\Psi^r\bigl(\{\mu\boldsymbol w_s:0\leq\mu\leq a\}\bigr)$,
		an arc of the unstable manifold of $(0,0)$ that starts at $(0,0)$ on the side $u>0$.
		Let $\vartheta:[0,1]\to\mathcal C_r$ be a homeomorphism.
		For every $\varepsilon>0$ and every sufficiently large $j$,
		the spargel $\Sigma_j$ contains an arc with a homeomorphism $\vartheta_j$ from $[0,1]$ onto it such that
		$|\vartheta_j(\tau)-\vartheta(\tau)|<\varepsilon$ for all $\tau\in[0,1]$.
	\end{enumerate}
\end{proposition}

The proof is given in Appendix~\ref{app:pocket}.
Parts (a) and (b) are the last and first crossings in the proof of Theorem~\ref{thm:all_front_existence},
seen from the period map.
Part (e) is proved by placing separated fronts of Theorem~\ref{thm:separated_fronts} on a given spargel,
and part (f) by the linear action of $\Psi$ near $(0,0)$.

\paragraph{How the spargels fold.}
By Proposition~\ref{prop:pocket}(b), the interior of $\Sigma_0$ consists of states whose solutions return to the threshold.
Near $(0,0)$, as long as the trajectory stays below the threshold,
$\Psi$ multiplies $\mu_+$ by $\lambda$ and $\mu_-$ by $\lambda^{-1}$.
The successive spargels are therefore stretched along the unstable eigendirection of $(0,0)$,
and by Proposition~\ref{prop:pocket}(f) the late spargels follow each of the arcs $\mathcal C_r$, $r\geq0$, of the unstable manifold of $(0,0)$,
as in Figure~\ref{fig:global_first_return}(b).
These arcs increase with $r$ and exhaust the branch of this manifold that starts at $(0,0)$ on the side $u>0$,
since the backward samples of each of its states converge to $(0,0)$ and therefore eventually lie on the segment $\{\mu\boldsymbol w_s:0<\mu\leq a\}$.
The two unstable manifolds never meet,
since a common state would have two different limits as $x\to-\infty$.
By the construction at threshold $1-a$, complemented by $u\mapsto1-u$,
the late arcs of the unstable manifold of $(0,0)$ follow the unstable manifold of $(1,0)$ in the same way.
A late spargel of $(1,0)$ therefore follows arcs of the unstable manifold of $(0,0)$,
which themselves follow earlier spargels of $(1,0)$ and fold with them.
These repeated passages near the two equilibria suggest increasingly numerous and elongated folds.
A fold of $\Sigma_j$, $j\geq1$, that enters the pocket meets the base once on the way in and once on the way out,
so the folding explains why the fronts carried by a spargel appear in pairs when the contacts with the base are transverse.
The next paragraph makes this precise,
and Proposition~\ref{prop:pocket}(e) shows that the number of points of $\Sigma_j$ in the open base tends to infinity as $j\to\infty$.

\paragraph{The pocket and pairs of fronts.}
A spargel meets the open base where it carries a nonmonotone front (Proposition~\ref{prop:pocket}(a)).
For $j\geq1$, both endpoints $\lambda^{-j}Q_0$ and $P_{j+1}$ of the spargel $\Sigma_j$ lie between $(0,0)$ and $P_1$,
outside the closed pocket.
The unstable manifold of $(1,0)$ does not intersect itself, so $\Sigma_j$ cannot cross the curved side $E$ of the pocket.
Hence every part of $\Sigma_j$ that lies inside the pocket is an arc
whose two ends are on the open base.
The spargel enters the pocket through the base and leaves it through the base again.
When these contacts with the base are transverse crossings, each visit of the pocket gives two fronts.
By Proposition~\ref{prop:pocket}(d), a separated front with last letter $P$ lies near the end $P_0$ of the base,
and one with last letter $Q$ near the end $Q_0$.
In the computed example of Figure~\ref{fig:global_first_return}(c), the six visits of $\Sigma_5$ are nested,
that is, along the base the two crossings of each visit lie between the two crossings of the next outer visit.
The two crossings of each visit belong to fronts whose words differ only in their last letter,
and the one ending with $P$ lies nearer to $P_0$.

\paragraph{Collisions.}
As $a$ increases to $a_+$, the point $Q_0$ moves to $P_0$ and the pocket shrinks to a point.
For each fixed $j\geq1$, the spargel $\Sigma_j$ at $a=a_+$ is a compact arc that does not contain $P_0$.
By continuity, $\Sigma_j$ no longer meets the base when $a$ is sufficiently close to $a_+$.
Its visits to the pocket must therefore end.
Figure~\ref{fig:global_many_crossings}(b) shows how this happens for one computed pair of fronts
whose words differ only in their last letter.
As $a$ increases, their last crossings merge at a coefficient interface.
For every fixed interior threshold, later spargels still meet the base
(Proposition~\ref{prop:pocket}(e)).

\paragraph{The pinning interval in this picture.}
The interior of $\Sigma_0$ comes from the forward crossing phases $\sigma$ with $b(\sigma)<a$,
and the interior of the curved side $E$ of the pocket from those with $b(\sigma)>a$.
Among the stable samples in $(P_1,P_0]$, the base consists precisely of those
whose backward crossing time $t\in(s-L,s]$ satisfies $b(t)\geq a$
(Appendix~\ref{app:pocket}).
A nonmonotone front therefore needs both a nondegenerate arc $\Sigma_0$ and a nondegenerate base.
At $a=a_-$, sampled at $s=L_1$, every first downward crossing other than that of the decreasing front escapes,
and $\Sigma_0$ degenerates to the point $Q_0=P_1$ (Figure~\ref{fig:global_endpoints}(a)).
At $a=a_+$, sampled at $s=0$, every other first downward crossing returns,
but $Q_0=P_0$, and the base and the pocket collapse to this point (Figure~\ref{fig:global_endpoints}(c)).
In both cases the samples of the single decreasing front are the only common points of the two manifolds,
as Theorem~\ref{thm:all_front_existence} asserts.

\paragraph{The balanced case.}
At $a=1/2$ the sampling phase is $s=L_1/2$, and the reflection $x\mapsto L_1-x$ preserves $p$.
If $u$ is a solution, so is $x\mapsto1-u(L_1-x)$, and its flux at $x$ equals the flux of $u$ at $L_1-x$.
Hence $R(u,v):=(1-u,v)$ satisfies $R\circ\Psi=\Psi^{-1}\circ R$
and maps the unstable manifold of $(1,0)$ onto the stable manifold of $(0,0)$.
The spargels $\Sigma_j$ and their mirror images $R(\Sigma_j)$ are therefore arcs of the two manifolds,
and their common points are fronts.
On the threshold line $u=1/2$, which $R$ fixes,
every point of a spargel is also a point of its mirror image,
so it is the sample at $s$ of a front satisfying $u(x)=1-u(L_1-x)$.
Conversely, every nonmonotone front with this symmetry gives a point of this kind.
Reflection reverses the order of the crossings and maps the letters $P_k$ and $Q_k$ to $P_{-k}$ and $Q_{-k-1}$.
For large $N$, Theorem~\ref{thm:separated_fronts} gives four fronts with words
$P_{-N}P_0P_N$, $Q_{-N}P_0Q_{N-1}$, $P_{-N}P_0Q_{N-1}$ and $Q_{-N}P_0P_N$
whose samples at $s$ lie in $\Sigma_{N-1}\cap R(\Sigma_{N-1})$.
The first two are symmetric and lie on the threshold line,
while the last two are exchanged by $R$ and lie off the line.
Figure~\ref{fig:global_first_return}(d) shows computed intersections with these four words for $N=3$,
after shifting the slot indices so that the first is zero,
and Figure~\ref{fig:global_endpoints}(b) shows the trajectories of the two symmetric fronts.
Proposition~\ref{prop:balanced_minimum_gap} determines exactly the $N$ for which the symmetric front $P_{-N}P_0P_N$ exists.
Earlier spargels can meet the threshold line at samples of symmetric fronts
whose exterior crossings are at coefficient interfaces,
so the smallest $j$ for which $\Sigma_j$ meets the line can be less than $N_*-1$.

\paragraph{The computed picture.}
In the computed example of Figure~\ref{fig:global_first_return}(a), where $a=1/2$, the interior of the computed $\Sigma_0$ lies on the inner side below the threshold,
and the arc meets the base only at its endpoint $Q_0$.
Since the McKean reaction is affine on each side of the threshold,
the computations for Figures~\ref{fig:global_first_return} and~\ref{fig:global_endpoints} involve no discretization of the equation.
The segment $I$ in \eqref{eq:spargel_seed} is an exact part of the unstable manifold of $(1,0)$.
Between threshold crossings, a state is advanced on each coefficient interval by the explicit affine flows built from \eqref{time_t} and \eqref{time_t_alpha},
and the distance $t$ to the next crossing is found in closed form from a quadratic equation in $e^{t}$ or $e^{\sqrt\alpha t}$.
The computed spargels are the images of finely and equally spaced points of $I$.
Each intersection with the base in Figure~\ref{fig:global_first_return}(c), and of $\Sigma_2$ with $R(\Sigma_2)$ in~(d), is detected between the images of neighbouring points
and then located by a root finder in the parameters $c$ of the corresponding points of $I$, using the same maps.
These are floating-point computations of finitely many arcs.
Any exact intersection of an arc with the base represents a front,
but finitely many computed images do not determine all such intersections.
For a smooth reaction, the local manifolds and the period map are in general not known in closed form and must themselves be approximated.
For a small periodic perturbation of the Allen--Cahn equation, for example, Bastiaansen, Doelman, and Kaper approximate the stable and unstable manifolds of the state near $u=-1$ to first order in the perturbation amplitude,
and they do not carry the approximation beyond the first passage near $u=1$ \cite[Section~6.2]{bastiaansenMultifrontDynamicsSpatially2026}.

\FloatBarrier
\section{Discussion}
\label{sec:discussion}

The pinning interval $[a_-,a_+]$ governs all stationary fronts of the scalar problem,
and hence all necklace fronts with equal profiles on each pair of parallel edges.
Up to translation by integer multiples of $L$, each endpoint admits one front, which is decreasing,
whereas each interior threshold admits two decreasing fronts and infinitely many nonmonotone ones.
Thus the matching function, which locates the crossing phases of the two decreasing fronts,
also identifies precisely the thresholds at which nonmonotone fronts occur.

For sufficiently separated crossings, the two isolated phases give independent choices,
and the interaction of the crossings moves them by exponentially small amounts.
The resulting count is exhaustive within the prescribed slots.
Near an endpoint, the interaction can bring a crossing to a coefficient interface.
Figure~\ref{fig:global_many_crossings}(b) shows two computed fronts whose words differ only in their last letter
meeting in this way and disappearing as $a$ increases.
At the same threshold, Theorem~\ref{thm:separated_fronts} still gives fronts with both words
in more widely separated slots.
Which pairs meet, and at which thresholds, is not analyzed here.

The uniform gap bound of Corollary~\ref{cor:uniform_central_fronts}
applies to every odd crossing number
outside layers of width $\mathcal O(\lambda^{-N})$ at the endpoints of the pinning interval.
It does not classify fronts with shorter gaps.
Even at $a=1/2$, Proposition~\ref{prop:balanced_minimum_gap}
gives a nontrivial minimum gap for the reflection-symmetric
family $P_{-N}P_0P_N$, and this minimum diverges as the diffusion contrast disappears.
Which other words occur at short gaps remains open.
The repeated returns of the computed unstable-manifold arcs to the stable segment
in Figure~\ref{fig:global_first_return}
also raise the question of the structure of the set of all bounded spatial trajectories,
for example whether it admits a Cantor-set description.
Computing finitely many manifold arcs cannot settle this question.

The temporal results concern the two decreasing fronts,
among mild solutions with a single threshold crossing whose location is continuous in time.
The coordinate change and the norm equivalence of Section~\ref{subsec:graph_transfer}
transfer these results to the corresponding class of necklace solutions
with equal profiles on each pair of parallel edges.
Stability of nonmonotone fronts, and stability under perturbations
that break this equality, require an evolution theory and comparison
or spectral arguments for those larger classes.
For a small, smooth, periodic topographic perturbation of the Allen--Cahn equation,
Bastiaansen, Doelman, and Kaper \cite{bastiaansenMultifrontDynamicsSpatially2026}
construct multi-front patterns with sufficiently separated fronts
near simple zeros of a Melnikov function, which give the single-front pinning positions to leading order.
The derivatives at the chosen zeros determine the leading terms of the small eigenvalues,
and patterns assembled from stable single-front pinning positions are spectrally stable for sufficiently large gaps.
This suggests that the separated fronts of Theorem~\ref{thm:separated_fronts}
whose words contain only the letter $P$ are stable.
Two features of our construction support this expectation.
Each crossing selected by the letter $P$ lies where $b$ is decreasing.
This is the sign used in Section~\ref{subsec:stability} to build the comparison profiles that make $u^P$ stable
(see \eqref{eq:g_sign}), and the letter $Q$ gives the opposite sign.
Moreover, the interactions between the crossings are exponentially small in the slot gap,
so for sufficiently large gaps they do not reverse the signs supplied by the isolated crossings.
A proof along the lines of Section~\ref{sec:stability} would need local existence, uniqueness and comparison
for mild solutions with several continuously moving threshold crossings,
and comparison profiles that control all the crossings at once.

\paragraph{Code availability.}
The code that computes Figures~\ref{fig:monotone_summary}--\ref{fig:global_endpoints}
is available at \url{https://github.com/mathmd/necklace-figures}.

\paragraph{Use of large language models.}
The authors used large language models (Claude, Anthropic; GPT models, OpenAI)
to improve the wording and organization of the manuscript,
to assist in writing the code that produces the figures,
and to check the consistency of arguments and cross-references.
The authors reviewed all suggested text and code
and take full responsibility for the content of the paper.

\appendix
\section{The McKean switching flow and local intersections}
\label{app:mckean}

\subsection{The switching flow}
\label{app:mckean_switching}

\begin{proof}[Proof of Lemma~\ref{lem:mckean_switching}]
	Fix $\kappa\in\{1,\alpha\}$
	and write $Z_a:=(a,0)$.

	Away from $u=a$,
	local existence and uniqueness follow from the linear equations.
	At a point $(a,v_0)$ with $v_0\neq0$,
	the identity $u'=\kappa v$ and continuity of $v$
	determine the crossing direction.
	For increasing $x$,
	the trajectory passes into $u>a$ if $v_0>0$
	and into $u<a$ if $v_0<0$.
	Concatenating the corresponding linear solutions
	gives a unique local continuation in both time directions.
	Thus local existence and uniqueness hold throughout $\Omega_a$.

	The primitive $F_a(u):=\int_0^uf_{\mathrm{MK}}(w,a)\,dw$ of the McKean reaction is
	\begin{equation*}
		F_a(u)=-\frac{u^2}{2}+\max\{u-a,0\},
		\qquad
		H_\kappa(u,v):=\frac{\kappa v^2}{2}+F_a(u).
	\end{equation*}
	On each segment contained in $u<a$ or $u>a$,
	differentiation gives
	\begin{equation*}
		\frac{d}{dx}H_\kappa(u,v)
		=\kappa vv'+f_{\mathrm{MK}}(u,a)u'
		=0.
	\end{equation*}
	At a threshold crossing with $v\neq0$,
	continuity of the trajectory and of $H_\kappa$
	ensures that the conserved values on the adjacent segments agree.
	Consequently,
	for any solution $Z:I\to\mathbb{R}^2$,
	the function $H_\kappa(Z(x))$ is locally constant wherever $Z(x)\neq Z_a$,
	and hence constant on each connected component of
	$\{x\in I:Z(x)\neq Z_a\}$.

	Since $F_a'(a^-)=-a<0<F_a'(a^+)=1-a$,
	$Z_a$ is a strict local minimum of $H_\kappa$.
	Thus $Z_a$ is isolated in the level set
	$\{H_\kappa=H_\kappa(Z_a)\}$.
	On any connected interval where a trajectory avoids $Z_a$,
	conservation and continuity imply that approach to $Z_a$
	would force its energy to equal $H_\kappa(Z_a)$,
	contradicting this isolation.
	A solution attaining $Z_a$ must therefore be identically $Z_a$,
	but the constant curve does not satisfy
	$v'=-f_{\mathrm{MK}}(a,a)=a-1\neq0$.
	Therefore, no solution passes through $Z_a$,
	and every threshold crossing is transverse.

	For global existence, note that the vector field satisfies
	$|Z'|\leq C_\kappa(1+|Z|)$ almost everywhere
	for some $C_\kappa>0$.
	If a maximal trajectory had a finite right endpoint $T$,
	Gr\"onwall's inequality would give boundedness of $Z$ as $x\to T^-$.
	The preceding bound then gives uniform boundedness of $Z'$,
	so $Q:=\lim_{x\to T^-}Z(x)$ exists.
	The energy argument excludes $Q=Z_a$.
	Local existence at $Q\in\Omega_a$
	would therefore extend the trajectory beyond $T$,
	contradicting maximality.
	The same argument applies at a finite left endpoint.
	Hence every trajectory starting in $\Omega_a$
	exists for all positive and negative times.

	Since $u'=\kappa v$ with $v$ continuous,
	an accumulation of threshold crossings at a finite time $t_*$
	would imply $u(t_*)=a$ and $v(t_*)=0$,
	contradicting the exclusion of $Z_a$.
	Thus every compact time interval contains only finitely many crossings.
	It remains to prove continuous dependence through one crossing.

	Let $\varphi_{\mathrm{in}}$ and $\varphi_{\mathrm{out}}$
	be the affine flows before and after the crossing,
	extended to $\mathbb{R}^2$.
	Suppose a trajectory starting off the threshold at $z_0$
	first crosses at time $t_0>0$.
	At $(t_0,z_0)$, the equation
	\begin{equation*}
		\bigl(\varphi_{\mathrm{in}}(t,z)\bigr)_1=a
	\end{equation*}
	has nonzero time derivative $\kappa v(t_0)$.
	The implicit function theorem therefore gives a continuous crossing time
	$t=\tau(z)$ near $t_0$.
	For a fixed small $\varepsilon>0$,
	separation of the reference trajectory from the threshold on $[0,t_0-\varepsilon]$
	excludes earlier crossings for nearby initial points.
	The crossing state
	$Q(z):=\varphi_{\mathrm{in}}(\tau(z),z)$ is continuous as well.
	After the crossing, the trajectory is
	$\varphi_{\mathrm{out}}(t-\tau(z),Q(z))$,
	which is jointly continuous in $(t,z)$ and agrees with the incoming trajectory
	at $t=\tau(z)$.

	For an initial point $z_0=(a,v_0)$ with $v_0\neq0$,
	apply the same implicit function theorem to the affine flow on each side,
	allowing positive and negative times near zero.
	The resulting local crossing time $\tau(z)$ is positive on the incoming side
	and negative on the outgoing side.
	For points on the threshold, set $\tau(z)=0$ and $Q(z)=z$.
	The one-sided constructions have the common limits
	$\tau(z)\to0$ and $Q(z)\to z_0$,
	so these definitions are continuous near $z_0$.
	For $(t,z)$ near $(0,z_0)$, the switched trajectory is
	\begin{equation*}
		\Phi(t,z)=
		\begin{cases}
			\varphi_{\mathrm{in}}\bigl(t-\tau(z),Q(z)\bigr),
			 & t\leq\tau(z), \\
			\varphi_{\mathrm{out}}\bigl(t-\tau(z),Q(z)\bigr),
			 & t\geq\tau(z).
		\end{cases}
	\end{equation*}
	Both expressions are jointly continuous and agree when $t=\tau(z)$.
	Between crossings, continuity follows from the affine flows.
	Composing across the finitely many crossings on a compact time interval
	proves joint continuity of the switched flow.
	Backward evolution is treated in the same way,
	and autonomy and uniqueness give the flow property.
\end{proof}

\subsection{Local invariant-manifold intersections}
\label{app:mckean_geometry}

\begin{proof}[Proof of Lemma~\ref{lem:local_branch_ordering}]
	Choose $j\in\{1,2\}$ according to the coefficient interval containing $s_0$.
	For $\theta$ close to $s_0$ in the same interval,
	write
	\begin{equation*}
		T^\eta(\theta):=\bigl(a,v^{\mathrm \eta}(\theta,a)\bigr)^{\mathsf T},
		\qquad \eta\in\{u,s\},
	\end{equation*}
	with the threshold fluxes \eqref{eq:threshold_flux_asymmetric}.
	For $\eta\in\{u,s\}$,
	let $Z^\eta_\theta$ be the solution of the switching system with
	$Z^\eta_\theta(\theta)=T^\eta(\theta)$,
	obtained by extending the corresponding equilibrium tail.
	We evaluate these solutions at the fixed spatial point $x=s_0$.
	Since $v_*<0$,
	for $\theta$ sufficiently close to $s_0$ we have
	\begin{equation*}
		\gamma^\eta(\theta):=Z^\eta_\theta(s_0)=
		\begin{cases}
			\phi_l^j(s_0-\theta)T^\eta(\theta), & \theta\leq s_0, \\
			\phi_r^j(s_0-\theta)T^\eta(\theta), & \theta\geq s_0.
		\end{cases}
	\end{equation*}
	The solution $Z^u_\theta$ approaches $(1,0)$ as $x\to-\infty$,
	whereas $Z^s_\theta$ approaches $(0,0)$ as $x\to+\infty$.
	Hence $\gamma^u(\theta)$ and $\gamma^s(\theta)$ belong, respectively,
	to the unstable and stable manifolds of these equilibria for $\Psi^{s_0}$.
	Varying $\theta$ parametrizes arcs through $(a,v_*)$,
	which corresponds to $\theta=s_0$.

	The affine flows and the threshold fluxes are $C^1$
	on each side of $\theta=s_0$,
	and differentiation gives
	\begin{equation*}
		(\gamma_1^\eta)'(s_0^-)
		=(\gamma_1^\eta)'(s_0^+)
		=-\kappa v_*>0.
	\end{equation*}
	Their first components are therefore strictly increasing near $s_0$.
	Using these components as the coordinate $\xi$
	gives the stated graphs on a common neighborhood of $a$.
	Since the second components of the affine vector fields
	at $(a,v_*)$ are $a$ and $a-1$,
	the one-sided graph derivatives are
	\begin{equation*}
		(h^\eta)'(a^-)
		=\frac{\partial_sv^{\mathrm \eta}(s_0,a)-a}{-\kappa v_*},
		\qquad
		(h^\eta)'(a^+)
		=\frac{\partial_sv^{\mathrm \eta}(s_0,a)-a+1}{-\kappa v_*}.
	\end{equation*}
	Subtracting the unstable derivative from the stable derivative
	gives the same value on both sides of $a$,
	namely the first expression in \eqref{eq:local_branch_gap}.
	The difference $h^s-h^u$ is therefore $C^1$ across $a$.

	Since $b(s_0)=a$, differentiating \eqref{eq:flux_mismatch} at fixed $a$ gives
	\begin{equation*}
		\partial_s(v^{\mathrm u}-v^{\mathrm s})(s_0,a)
		=-\bigl(r_+(s_0)+r_-(s_0)\bigr)b'(s_0),
	\end{equation*}
	which is the second expression in \eqref{eq:local_branch_gap}.
	By Lemma~\ref{lem:matching_function}, this is positive at $\rho_P\in(0,L_1)$
	and negative at $\rho_Q\in(L_1,L)$.
	Since $\kappa v_*<0$,
	the derivative in \eqref{eq:local_branch_gap}
	is negative at $\rho_P$ and positive at $\rho_Q$.
	As $h^s(a)=h^u(a)$,
	this proves \eqref{eq:local_branch_ordering}.

	Finally, restrict $|\theta-s_0|<L$,
	so that $s_0-L<\theta<s_0+L$.
	The states
	$\Psi^{s_0}\gamma^s(\theta)=Z^s_\theta(s_0+L)$
	and
	$(\Psi^{s_0})^{-1}\gamma^u(\theta)=Z^u_\theta(s_0-L)$
	therefore determine, respectively,
	a stable tail lying in $u<a$ and an unstable tail lying in $u>a$.
	By periodicity and Lemma~\ref{lem:equilibrium_tails},
	these states belong to the corresponding straight portions
	of the invariant manifolds of $\Psi^{s_0}$.
	This proves the iteration assertion.
\end{proof}

\begin{proof}[Proof of Lemma~\ref{lem:boundary_contact}]
	Write
	\begin{equation*}
		\kappa_-:=p(s_0^-)^{-1},
		\qquad
		\kappa_+:=p(s_0^+)^{-1},
		\qquad
		\Delta(s,a):=v^{\mathrm u}(s,a)-v^{\mathrm s}(s,a).
	\end{equation*}
	At $s_0=0$,
	use the periodic extensions of the threshold fluxes for negative phases.
	For $\eta\in\{u,s\}$,
	construct the solutions $Z^\eta_\theta$ of the switching system
	and the curves $\gamma^\eta(\theta):=Z^\eta_\theta(s_0)$
	as in the proof of Lemma~\ref{lem:local_branch_ordering},
	using the coefficient interval to the left of $s_0$ when $\theta<s_0$
	and the one to its right when $\theta>s_0$.
	The one-sided derivatives of $\gamma^\eta_1$ at $\theta=s_0$
	are now $-\kappa_-v_*>0$ and $-\kappa_+v_*>0$.
	Thus these curves parametrize the stated graphs.
	Subtracting their one-sided graph derivatives gives
	\begin{equation}
		\label{eq:boundary_branch_derivatives}
		(h^s-h^u)'(a^-)
		=\frac{\partial_s\Delta(s_0^-,a)}{\kappa_-v_*},
		\qquad
		(h^s-h^u)'(a^+)
		=\frac{\partial_s\Delta(s_0^+,a)}{\kappa_+v_*}.
	\end{equation}

	Since $b(s_0)=a$, the identity \eqref{eq:flux_mismatch} gives
	\begin{equation*}
		\partial_s\Delta(s_0^\pm,a)
		=-\bigl(r_+(s_0)+r_-(s_0)\bigr)b'(s_0^\pm).
	\end{equation*}
	By Lemma~\ref{lem:matching_function}, $b'$ is negative to the left of $L_1$ and to the right of $0$,
	and positive to the right of $L_1$ and to the left of $0$,
	with nonzero one-sided limits.
	Thus $\partial_s\Delta(s_0^\pm,a)$ is positive on the side where $\kappa_\pm=1$
	and negative on the side where $\kappa_\pm=\alpha$.

	For $a=a_-$,
	we have $s_0=L_1$ and $(\kappa_-,\kappa_+)=(1,\alpha)$.
	Since $v_*<0$,
	\eqref{eq:boundary_branch_derivatives} gives
	$(h^s-h^u)'(a^-)<0<(h^s-h^u)'(a^+)$.
	For $a=a_+$,
	we have $s_0=0$ and $(\kappa_-,\kappa_+)=(\alpha,1)$,
	so both signs are reversed.
	As $h^s(a)=h^u(a)$,
	the asserted inequalities follow.
\end{proof}

\section{Proofs for the local McKean evolution}
\label{app:temporal}

We prove the semigroup estimates, single-crossing lemma, positivity preservation,
and comparison principle stated in Section~\ref{sec:stability}.
The notation and hypotheses are those of the corresponding statements.

\begin{proof}[Proof of Lemma~\ref{lem:semigroup_est}]
	Kato's representation theorem and the Lax--Milgram theorem,
	applied to the bilinear form
	$\mathfrak k[u,v]:=\int_{\mathbb{R}}pu_xv_x + \int_{\mathbb{R}}uv$
	on $H^1(\mathbb{R})$,
	give $K$ as a self-adjoint operator satisfying $K\geq1$ on $L^2(\mathbb{R})$,
	with resolvents $(\mu+K)^{-1}$ satisfying, for $\mu>-1$,
	\begin{equation}
		\left\|(\mu+K)^{-1}\right\|_{2\to 2} \leq (1+\mu)^{-1}.
	\end{equation}
	The Hille--Yosida theorem implies that $-K$ generates a strongly continuous semigroup $\mathsf S(t)$
	on $L^2(\mathbb{R})$ with
	\begin{equation}
		\|\mathsf S(t)\|_{2\to 2}\leq e^{-t}.
	\end{equation}
	Write $\mathsf S(t)=e^{-t}\mathcal T(t)$,
	where $\mathcal T(t)=e^{-tK_0}$ and $K_0:=K-1$ on $\mathcal D(K_0)=\mathcal D(K)$.

	\textit{Step 1: extensions and contractivity.}
	For \eqref{eq:S_contract}, we show that
	$\mathcal{T}(t)$ extends to a contraction on $L^{\infty}(\mathbb{R})$.
	Let $f\in L^2(\mathbb{R})\cap L^\infty (\mathbb{R})$ and
	$u = (1+\tau K_0)^{-1}f$ for $\tau>0$.
	Equivalently,
	\begin{equation}
		\int u\varphi + \tau\int pu_x \varphi_x  = \int f\varphi,\quad
		\forall \varphi\in H^1(\mathbb{R}).
	\end{equation}
	Set $\theta:=\|f\|_\infty$ and choose
	$\varphi:=(u-\theta)^+\in H^1(\mathbb{R})$ in the weak equation.
	Since $u=\theta+\varphi$ wherever $\varphi>0$, we obtain
	\begin{equation*}
		\|\varphi\|_2^2
		+\tau\int_{\mathbb{R}}p|\varphi_x|^2\,dx
		=\int_{\mathbb{R}}(f-\theta)\varphi\,dx\leq0.
	\end{equation*}
	Both terms on the left are nonnegative, so $\varphi=0$.
	Thus $(1+\tau K_0)^{-1}f\leq\|f\|_\infty$ almost everywhere.
	The same argument applied to $-u$ yields
	$(1+\tau K_0)^{-1}f\geq-\|f\|_\infty$.
	Therefore,
	\begin{equation}\label{eq:approx_Linf_bdd}
		\left|(1+\tau K_0)^{-1}f\right|\leq \|f\|_\infty,\quad
		\forall f\in L^2(\mathbb{R})\cap L^\infty(\mathbb{R}).
	\end{equation}
	The resolvent approximation gives
	\begin{equation*}
		\mathcal T(t)f
		=\lim_{n\to\infty}\left(1+\frac{t}{n}K_0\right)^{-n}f
		\quad\text{in }L^2(\mathbb{R}).
	\end{equation*}
	Applying \eqref{eq:approx_Linf_bdd} successively shows that
	each approximant has absolute value at most $\|f\|_\infty$ almost everywhere.
	Passing to a subsequence that converges almost everywhere gives
	\begin{equation*}
		\|\mathcal T(t)f\|_\infty\leq\|f\|_\infty,
		\qquad f\in L^2(\mathbb{R})\cap L^\infty(\mathbb{R}).
	\end{equation*}
	For $g\in L^1(\mathbb{R})\cap L^2(\mathbb{R})$, we have
	\begin{equation}\label{eq:T1_subest}
		\left\|\mathcal{T}(t)g\right\|_1
		= \sup_{\substack{\|f\|_\infty\leq 1\\f\in L^2(\mathbb{R})\cap L^\infty(\mathbb{R})}} \left|\int  \mathcal{T}(t)gf\right|
		= \sup_{\substack{\|f\|_\infty\leq 1\\f\in L^2(\mathbb{R})\cap L^\infty(\mathbb{R})}} \left|\int g \mathcal{T}(t)f\right|
		\leq \|g\|_1,
	\end{equation}
	by self-adjointness of $\mathcal{T}(t)$ on $L^2(\mathbb{R})$.
	Equation \eqref{eq:T1_subest} extends to all $g\in L^1(\mathbb{R})$ by density:
	\begin{equation}
		\left\|\mathcal{T}_1(t)g\right\|_1 \leq \|g\|_1,\quad
		\forall g\in L^1(\mathbb{R}),
	\end{equation}
	where $\mathcal{T}_1 = \mathcal{T}$ in $L^2(\mathbb{R})\cap L^1(\mathbb{R})$.
	By duality, $\mathcal{T}_\infty(t):=\mathcal{T}_1(t)^*$
	coincides with $\mathcal{T}(t)$ on $L^2(\mathbb{R})\cap L^\infty(\mathbb{R})$
	and is its unique weak-$*$ continuous extension to $L^\infty(\mathbb{R})$.
	Moreover,
	\begin{equation}
		\left\|\mathcal{T}_\infty(t)\right\|_{\infty \to\infty}
		= \left\|\mathcal{T}_1(t)\right\|_{1\to 1} \leq 1.
	\end{equation}
	This proves \eqref{eq:S_contract}.

	To prove strong continuity on $L^1(\mathbb{R})$,
	let $f\in L^2(\mathbb R)$ be supported in a bounded interval $I_0$.
	Then $f\in L^1(\mathbb{R})$, and contractivity on $L^1(\mathbb{R})$ gives
	\begin{align*}
		\int_{\mathbb R\setminus I_0}|\mathcal T_1(t)f|\,dx
		 & \leq\int_{I_0}\bigl(|f|-|\mathcal T_1(t)f|\bigr)\,dx \\
		 & \leq\int_{I_0}|\mathcal T_1(t)f-f|\,dx.
	\end{align*}
	Consequently,
	\begin{equation*}
		\|\mathcal T_1(t)f-f\|_1
		\leq2|I_0|^{1/2}\|\mathcal T(t)f-f\|_2
		\to 0
		\qquad(t\downarrow0).
	\end{equation*}
	Compactly supported $L^2(\mathbb{R})$ functions are dense in $L^1(\mathbb{R})$.
	Contractivity therefore extends this convergence to every $f\in L^1(\mathbb{R})$,
	and the semigroup property gives continuity at every $t\geq0$.
	The $L^1(\mathbb{R})$ extension of $\mathsf S(t)=e^{-t}\mathcal T(t)$
	is consequently strongly continuous as well.

	Set $\mathsf S_1(t):=e^{-t}\mathcal T_1(t)$
	and $\mathsf S_\infty(t):=e^{-t}\mathcal T_\infty(t)=\mathsf S_1(t)^*$.
	For $g\in L^\infty(\mathbb R)$, $\varphi\in L^1(\mathbb R)$,
	and $r\geq0$, duality gives
	\begin{equation*}
		\left|\int_{\mathbb R}
		\bigl(\mathsf S_\infty(t)g-\mathsf S_\infty(r)g\bigr)\varphi\,dx\right|
		\leq\|g\|_\infty\|\mathsf S_1(t)\varphi-\mathsf S_1(r)\varphi\|_1
		\longrightarrow0
		\qquad(t\to r).
	\end{equation*}
	This proves the asserted weak-$*$ continuity in time.

	\textit{Step 2: the flux estimate.}
	For \eqref{eq:DxS_contract}, set
	\begin{equation*}
		H_0:=L^2(\mathbb R,dx),
		\qquad
		H_1:=L^2(\mathbb R,p^{-1}\,dx),
	\end{equation*}
	and let $D:H^1(\mathbb R)\subset H_0\to H_1$ be given by $Du=pu_x$.
	Since $p$ and $p^{-1}$ are bounded, $D$ is closed and densely defined.
	Integration by parts gives
	\begin{equation*}
		\langle Du,v\rangle_{H_1}=-\langle u,v_x\rangle_{H_0},
		\qquad u,v\in H^1(\mathbb R).
	\end{equation*}
	Thus $\mathcal D(D^*)=H^1(\mathbb R)$ and $D^*v=-v_x$.
	The operators $K_0=D^*D$ and $B:=DD^*$
	are nonnegative and self-adjoint on $H_0$ and $H_1$, respectively.

	Fix $\mu>0$ and $u\in\mathcal D(D)$,
	and set $w:=(\mu+K_0)^{-1}u$.
	Then $w\in\mathcal D(K_0)\subset\mathcal D(D)$
	and $K_0w=u-\mu w\in\mathcal D(D)$.
	Hence $Dw\in\mathcal D(B)$, and
	\begin{equation*}
		(\mu+B)Dw=D(\mu w+K_0w)=Du.
	\end{equation*}
	Iterating the resulting resolvent identity yields
	\begin{equation*}
		D\left(1+\frac{t}{n}K_0\right)^{-n}u
		=\left(1+\frac{t}{n}B\right)^{-n}Du.
	\end{equation*}
	The resolvent approximations converge in $H_0$ and $H_1$, respectively.
	Closedness of $D$ therefore gives
	\begin{equation}
		\label{eq:intertwining_id}
		De^{-tK_0}u=e^{-tB}Du,
		\qquad u\in H^1(\mathbb R).
	\end{equation}

	The quadratic form of $B$ is
	\begin{equation*}
		\mathfrak b[v,v]=\|D^*v\|_{H_0}^2
		=\int_{\mathbb R}|v_x|^2\,dx,
		\qquad \mathcal D(\mathfrak b)=H^1(\mathbb R).
	\end{equation*}
	Applying the truncation and duality argument of Step 1
	with the measure $p^{-1}\,dx$
	gives a weak-$*$ continuous contraction extension of $e^{-tB}$ to $L^\infty(\mathbb{R})$.

	To extend \eqref{eq:intertwining_id} to $g\in W^{1,\infty}(\mathbb R)$,
	choose $\chi\in C_c^\infty(\mathbb R)$ with $0\leq\chi\leq1$
	and $\chi=1$ near zero, and set $g_R(x):=\chi(x/R)g(x)$.
	Then $g_R\in H^1(\mathbb R)$,
	while $g_R\to g$ and $Dg_R\to pg_x$ weak-$*$ in $L^\infty(\mathbb{R})$.
	Using the consistent $L^\infty(\mathbb{R})$ extensions of both semigroups,
	pass to the limit in \eqref{eq:intertwining_id} in distributions to obtain
	\begin{equation*}
		\partial_x e^{-tK_0}g=p^{-1}e^{-tB}(pg_x).
	\end{equation*}
	Since $\mathsf S(t)=e^{-t}e^{-tK_0}$,
	the contraction estimate for $e^{-tB}$ gives
	\begin{equation*}
		\|p\partial_x \mathsf S(t)g\|_\infty
		\leq e^{-t}\|pg_x\|_\infty,
	\end{equation*}
	which proves \eqref{eq:DxS_contract}.

	\textit{Step 3: the smoothing estimate.}
	To prove \eqref{eq:S_bdd}, we first estimate $\mathcal T(t)$
	from $L^1(\mathbb{R})$ to $L^2(\mathbb{R})$.
	Set $p_{\min}:=\|p^{-1}\|_\infty^{-1}>0$.
	For $0\ne g\in L^1(\mathbb{R})\cap L^2(\mathbb{R})$,
	let $u(t):=\mathcal T(t)g$.
	For $t>0$, the energy identity gives
	\begin{equation*}
		\frac12\frac{d}{dt}\|u(t)\|_2^2
		=-\int_{\mathbb{R}}p|u_x(t)|^2\,dx
		\leq-p_{\min}\|u_x(t)\|_2^2.
	\end{equation*}
	The one-dimensional Nash inequality
	\begin{equation*}
		\|v\|_2^6\leq C_{\mathrm N}\|v_x\|_2^2\|v\|_1^4,
		\qquad v\in H^1(\mathbb{R})\cap L^1(\mathbb{R}),
	\end{equation*}
	and $\|u(t)\|_1\leq\|g\|_1$ imply that $y(t):=\|u(t)\|_2^2$ satisfies
	\begin{equation*}
		y'(t)\leq-\frac{2p_{\min}}{C_{\mathrm N}\|g\|_1^4}y(t)^3.
	\end{equation*}
	Integrating this inequality and using continuity at $t=0$
	in $L^2(\mathbb{R})$ gives
	\begin{equation*}
		y(t)\leq
		\left(y(0)^{-2}+\frac{4p_{\min}t}{C_{\mathrm N}\|g\|_1^4}\right)^{-1/2}
		\leq\left(\frac{C_{\mathrm N}}{4p_{\min}t}\right)^{1/2}\|g\|_1^2.
	\end{equation*}
	Consequently, with $C_*:=(C_{\mathrm N}/(4p_{\min}))^{1/4}$,
	\begin{equation*}
		\|\mathcal T(t)g\|_2\leq C_*t^{-1/4}\|g\|_1.
	\end{equation*}
	The estimate also holds for $g=0$ and extends by density
	to $\mathcal T_1(t)$ on $L^1(\mathbb{R})$.
	Duality and self-adjointness of $\mathcal T(t)$ give
	\begin{equation*}
		\|\mathcal T(t)\|_{2\to\infty}
		=\|\mathcal T_1(t)\|_{1\to2}
		\leq C_*t^{-1/4}.
	\end{equation*}
	For $g\in L^1(\mathbb{R})$, the semigroup property therefore yields
	\begin{align*}
		\|\mathcal T_1(t)g\|_\infty
		 & \leq\|\mathcal T(t/2)\|_{2\to\infty}
		\|\mathcal T_1(t/2)g\|_2                \\
		 & \leq\sqrt{2}\,C_*^2t^{-1/2}\|g\|_1.
	\end{align*}
	Since $\mathsf S(t)=e^{-t}\mathcal T(t)$ on the consistent realizations,
	this proves \eqref{eq:S_bdd} with $C_0=\sqrt{2}\,C_*^2$.

	\textit{Step 4: the derivative estimate.}
	Let $g\in L^2(\mathbb R)$ and $t>0$,
	and set $u:=\mathcal T(t)g$ and $q:=pu_x$.
	Since $K_0$ is nonnegative and self-adjoint,
	and $\mathcal T(t)=e^{-tK_0}$, spectral calculus gives
	\begin{equation*}
		\|K_0^\sigma\mathcal T(t)\|_{2\to2}
		\leq C_\sigma t^{-\sigma},
		\qquad \sigma\in\{1/2,1\}.
	\end{equation*}
	In particular, $u\in\mathcal D(K_0)$ and $q\in H^1(\mathbb R)$.
	The form identity $\|K_0^{1/2}u\|_2^2=\int_{\mathbb R}p|u_x|^2\,dx$
	and $q_x=-K_0u$ imply
	\begin{align*}
		\|q\|_2
		 & \leq \|p\|_\infty^{1/2}\|K_0^{1/2}\mathcal T(t)g\|_2
		\leq Ct^{-1/2}\|g\|_2,                                \\
		\|q_x\|_2
		 & =\|K_0\mathcal T(t)g\|_2
		\leq Ct^{-1}\|g\|_2.
	\end{align*}
	The one-dimensional interpolation inequality
	$\|q\|_\infty^2\leq2\|q\|_2\|q_x\|_2$
	therefore gives $\|q\|_\infty\leq Ct^{-3/4}\|g\|_2$.
	Since $\partial_x \mathsf S(t)g=e^{-t}p^{-1}q$, we obtain
	\begin{equation*}
		\|\partial_x \mathsf S(t)g\|_\infty
		\leq C_1 e^{-t}t^{-3/4}\|g\|_2,
	\end{equation*}
	which proves \eqref{eq:DxS_bdd}.

	Finally, for the forcing in the mild formulation \eqref{eq:mild_w},
	\begin{equation*}
		\|h_\rho(r)\|_1=|\rho(r)-\rho_*|,
		\qquad
		\|h_\rho(r)\|_2=|\rho(r)-\rho_*|^{1/2}.
	\end{equation*}
	Applying \eqref{eq:S_contract} to the initial-data term
	and \eqref{eq:S_bdd} to the time integral gives \eqref{eq:w_bdd}.
	For the derivative, \eqref{eq:DxS_contract} gives
	\begin{equation*}
		\|\partial_x\mathsf S(t)w_0\|_\infty
		\leq e^{-t}\|p^{-1}\|_\infty\|p\partial_xw_0\|_\infty.
	\end{equation*}
	Apply \eqref{eq:DxS_bdd} to the time integral.
	Differentiation under this integral is justified in distributions
	by boundedness of $h_\rho$ in $L^2(\mathbb{R})$
	and integrability of $(t-r)^{-3/4}$.
	Using $|\rho(r)-\rho_*|\leq d_\rho(t)$ for $0\leq r\leq t$
	and changing variables in the time integral gives \eqref{eq:wx_bdd}.
\end{proof}

\begin{proof}[Proof of Lemma~\ref{lem:single_crossing}]
	Set $J:=[\rho_*-R,\rho_*+R]$ and $m:=M-M_*>0$.
	Fix $\rho\in X_{R,T}$.
	Since $d_\rho(t)\leq R$ and $I_\gamma(t)\leq I_\gamma(T)$,
	\eqref{eq:w_bdd}--\eqref{eq:wx_bdd} give
	\begin{align*}
		\|u^\rho(\cdot,t)-u_*\|_\infty
		 & \leq\|w_0\|_\infty+C_0R I_{1/2}(T)<\delta,     \\
		\|\partial_x(u^\rho(\cdot,t)-u_*)\|_\infty
		 & \leq\|p^{-1}\|_\infty\|p\partial_xw_0\|_\infty
		+C_1R^{1/2}I_{3/4}(T)\leq M_*
	\end{align*}
	for every $t\in[0,T]$.
	Combining these estimates with \eqref{eq:R_restriction} yields
	\begin{equation*}
		\begin{aligned}
			u^\rho(x,t)           & >a+\delta     &                         & (x\leq\rho_*-R), \\
			u^\rho(x,t)           & <a-\delta     &                         & (x\geq\rho_*+R), \\
			\partial_xu^\rho(x,t) & \leq-M+M_*=-m
			                      &               & \text{for a.e. }x\in J.
		\end{aligned}
	\end{equation*}
	For $x<y$ in $J$, integration of the last inequality gives
	\begin{equation}
		\label{eq:trial_spatial_separation}
		u^\rho(x,t)-u^\rho(y,t)
		=-\int_x^y\partial_z u^\rho(z,t)\,dz
		\geq m(y-x).
	\end{equation}
	Thus the continuous representative of $u^\rho(\cdot,t)$
	has a unique threshold crossing $\widetilde\rho(t)$ in the interior of $J$.

	We next prove that $u^\rho\in C([0,T];C(J))$.
	For $\varphi\in C_c^\infty(\mathbb R)$,
	the duality construction of the $L^\infty(\mathbb{R})$ semigroup gives
	\begin{equation*}
		\int_{\mathbb R}\mathsf S(t)w_0\,\varphi\,dx
		=\int_{\mathbb R}w_0\,\mathsf S(t)\varphi\,dx.
	\end{equation*}
	The right-hand side is continuous in $t$,
	since $\mathsf S(t)\varphi$ is continuous in $L^1(\mathbb{R})$
	by the proof of Lemma~\ref{lem:semigroup_est}.
	For the other term in \eqref{eq:mild_w}, write
	\begin{equation*}
		z(t):=\int_0^t \mathsf S(t-r)h_\rho(r)\,dr.
	\end{equation*}
	Regard the integrand as zero for $r>t$.
	When $t_n\to t$, these extended integrands converge in $L^2(\mathbb{R})$
	for a.e. $r\in[0,T]$, by strong continuity of $\mathsf S(t)$ on $L^2(\mathbb{R})$.
	Their $L^2(\mathbb{R})$ norms are bounded by $R^{1/2}$.
	Dominated convergence therefore gives $z(t_n)\to z(t)$ in $L^2(\mathbb{R})$.
	It follows that
	\begin{equation*}
		t\longmapsto\int_{\mathbb R}u^\rho(x,t)\varphi(x)\,dx
	\end{equation*}
	is continuous for every $\varphi\in C_c^\infty(\mathbb R)$.

	For any sequence $t_n\to t$ in $[0,T]$,
	the uniform $W^{1,\infty}(\mathbb{R})$ bounds and the Arzel\`a--Ascoli theorem
	make $u^\rho(\cdot,t_n)|_J$ relatively compact in $C(J)$.
	The preceding convergence of integrals identifies every subsequential limit
	with $u^\rho(\cdot,t)|_J$,
	first on the interior of $J$ and then at its endpoints by continuity.
	Hence the whole sequence converges in $C(J)$,
	proving the claimed time continuity.

	Fix $t,r\in[0,T]$ and apply \eqref{eq:trial_spatial_separation}
	to the profile $u^\rho(\cdot,t)$ at the two spatial points
	$\widetilde\rho(t)$ and $\widetilde\rho(r)$.
	Since $u^\rho(\widetilde\rho(t),t)
		=u^\rho(\widetilde\rho(r),r)=a$, we obtain
	\begin{align*}
		m|\widetilde\rho(t)-\widetilde\rho(r)|
		 & \leq|u^\rho(\widetilde\rho(r),t)-a|             \\
		 & =|u^\rho(\widetilde\rho(r),t)
		-u^\rho(\widetilde\rho(r),r)|                      \\
		 & \leq\|u^\rho(\cdot,t)-u^\rho(\cdot,r)\|_{C(J)}.
	\end{align*}
	Thus $\widetilde\rho$ is continuous and belongs to $X_{R,T}$.

	It remains to prove the contraction estimate.
	Let $\eta\in X_{R,T}$ and denote the crossing path of $u^\eta$
	by $\widetilde\eta:=\Phi_{u_0}(\eta)$.
	Since both trial solutions have initial datum $u_0$,
	subtracting their mild formulas gives
	\begin{equation*}
		u^\rho(\cdot,t)-u^\eta(\cdot,t)
		=\int_0^t \mathsf S(t-r)
		\bigl(\mathbf{1}_{\{x<\rho(r)\}}-\mathbf{1}_{\{x<\eta(r)\}}\bigr)\,dr.
	\end{equation*}
	Using \eqref{eq:S_bdd} and
	$\|\mathbf{1}_{\{x<\rho(r)\}}-\mathbf{1}_{\{x<\eta(r)\}}\|_1
		=|\rho(r)-\eta(r)|$, we obtain
	\begin{equation*}
		\|u^\rho(\cdot,t)-u^\eta(\cdot,t)\|_\infty
		\leq C_0 I_{1/2}(t)\|\rho-\eta\|_{C([0,t])}.
	\end{equation*}
	Apply \eqref{eq:trial_spatial_separation} to $u^\rho(\cdot,t)$
	at $\widetilde\rho(t)$ and $\widetilde\eta(t)$.
	The equalities $u^\rho(\widetilde\rho(t),t)
		=u^\eta(\widetilde\eta(t),t)=a$ then give
	\begin{align*}
		m|\widetilde\rho(t)-\widetilde\eta(t)|
		 & \leq|u^\rho(\widetilde\eta(t),t)-a|             \\
		 & =|u^\rho(\widetilde\eta(t),t)
		-u^\eta(\widetilde\eta(t),t)|                      \\
		 & \leq\|u^\rho(\cdot,t)-u^\eta(\cdot,t)\|_\infty.
	\end{align*}
	Taking the supremum over $t\in[0,T]$ gives
	\begin{equation*}
		\|\Phi_{u_0}(\rho)-\Phi_{u_0}(\eta)\|_{C([0,T])}
		\leq\frac{C_0 I_{1/2}(T)}{M-M_*}
		\|\rho-\eta\|_{C([0,T])}.
	\end{equation*}
	Condition \eqref{eq:contract_est} therefore makes $\Phi_{u_0}$ a contraction.
\end{proof}

\begin{proof}[Proof of Lemma~\ref{lem:semigroup_positive}]
	Let $f\in L^2(\mathbb{R})$ be nonnegative,
	fix $\mu>0$, and set $v:=(\mu+K)^{-1}f$.
	The resolvent equation gives
	\begin{equation*}
		\int_{\mathbb{R}}p v_x\varphi_x\,dx
		+(\mu+1)\int_{\mathbb{R}}v\varphi\,dx
		=\int_{\mathbb{R}}f\varphi\,dx,
		\qquad \varphi\in H^1(\mathbb{R}).
	\end{equation*}
	Taking $\varphi=-v^-$, where $v^-:=\max\{-v,0\}$, yields
	\begin{equation*}
		\int_{\mathbb{R}}p|\partial_xv^-|^2\,dx
		+(\mu+1)\|v^-\|_2^2
		=-\int_{\mathbb{R}}fv^-\,dx\leq0.
	\end{equation*}
	The left-hand side is nonnegative, so $v^-=0$.
	Thus $(\mu+K)^{-1}$ preserves nonnegativity for every $\mu>0$.
	For $t>0$, the resolvent approximation therefore gives
	\begin{equation*}
		\left(1+\frac{t}{n}K\right)^{-n}f\geq0,
		\qquad
		\left(1+\frac{t}{n}K\right)^{-n}f
		\longrightarrow \mathsf S_2(t)f
		\quad\text{in }L^2(\mathbb{R}).
	\end{equation*}
	Since the set of nonnegative functions is closed in $L^2(\mathbb{R})$,
	we obtain $\mathsf S_2(t)f\geq0$.

	For nonnegative $f\in L^1(\mathbb{R})$, set
	$f_n:=\min\{f,n\}\mathbf{1}_{[-n,n]}$.
	Then $f_n\in L^1(\mathbb{R})\cap L^2(\mathbb{R})$ and $f_n\to f$
	in $L^1(\mathbb{R})$.
	Consistency and contractivity give
	$\mathsf S_1(t)f_n=\mathsf S_2(t)f_n\geq0$ and $\mathsf S_1(t)f_n\to \mathsf S_1(t)f$
	in $L^1(\mathbb{R})$.
	Hence $\mathsf S_1(t)f\geq0$.

	Finally, for nonnegative $f\in L^\infty(\mathbb{R})$
	and nonnegative $\varphi\in L^1(\mathbb{R})$,
	the adjoint identity gives
	\begin{equation*}
		\int_{\mathbb{R}}\mathsf S_\infty(t)f\,\varphi\,dx
		=\int_{\mathbb{R}}f\,\mathsf S_1(t)\varphi\,dx\geq0.
	\end{equation*}
	Thus $\mathsf S_\infty(t)f\geq0$ almost everywhere.
	The assertion at $t=0$ follows from $\mathsf S(0)=1$.
	Applying positivity to $g-f$ and using linearity
	gives the stated order preservation.
\end{proof}

\begin{proof}[Proof of Proposition~\ref{prop:loc_comparison}]
	All inequalities between paths are understood pointwise on $[0,T]$.
	For initial data $f\leq g$ and paths $\rho\leq\eta$,
	positivity of $\mathsf S(t)$ gives
	\begin{align*}
		\mathcal{M}[\eta,g](t)-\mathcal{M}[\rho,f](t)
		 & =\mathsf S(t)(g-f)            \\
		 & \quad+\int_0^t \mathsf S(t-r)
		\bigl(\mathbf{1}_{\{x<\eta(r)\}}-\mathbf{1}_{\{x<\rho(r)\}}\bigr)\,dr
		\geq0.
	\end{align*}
	The trial solutions have continuous spatial representatives,
	so this inequality holds at every spatial point.
	Set $\Phi:=\Phi_{u_0^-}$.
	By Lemma~\ref{lem:single_crossing},
	$\Phi$ is a contraction from $X_{R,T}$ into itself.
	It is also order preserving.
	Indeed, if $\rho\leq\eta$, then
	$\mathcal{M}[\rho,u_0^-](t)\leq\mathcal{M}[\eta,u_0^-](t)$.
	At the crossing of $\mathcal{M}[\eta,u_0^-](t)$,
	$\mathcal{M}[\rho,u_0^-](t)$ is at most $a$.
	Their sign properties therefore give
	$\Phi(\rho)(t)\leq\Phi(\eta)(t)$.

	The sub- and supersolution inequalities, together with $u_0^-\leq u_0^+$, give
	\begin{equation*}
		u^-(\cdot,t)\leq\mathcal{M}[\rho^-,u_0^-](t),
		\qquad
		\mathcal{M}[\rho^+,u_0^-](t)
		\leq\mathcal{M}[\rho^+,u_0^+](t)\leq u^+(\cdot,t).
	\end{equation*}
	At $x=\rho^-(t)$, $\mathcal{M}[\rho^-,u_0^-](t)$ is at least $a$,
	so its crossing lies at or to the right of $\rho^-(t)$.
	At $x=\rho^+(t)$, the trial solution with initial datum $u_0^-$
	is at most $a$, so its crossing lies at or to the left of $\rho^+(t)$.
	Hence
	\begin{equation*}
		\rho^-\leq\Phi(\rho^-),
		\qquad \Phi(\rho^+)\leq\rho^+.
	\end{equation*}
	Iterating these inequalities using order preservation yields
	\begin{equation*}
		\rho^-\leq\Phi^n(\rho^-),
		\qquad \Phi^n(\rho^+)\leq\rho^+,
		\qquad n\geq1.
	\end{equation*}
	Since $X_{R,T}$ is complete and $\Phi$ is a contraction,
	both sequences converge in $C([0,T])$ to its unique fixed point $\rho_0$.
	Passing to the limit gives $\rho^-\leq\rho_0\leq\rho^+$.
	In particular, $\rho^-\leq\rho^+$.
	Applying the order property of $\mathcal{M}$ once more, we obtain
	\begin{equation*}
		u^-(\cdot,t)
		\leq\mathcal{M}[\rho^-,u_0^-](t)
		\leq\mathcal{M}[\rho^+,u_0^+](t)
		\leq u^+(\cdot,t),
		\qquad 0\leq t\leq T.
	\end{equation*}
\end{proof}

\begin{proof}[Proof of Lemma~\ref{lem:unstable_subsuper}]
	We construct each profile by joining two stationary half-line solutions
	at $x=\rho_Q$ with a common value $\beta$ close to $a$.
	The local ordering of the invariant-manifold arcs determines the flux jump
	at the joining point, and its sign will give the mild inequality.

	Set $Z_*(x):=(u^Q(x),p(x)u_x^Q(x))^{\mathsf T}$,
	using the continuous representative of the flux.
	Let $v_Q:=(pu_x^Q)(\rho_Q)<0$.
	Use the local graphs $h^u,h^s$ through $(a,v_Q)$
	from Lemma~\ref{lem:local_branch_ordering} for $\Psi^{\rho_Q}$.
	For $\beta$ near $a$, let $Z_\beta^u=(u_\beta^u,v_\beta^u)^{\mathsf T}$
	and $Z_\beta^s=(u_\beta^s,v_\beta^s)^{\mathsf T}$ be the solutions of the switching system with
	\begin{equation*}
		Z_\beta^u(\rho_Q)=(\beta,h^u(\beta))^{\mathsf T},
		\qquad
		Z_\beta^s(\rho_Q)=(\beta,h^s(\beta))^{\mathsf T}.
	\end{equation*}
	Define the continuous profile
	\begin{equation*}
		W_\beta(x):=
		\begin{cases}
			u_\beta^u(x), & x\leq\rho_Q, \\
			u_\beta^s(x), & x\geq\rho_Q.
		\end{cases}
	\end{equation*}
	In particular, $W_\beta(\rho_Q)=\beta$ and $W_a=u^Q$.

	To control the half-line profiles as $\beta\to a$,
	set $c_-:=\rho_Q-L$ and $c_+:=\rho_Q+L$.
	Continuous dependence for the switching system gives
	\begin{equation}
		\label{eq:orbit_unif_conv}
		\sup_{[c_-,\rho_Q]}|Z_\beta^u-Z_*|
		+\sup_{[\rho_Q,c_+]}|Z_\beta^s-Z_*|
		\longrightarrow0\qquad(\beta\to a).
	\end{equation}
	This follows from Lemma~\ref{lem:mckean_switching}
	by composing the flows across the fixed coefficient intervals.
	By the last assertion of Lemma~\ref{lem:local_branch_ordering},
	$Z_\beta^u(c_-)$ and $Z_\beta^s(c_+)$ lie on the corresponding straight portions
	of the unstable and stable manifolds.
	For $\beta$ sufficiently close to $a$, their state coordinates satisfy
	$u_\beta^u(c_-)\in(a,1)$ and $u_\beta^s(c_+)\in(0,a)$.
	Lemma~\ref{lem:equilibrium_tails} therefore gives decreasing tails
	on $(-\infty,c_-]$ and $[c_+,\infty)$,
	with limits $1$ and $0$, respectively, and negative flux.

	The differences from $Z_*$ on these tails solve \eqref{eq:linear_tail_system}.
	Their states at $c_-$ and $c_+$ lie in the unstable and stable eigenspaces
	of $A_{\rho_Q}$, respectively.
	Thus their states at $c_--nL$ and $c_++nL$ are the respective endpoint differences
	each multiplied by $\lambda^{-n}$.
	The linear evolution within one period is bounded.
	Together with \eqref{eq:orbit_unif_conv}, these estimates yield
	\begin{equation}
		\label{eq:interpolant_unif_conv}
		\sup_{x\leq\rho_Q}|Z_\beta^u(x)-Z_*(x)|
		+\sup_{x\geq\rho_Q}|Z_\beta^s(x)-Z_*(x)|
		\longrightarrow0.
	\end{equation}
	Since $\partial_xW_\beta=p^{-1}v_\beta^u$ on the left
	and $\partial_xW_\beta=p^{-1}v_\beta^s$ on the right,
	we obtain $W_\beta\to u^Q$ in $W^{1,\infty}(\mathbb{R})$.
	The same eigenvalue estimates give exponential decay of
	$W_\beta-u^Q$ and its derivative at both infinities;
	in particular, $W_\beta-u^Q\in H^1(\mathbb{R})$.

	On the compact interval $[c_-,c_+]$, the flux of $u^Q$
	is continuous and strictly negative.
	By \eqref{eq:orbit_unif_conv}, the fluxes of the joined pieces
	are also negative there when $\beta$ is sufficiently close to $a$.
	Together with the negative flux on the tails, this proves that $W_\beta$
	is strictly decreasing on $\mathbb{R}$.
	Its limits $1$ and $0$ give a unique threshold crossing, denoted by $\chi_\beta$.
	Since $W_\beta(\rho_Q)=\beta$,
	\begin{equation*}
		\chi_\beta>\rho_Q\quad(\beta>a),
		\qquad
		\chi_\beta<\rho_Q\quad(\beta<a).
	\end{equation*}
	Uniform convergence to $u^Q$ and its strict crossing at $\rho_Q$
	also give $\chi_\beta\to\rho_Q$ as $\beta\to a$.
	Set
	\begin{equation*}
		\underline u_\sigma:=W_{a+\sigma},
		\qquad \overline u_\sigma:=W_{a-\sigma},
		\qquad \rho_\sigma^\pm:=\chi_{a\pm\sigma}.
	\end{equation*}

	We next establish their pointwise ordering.
	For $\zeta:=\underline u_\sigma-u^Q$, subtraction of the stationary equations
	on the two open half-lines gives
	\begin{equation*}
		K\zeta=\mathbf{1}_{(\rho_Q,\rho_\sigma^+)},
		\qquad x\ne\rho_Q,
		\qquad \zeta(\rho_Q)=\sigma,
		\qquad \zeta(\pm\infty)=0.
	\end{equation*}
	On $(-\infty,\rho_Q)$, the homogeneous equation and the boundary values give
	$\zeta(x)=\sigma y_+(x)/y_+(\rho_Q)$,
	where $y_+$ is the positive scalar component of the normalized Floquet solution
	in \eqref{eq:normalized_floquet_solutions}, which decays at $-\infty$.
	Hence $\zeta>0$ on this half-line.
	Let $y$ be the stable equilibrium tail with value $\sigma$ at $\rho_Q$.
	Then it is a positive solution of
	$Ky=0$ on $(\rho_Q,\infty)$ with limit zero at infinity.
	The half-line weak maximum principle yields $\zeta-y\geq 0$.
	Positivity of $y$ gives strict inequality $\zeta>0$ everywhere.
	Likewise, for $\zeta:=u^Q-\overline u_\sigma$,
	\begin{equation*}
		K\zeta=\mathbf{1}_{(\rho_\sigma^-,\rho_Q)},
		\qquad x\ne\rho_Q,
		\qquad \zeta(\rho_Q)=\sigma,
		\qquad \zeta(\pm\infty)=0,
	\end{equation*}
	so $\overline u_\sigma<u^Q<\underline u_\sigma$ on $\mathbb{R}$.

	It remains to prove the mild inequalities.
	The only flux jump of $W_\beta$ occurs at $\rho_Q$ and equals
	$J_\beta:=h^s(\beta)-h^u(\beta)$.
	Thus, in distributions,
	\begin{equation*}
		-KW_\beta+\mathbf{1}_{\{x<\chi_\beta\}}=J_\beta\delta_{\rho_Q},
	\end{equation*}
	where $\delta_{\rho_Q}$ is the unit point mass at $\rho_Q$.
	Since $\rho_Q\in(L_1,L)$,
	Lemma~\ref{lem:local_branch_ordering} gives
	\begin{align}
		-K\underline u_\sigma+\mathbf{1}_{\{x<\rho_\sigma^+\}}
		 & =J_{a+\sigma}\delta_{\rho_Q}\geq0,
		\label{eq:U_sub}                      \\
		-K\overline u_\sigma+\mathbf{1}_{\{x<\rho_\sigma^-\}}
		 & =J_{a-\sigma}\delta_{\rho_Q}\leq0.
		\label{eq:U_sup}
	\end{align}

	Use $u_*=u^Q$ and $\rho_*=\rho_Q$ in \eqref{eq:prescribed_duhamel},
	and regard $\chi_\beta$ as a constant prescribed path.
	Write $w_\beta:=W_\beta-u^Q\in H^1(\mathbb{R})$ and
	$k_\beta:=\mathbf{1}_{\{x<\chi_\beta\}}-\mathbf{1}_{\{x<\rho_Q\}}
		\in L^2(\mathbb{R})$.
	Subtracting the stationary equation for $u^Q$ gives, in distributions,
	\begin{equation*}
		-Kw_\beta+k_\beta=J_\beta\delta_{\rho_Q}.
	\end{equation*}
	Point evaluation is continuous on $H^1(\mathbb{R})$:
	$|\psi(\rho_Q)|\leq C\|\psi\|_{H^1(\mathbb{R})}$.
	Since $w_\beta\in H^1(\mathbb{R})$ and $k_\beta\in L^2(\mathbb{R})$,
	density extends this identity to every $\psi\in H^1(\mathbb{R})$:
	\begin{equation*}
		-\int_{\mathbb{R}}
		\bigl(p\partial_xw_\beta\partial_x\psi+w_\beta\psi\bigr)\,dx
		+\int_{\mathbb{R}}k_\beta\psi\,dx
		=J_\beta\psi(\rho_Q).
	\end{equation*}

	To determine the sign of $\mathcal M[\chi_\beta,W_\beta](t)-W_\beta$,
	fix $t>0$ and a nonnegative $\varphi\in C_c^\infty(\mathbb{R})$,
	and define
	\begin{equation*}
		\psi_t:=\int_0^t \mathsf S(r)\varphi\,dr,
	\end{equation*}
	where the integral is taken in $L^2(\mathbb{R})$.
	Then $\psi_t\in\mathcal D(K)$ and $-K\psi_t=\mathsf S(t)\varphi-\varphi$.
	Indeed, for $s>0$ the semigroup property gives
	\begin{equation*}
		\frac{\mathsf S(s)\psi_t-\psi_t}{s}
		=\frac1{s}\int_t^{t+s}\mathsf S(r)\varphi\,dr
		-\frac1{s}\int_0^{s} \mathsf S(r)\varphi\,dr.
	\end{equation*}
	Strong continuity on $L^2(\mathbb{R})$ makes the right-hand side converge
	to $\mathsf S(t)\varphi-\varphi$ in $L^2(\mathbb{R})$ as $s\downarrow0$.
	The claim follows from the definition of the generator $-K$.
	In particular, $\psi_t\in H^1(\mathbb{R})$ and is an admissible test function.
	Positivity of $\mathsf S(r)$ gives $\psi_t\geq0$ almost everywhere,
	so its continuous representative satisfies $\psi_t(\rho_Q)\geq0$.

	Since $\chi_\beta$ is constant, changing the integration variable
	in \eqref{eq:prescribed_duhamel} gives
	\begin{equation*}
		\mathcal M[\chi_\beta,W_\beta](t)-W_\beta
		=\mathsf S(t)w_\beta-w_\beta+\int_0^t \mathsf S(r)k_\beta\,dr.
	\end{equation*}
	Testing against $\varphi$, using self-adjointness of $\mathsf S(r)$,
	and then using $-K\psi_t=\mathsf S(t)\varphi-\varphi$, we obtain
	\begin{align*}
		\int_{\mathbb{R}}
		\bigl(\mathcal M[\chi_\beta,W_\beta](t)-W_\beta\bigr)\varphi\,dx
		 & =-\int_{\mathbb{R}}w_\beta K\psi_t\,dx
		+\int_{\mathbb{R}}k_\beta\psi_t\,dx      \\
		 & =-\int_{\mathbb{R}}
		\bigl(p\partial_xw_\beta\partial_x\psi_t+w_\beta\psi_t\bigr)\,dx
		+\int_{\mathbb{R}}k_\beta\psi_t\,dx      \\
		 & =J_\beta\psi_t(\rho_Q).
	\end{align*}
	The second equality is weak integration by parts,
	valid because $w_\beta\in H^1(\mathbb{R})$ and $\psi_t\in\mathcal D(K)$.

	Since $\varphi$ is an arbitrary nonnegative test function,
	$J_{a+\sigma}>0$ and $J_{a-\sigma}<0$ imply
	\begin{equation*}
		\underline u_\sigma
		\leq\mathcal M[\rho_\sigma^+,\underline u_\sigma](t),
		\qquad
		\mathcal M[\rho_\sigma^-,\overline u_\sigma](t)
		\leq\overline u_\sigma,
		\qquad t\geq0.
	\end{equation*}
	These inequalities hold for the continuous spatial representatives.
	Together with strict decrease and the unique threshold crossings,
	they prove the asserted stationary mild sub- and supersolution properties.
\end{proof}

\section{Proof of the balanced minimum-gap criterion}
\label{app:balanced_minimum_gap}

\begin{proof}[Proof of Proposition~\ref{prop:balanced_minimum_gap}]
    Reflection symmetry fixes the middle crossing at $s_c$ and makes the exterior crossings symmetric about $s_c$.
    The strict inequality at the first crossing in the proof of Theorem~\ref{thm:all_front_existence} gives $b(\rho_1)<1/2$.
    Since $b$ is strictly decreasing on $(0,L_1)$ and $b(s_c)=1/2$, its first crossing phase must lie in $(s_c,L_1)$.
    Thus every front in the statement has the crossing locations asserted there, with $0<\delta<s_c$.

    Reflection of the coefficient about $s_c$ implies that the stable eigendirection there has slope $-r_c$.
    Hence the monodromy matrix based at $s_c$ is
    \begin{equation*}
        A_{s_c}=\begin{pmatrix}
            \cosh\ell&r_c^{-1}\sinh\ell\\
            r_c\sinh\ell&\cosh\ell
        \end{pmatrix},
        \qquad \ell:=\log\lambda.
    \end{equation*}
    Put $q:=\lambda^{-N}$, and recall $Y(\delta)=\cosh\delta+r_c\sinh\delta$ for $0\leq\delta\leq s_c$.
    At the first crossing, the bounded left tail has value $1/2$ and flux
    $-\tfrac12Y'(\delta)/Y(\delta)$.
    Evolve this state to $s_c$ using
    \begin{equation*}
        A_{s_c}^N
        \begin{pmatrix}\cosh\delta&-\sinh\delta\\
                       -\sinh\delta&\cosh\delta\end{pmatrix}.
    \end{equation*}
    Its value at $s_c$ equals $\tfrac12 R_q(\delta)$, where
    \begin{equation*}
        R_q(\delta)=\frac{q\cosh(2\delta)+D(q)\sinh(2\delta)}{Y(\delta)},
        \qquad
        D(q):=\frac{r_c-r_c^{-1}}{2q}+\frac{r_c+r_c^{-1}}{2}q.
    \end{equation*}
    Since $B(\delta)Y(\delta)=\cosh(2\delta)+\tfrac12(r_c+r_c^{-1})\sinh(2\delta)$,
    this agrees with \eqref{eq:balanced_scalar_matching}.
    Requiring $R_q(\delta)=1$ matches the value at the middle crossing.
    Reflection and complementation then supply the right half of the profile with continuous flux.
    This bounded function solves the linear equation with the prescribed alternating forcing and has value $1/2$ at all three interfaces.
    Lemma~\ref{lem:alternating_matching} therefore makes it a stationary front with exactly those crossings.
    The shooting calculation also proves necessity.

    We have $R_q(0)=q<1$ and $R_q'(\delta)>0$ for $0\leq\delta\leq s_c$.
    To verify the derivative sign, set $t:=\tanh\delta$ and suppress the positive common factor in its numerator.
    The remaining expression is
    \begin{equation*}
        q(3t+3r_ct^2-t^3-r_c)+D(q)(2+2r_ct^3).
    \end{equation*}
    Since $D(q)>qr_c/2$, this is strictly greater than
    $q\bigl(3t+3r_ct^2+(r_c^2-1)t^3\bigr)\geq0$.
    Finally, $R_q(s_c)=A(s_c)/q+B(s_c)q$.
    Strict increase proves the existence and uniqueness criterion, including its equality case.

    It remains to prove that the interior matching root persists when $N$ increases.
    Recall $R_q(\delta)=A(\delta)/q+B(\delta)q$ from \eqref{eq:balanced_scalar_matching}.
    At a matching root, the same matrix calculation gives the middle flux
    \begin{equation*}
        v(s_c)=\frac{r_c}{2}
              \left(\frac{A(\delta)}q-B(\delta)q\right)>0.
    \end{equation*}
    The strict sign follows from Lemma~\ref{lem:alternating_matching}, since this crossing is upward.
    Thus $q<\sqrt{A(\delta)/B(\delta)}$.
    For every $0<\widetilde q<q$, we obtain
    $R_{\widetilde q}(\delta)>R_q(\delta)=1$.
    As $R_{\widetilde q}(0)<1$ and $R_{\widetilde q}$ is strictly increasing, its root lies in $(0,\delta)$.
    In particular, increasing $N$ preserves the interior root and strictly decreases $\delta_N$.
    Since $A(s_c)>0$, \eqref{eq:balanced_gap_criterion} holds for all sufficiently large $N$.
    The roots tend to zero, because $R_q(\delta)\to\infty$ at every fixed $\delta>0$ as $q\downarrow0$.

    Finally, as $\alpha\uparrow1$ with $L_1,L_2$ fixed,
    the coefficient tends to one and the monodromy matrix tends to that of the homogeneous equation.
    Therefore $r_c\to1$, $\lambda\to e^L$, $A(s_c)\to0$, and $B(s_c)\to e^{s_c}$.
    For every fixed $N\geq1$, the left-hand side of \eqref{eq:balanced_gap_criterion} tends to $e^{s_c-NL}<1$.
    Every fixed finite set of gaps is consequently excluded when $\alpha$ is sufficiently close to one, proving $N_*\to\infty$.
\end{proof}

\section{Proof of Proposition~\ref{prop:pocket}}
\label{app:pocket}

Throughout, $a\in(a_-,a_+)$, $s=\rho_P$, $\Psi=\Psi^s$, and states are decomposed as in Section~\ref{subsec:global_manifold_illustration}.
By Lemma~\ref{lem:matching_function}, $b<a$ on $(\rho_P,\rho_Q)$ and $b>a$ on $(\rho_Q,\rho_P+L)$.

\begin{proof}[Proof of Proposition~\ref{prop:pocket}]
	\textit{(a)}
	Let $u$ be a front and $t$ its last threshold crossing.
	By Lemma~\ref{lem:equilibrium_tails}, its state is $aY_-(x)/y_-(t)$ for $x\geq t$.
	Choose the integer $j$ with $s+(j-1)L<t\leq s+jL$.
	Since $y_-$ is strictly decreasing and $y_-(x+L)=\lambda^{-1}y_-(x)$,
	the sample at $s+jL$ is $\mu\boldsymbol v_s$ with
	\begin{equation*}
		\mu=a\frac{y_-(s+jL)}{y_-(t)}\in(a/\lambda,a],
	\end{equation*}
	so it lies in the half-open segment $(P_1,P_0]$.
	Later samples have coefficients $\lambda^{-k}\mu\leq a/\lambda$ and lie outside $(P_1,P_0]$.
	An earlier sample in $(P_1,P_0]$ would be the state of a stable tail that stays below the threshold for all later $x$,
	contradicting the crossing at $t$.
	Hence each front has exactly one sample in $(P_1,P_0]$.

	Let $Z=\mu\boldsymbol v_s$ with $a/\lambda<\mu\leq a$.
	Following its stable tail backward, the solution first reaches the threshold at the phase $t\in(s-L,s]$ with
	$\mu=a\,y_-(s)/y_-(t)$, and $\mu$ increases with $t$.
	The points $P_1$, $Q_0$ and $P_0$ correspond to $t=s-L$, $t=\rho_Q-L$ and $t=s$.
	By the backward argument in the proof of Theorem~\ref{thm:all_front_existence},
	the solution escapes to $u\to+\infty$ as $x\to-\infty$ if $b(t)<a$,
	is the unstable tail if $b(t)=a$,
	and returns to the threshold if $b(t)>a$.
	For $t\in(s-L,\rho_Q-L)$ we have $b(t)<a$, so $Z$ is not on the unstable manifold of $(1,0)$.
	The equality $b(t)=a$ holds only at $Q_0$ and $P_0$, which are samples of $u^Q$ and $u^P$.
	For $t\in(\rho_Q-L,s)$ we have $b(t)>a$, and these are the points of the open base.
	Hence the unique sample of a front in $(P_1,P_0]$ lies in the base $[Q_0,P_0]$,
	and it lies in the open base exactly when the front is not monotone.
	Conversely, a common point of the unstable manifold of $(1,0)$ and the base
	is the state of a solution converging to $(1,0)$ as $x\to-\infty$ and,
	following the stable tail, to $(0,0)$ as $x\to+\infty$.
	By Proposition~\ref{prop:heteroclinic_representation} it is a front.
	$L$-translates have the same sample, and uniqueness of solutions makes the correspondence injective.

	\textit{(b)}
	The state $(1-c,-r_+(s)c)$ on $I$ is that of the unstable tail $E_1-cY_+(x)/y_+(s)$,
	which stays above the threshold until the phase $\sigma$ given by $c\,y_+(\sigma)/y_+(s)=1-a$.
	As $c$ decreases from $1-a$ to $(1-a)/\lambda$, the phase $\sigma$ increases from $s$ to $s+L$.
	The endpoint $P_0$ has $\sigma=s$, and $\Psi^{-1}(P_0)$ has $\sigma=s+L$.
	At $\sigma$ the state is $(a,-(1-a)r_+(\sigma))$, and, as in the proof of Theorem~\ref{thm:all_front_existence},
	as long as the solution stays below the threshold,
	\begin{equation*}
		Z(x)=\bigl(a-b(\sigma)\bigr)\frac{Y_+(x)}{y_+(\sigma)}+b(\sigma)\frac{Y_-(x)}{y_-(\sigma)}.
	\end{equation*}
	If $b(\sigma)>a$, both terms have negative flux, so the solution stays below the threshold and tends to $-\infty$.
	At $x=s+L$ its coefficients are
	\begin{equation*}
		\mu_+=\bigl(a-b(\sigma)\bigr)\frac{y_+(s+L)}{y_+(\sigma)}<0,
		\qquad
		\mu_-=b(\sigma)\frac{y_-(s+L)}{y_-(\sigma)}>0,
	\end{equation*}
	and $u<a$.
	This happens precisely for $\sigma\in(\rho_Q,s+L)$, whose images form the interior of $E$.
	If $b(\sigma)<a$, both coefficients are positive, $u$ stays positive, and the growing term forces a return to the threshold.
	This happens precisely for $\sigma\in(s,\rho_Q)$, whose images form the interior of $\Sigma_0$.
	The phase $\sigma=\rho_Q$ gives $u^Q$, whose image is $Q_0$.

	\textit{(c)}
	The segment $I$ is a simple arc and $\Psi$ is a homeomorphism, so $E$ is a simple arc.
	By (b), the interior of $E$ lies in $\{\mu_+<0\}$,
	and $E$ meets the stable line $\{\mu_+=0\}$ only at its endpoints.
	Hence $E$ and the base form a simple closed curve,
	which bounds a bounded open region, the pocket.
	This curve lies in the closed convex set $\{\mu_+\leq0,\ \mu_-\geq0,\ u\leq a\}$,
	so the pocket lies in its interior, where $\mu_+<0<\mu_-$ and $u<a$.
	Every such state has negative flux and escapes, by the decomposition in Remark~\ref{rem:spargel_departure},
	whose argument does not use $u>0$ when $\mu_+<0<\mu_-$.
	Hence the pocket contains no sample of a front.

	\textit{(d)}
	Let the last letter be $P$.
	By Theorem~\ref{thm:separated_fronts}, the last crossing is $\rho_m=N_mL+s+\delta$ with $|\delta|\leq C_m\lambda^{-n}$.
	The proof of Theorem~\ref{thm:all_front_existence} gives $b(\rho_m)>a$,
	and $b$ is decreasing near $s=\rho_P$, so $\delta<0$.
	The sample at $s+N_mL$ is $\mu\boldsymbol v_s$ with $\mu=a\,y_-(s)/y_-(s+\delta)$.
	If the last letter is $Q$, then $\rho_m=N_mL+\rho_Q+\delta$, and $b$ is increasing near $\rho_Q$, so $\delta>0$.
	The sample at $s+(N_m+1)L$ has coefficient $a\,y_-(s+L)/y_-(\rho_Q+\delta)$,
	compared with $q_0=a\,y_-(s+L)/y_-(\rho_Q)$ for $Q_0=q_0\boldsymbol v_s$.
	On each coefficient interval, $|y_-'/y_-|=p^{-1}r_-\leq\alpha^{-1/2}$ by Lemma~\ref{lem:floquet_slopes},
	and in both cases the ratio of the two values of $y_-$ involved is at most one.
	Hence the coefficients differ from $a$ and from $q_0$, respectively, by at most $\alpha^{-1/2}|\delta|$,
	and the samples differ from $P_0$ and $Q_0$ by at most $\alpha^{-1/2}|\boldsymbol v_s|\,|\delta|$,
	with $|\boldsymbol v_s|\leq(1+\alpha^{-1})^{1/2}$.
	Under the hypotheses of Corollary~\ref{cor:uniform_central_fronts}, let $n$ be the actual minimum slot gap.
	Since $n\geq N$, the range of thresholds allowed for $n$ contains the one for $N$.
	Applying the corollary with $N=n$ bounds $|\delta|$ by $C\lambda^{-n}$, with $C$ independent of $m$ and $a$.

	\textit{(e)}
	Fix an odd $m\geq3$ and an integer $N\geq n_m$ from Theorem~\ref{thm:separated_fronts}.
	For $j\geq(m-1)N-1$, take the slots $0,N,2N,\ldots,(m-2)N$ and $j+1$,
	whose gaps are at least $N$.
	Fix the first and the last letter to be $P$ and let the $m-2$ middle letters vary.
	This gives $2^{m-2}$ distinct fronts.
	The first crossing satisfies $b(\rho_1)<a$ by the proof of Theorem~\ref{thm:all_front_existence},
	so $\rho_1\in(s,\rho_Q)$, and the state of the front at $s$ lies on $I$ with a returning first crossing.
	By (b), its sample at $s+L$ lies on $\Sigma_0$, and its sample at $s+(j+1)L$ lies on $\Sigma_j$.
	By the proof of (d), this is its sample in the base, and it lies in the open base by (a).
	The fronts are distinct and all have their first occupied slot at zero,
	so none is a nontrivial period translate of another.
	By (a), their samples in the base are distinct.
	Take $j_m:=(m-1)N-1$.

	\textit{(f)}
	Use the coordinates $(\mu_+,\mu_-)$ at phase $s$.
	For $\sigma>s$ close to $s$, the state of $I$ with first crossing $\sigma$ is still below the threshold at $s+L$,
	and by the formula in (b) its sample there has
	\begin{equation*}
		\mu_+=\bigl(a-b(\sigma)\bigr)\frac{\lambda y_+(s)}{y_+(\sigma)},
		\qquad
		\mu_-=b(\sigma)\frac{y_-(s)}{\lambda y_-(\sigma)}.
	\end{equation*}
	At $\sigma=s$, the derivative of the first expression is $-\lambda b'(s)>0$ by Lemma~\ref{lem:matching_function}.
	Hence $\mu_+$ increases strictly from $0$ on a sufficiently short interval to the right of $s$.
	Hence the part of $\Sigma_0$ near $P_1$ is a graph $\mu_-=h(\mu_+)$, $0\leq\mu_+\leq\delta_0$,
	with $h$ continuous and $h(0)=a/\lambda$.
	Choose $\eta>0$ with $a/\lambda+\eta<a$, then $c\in(0,\delta_0]$ with $\lambda c+a/\lambda+\eta<a$ and $0<h\leq a/\lambda+\eta$ on $[0,c]$.
	A state with $0\leq\mu_+\leq c$ and $0<\mu_-\leq a/\lambda+\eta$ then stays below the threshold for one period,
	since its first component is at most $\lambda\mu_++\mu_-$ there,
	and $\Psi$ acts on it as $(\mu_+,\mu_-)\mapsto(\lambda\mu_+,\lambda^{-1}\mu_-)$.
	Applying $\Psi^k$ to the part of the graph with $0\leq\mu_+\leq c\lambda^{-k}$
	keeps every intermediate state in this region and shows that $\Sigma_k$ contains the graph
	\begin{equation*}
		\bigl\{(t,\lambda^{-k}h(\lambda^{-k}t)):0\leq t\leq c\bigr\}.
	\end{equation*}
	As $k\to\infty$, these graphs converge uniformly to the segment $\{(t,0):0\leq t\leq c\}$
	of the straight part of the unstable manifold of $(0,0)$.
	In these coordinates, $\mathcal C_r=\Psi^r(\{(t,0):0\leq t\leq a\})$.
	For $0\leq t\leq a/\lambda$, the first component $t\,y_+(x)/y_+(s)$ of the solution from $(t,0)$ increases to $\lambda t\leq a$ on $(s,s+L)$,
	so the solution stays below the threshold there and $\Psi(t,0)=(\lambda t,0)$.
	Hence $\mathcal C_r=\Psi^{r+i}(\{(t,0):0\leq t\leq\lambda^{-i}a\})$ for every $i\geq0$.
	Fix $i$ with $\lambda^{-i}a\leq c$, and put $n:=r+i$.
	Since $\Psi^n$ is continuous and injective on the compact segment $\{(t,0):0\leq t\leq\lambda^{-i}a\}$,
	the homeomorphism $\vartheta$ has the form $\vartheta(\tau)=\Psi^n(\phi(\tau),0)$ for a homeomorphism $\phi$ of $[0,1]$ onto $[0,\lambda^{-i}a]$.
	Put
	\begin{equation*}
		\vartheta_{k+n}(\tau):=\Psi^n\bigl(\phi(\tau),\lambda^{-k}h(\lambda^{-k}\phi(\tau))\bigr),
	\end{equation*}
	a homeomorphism of $[0,1]$ onto an arc of $\Sigma_{k+n}$.
	Since $0<\lambda^{-k}h\leq\lambda^{-k}(a/\lambda+\eta)$ and $\Psi^n$ is uniformly continuous on a compact neighbourhood of the segment,
	$\vartheta_{k+n}\to\vartheta$ uniformly on $[0,1]$ as $k\to\infty$.
	This proves (f).
\end{proof}

\printbibliography

\end{document}